%% file: FinRegModSt.tex
\documentclass[11pt]{amsart}

\usepackage[margin=1.05in]{geometry}

\usepackage[pdftex,pdfpagelabels,bookmarks,hyperindex,hyperfigures]{hyperref}
  \hypersetup{
    colorlinks,
    linkcolor={red!55!black},
    urlcolor={blue!55!black},
    citecolor={green!55!black}
  }

\input{./preamble}

\setenumerate{label=(\roman*)}

\usepackage[style=alphabetic]{biblatex}

\AtNextBibliography{\footnotesize}

\title{Finite regularity moduli stacks through diffeology}

\author{Mahmud Azam}
\email{mahmud.azam@usask.ca}

\author{Steven Rayan}
\email{rayan@math.usask.ca}

\address{Centre for Quantum Topology and Its Applications (quanTA) and
Department of Mathematics and Statistics, University of Saskatchewan, SK,
Canada~ S7N 5E6}

\begin{document}

\maketitle

\begin{abstract}
Diffeological non-Abelian Hodge theory motivates the study of $C^r$ families of 
Higgs bundles and flat bundles for finite $r$. In this article,
we provide a uniform language for treating moduli stacks parametrizing such objects
in terms of diffeology. More specifically, we develop a theory of stacks presentable by
diffeological groupoids over the site of $C^r$ manifolds for arbitrary $r$.
As an application, we construct moduli stacks over the
sites of $C^r$ manifolds parametrizing $C^r$ vector bundles over a smooth manifold, and $C^r$
families of Higgs bundles and of flat bundles over a complex manifold.
We also construct $C^r$ versions of the non-Abelian Hodge
filtration stack, interpolating between Higgs bundles and flat bundles through
$\lambda$--connections.
We show that all these stacks are presentable by diffeological groupoids.
Finally, we construct moduli stacks parametrizing $C^r$ families of
$\Delta^n$--shaped diagrams of
these objects, and show that they assemble into simplicial stacks whose levels are all
presentable by diffeological groupoids.
\end{abstract}

\tableofcontents

\section{Introduction}
\label{sec:introduction}

\subsection{Motivation}

Treating the complex non-Abelian Hodge correspondence at the level of moduli stacks requires
weakening some regularity hypotheses on families of Higgs bundles and flat bundles.
To see this, we recall an example from \cite[Proposition 14.4]{AR26b}:
consider an elliptic curve $X$, with a non-zero $1$--form
$\alpha \in H^0(X, \Omega^1_X)$ and a Higgs field on $\mcO_X^{\oplus 2}$ given by the following
matrix:
\[
\theta_z := \begin{bmatrix}
0 & z \\
z^2 & 0
\end{bmatrix} \alpha
\]
for $z \in \bC$. This defines are $C^\infty$ family of polystable Higgs bundles on $X$ parametrized
by $\bC$. For each $z \in \bC$, there is a harmonic metric defined by the matrix:
\[
h_z := \begin{bmatrix}
|z|^{1/2} & 0 \\
0 & |z|^{-1/2}
\end{bmatrix}
\]
The adjoint matrices of $\theta$ are:
\[
\begin{bmatrix}
0 & \bar{z}^2/|z| \\
|z|\bar{z} & 0
\end{bmatrix} \bar{\alpha} \text{ when $z \neq 0$ and $0$ otherwise}
\]
Hence, the flat connection family given by the non-Abelian Hodge correspondence is continuous
but not differentiable, even though the Higgs bundle family we started with is real analytic.
Thus, we are forced to consider moduli stacks parametrizing $C^r$ families of the respective objects
for finite $r$.
This need to work with $C^r$ families of objects required us to formulate moduli stacks over the
site of $C^r$ manifolds in \cite[\S 11]{AR26b}.
It is desirable to have precise foundations for this and to develop these foundations is the purpose
of this article.

\subsection{Overview and Main Results}

Most frameworks of moduli stacks in differential geometry, including differentiable stacks,
Fr\'echet differentiable stacks and diffeological stacks, are formulated
over the site $\Mfd$ of smooth manifolds with the Grothendieck topology generated by open covers.
However, in dealing with $C^r$ families, it is more natural to consider $C^r$ manifolds as the
parameter spaces, but to the best of
our knowledge, a treatment of geometric stacks over $C^r$ manifolds is difficult to find.
In this work, we develop a uniform language for this in terms of stacks presentable by diffeological
groupoids, in the sense of \cite{RV18}.
As an application, we produce several
moduli stacks, all presentable by diffeological groupoids, parametrizing objects of interest
in non-Abelian Hodge theory and complex geometry.
In \cref{sec:prest-st}, we prove some preliminary results on prestacks and stacks including
concrete descriptions of pushforward and pullback functors between $2$--categories
of prestacks $\PSt_\sC, \PSt_\sD$ over categories $\sC$ and $\sD$ induced by fully a faithful
functor $\sC \to \sD$. Using this concrete description, we prove a comparison
theorem (\cref{thm:comp-thm}) whose hypotheses are relatively simple to verify
for the inclusions of $C^r$ manifolds into diffeological spaces.

Using these results, in \cref{sec:dflg-st}, we show that stacks over the site $\Mfd^r$ of $C^r$
manifolds equipped with the Grothendieck topology generated by open covers are equivalent
to those over the site $\Dflg$ of diffeological spaces with the Grothendieck topology generated by
subductions, for all $r$.
This implies that stacks over $\Mfd^r$ are equivalent to those over
$\Mfd := \Mfd^\infty$ manifolds.
This is \cref{thm:St-Cr-Cinf-equiv} which can be read as: stacky $C^r$ geometry for finite $r$ is
equivalent to stacky $C^\infty$ geometry. Interestingly, this result requires no analysis.
In addition, we show that these equivalences of $2$--categories
of stacks restrict to stacks presentable by diffeological spaces in
\cref{prop:dflg-st-pullback-iff}.
This gives us a uniform language to treat families of objects with finite regularity in the
parameter direction, in terms of stacks over $\Mfd^r$ for various $r$, and then presenting
these stacks by diffeological groupoids.

In particular, for any smooth manifold $X$ and any $r = 0, 1, \dots, \infty$, using the fact that
the mapping stack $\HHom(X, BG)$ is presentable by a diffeological groupoid --- or simply,
diffeologically presentable --- for any diffeological group $G$,
we show that the moduli stack $\sM_{\Vect, r}(X)$ parametrizing $C^r$ vector bundles on $X$
is diffeologically presentable --- see \cref{thm:mod-st-Vec-Bun-dflg}.
For a complex manifold $X$, we then construct a moduli stack $\sM_{\Flt, r}(X)$ parametrizing $C^r$
families of bundles with generalized connections (in the sense of \cite[\S 2]{AR26a}) with respect
to differentials
obtained by complex linear combinations of holomorphic and antiholomorphic structures on vector
bundles --- see \cref{cns:M-Flt-r}.
This is a bidirectional version of the non-Abelian Hodge filtration stack from
\cite[\S 10]{NonAbHodgeFilt} in the context of differential geometry.
We show that this stack is diffeologically presentable in \cref{thm:M-Flt-r-dflg}.
This stack has a natural map to $\bC^2$, whose
pullback along the map $\bC \to \bC^2 : a \mapsto (1, a)$
yields the $C^r$ non-Abelian Hodge filtration stack $\sM_{\Hodge, r}(X)$.
Similarly, the fibres of $\sM_{\Flt, r}(X)$ above
$(0, 1)$ and $(1, 1)$ yield the moduli stacks
$\sM_{\Dol, r}(X)$ and $\sM_{\dR, r}(X)$ parametrizing $C^r$ families of Higgs bundles and flat
bundles on $X$ --- see \cref{defn:M-Dol-dR-Hodge}.
This immediately implies that these stacks are all diffeologically presentable --- see
\cref{thm:M-Dol-dR-Hodge-dflg}.
Finally, we construct simplicial stacks $S_\bullet\sM_{W, r}(X)$ for
$W = \Vect, \Flt, \Dol, \dR, \Hodge$ whose $n$--th level parametrizes $\Delta^n$--shaped diagrams
of the respective objects, where $\Delta^n$ is the standard $n$--simplex (the category associated to
the totally ordered set with $n + 1$ elements). We also show that these simplicial stacks
are diffeological at all levels --- see \cref{subsec:M-Diag-r}.

\subsection{Notation and Conventions}

We will use the concepts of abstract differentials and connections with respect to such
differentials developed in \cites{GR15, AR26a}.

\begin{conv}
All $C^r$ manifolds will be second countable and Hausdorff. In particular, they admit partitions
of unity.
\end{conv}

\begin{conv}
A $1$--category will refer to an ordinary category.
A $2$--category will refer to a bicategory with strictly associative and unital composition.
\end{conv}

\begin{notn}
We will use the following notation for various categories and $2$--categories:
\begin{figure}[H]
\begin{tabularx}{\textwidth}{l c l}
$\Set$ & : & sets and functions \\
$\Cat$ & : & categories, functors and natural transformations \\
$\Grpd$ & : & groupoids, functors and natural transformations \\
$\Top$ & : & topological spaces and continuous functions \\
$\Mfd^r$ & : & $C^r$ manifolds and $C^r$ functions \\
$\Cart^r$ & : & objects are open $U \subset \bR^n, n \in \bN \cup \set{0}$,
                and morphisms are $C^r$ functions \\
$\Cart$ & : & alternate notation for $\Cart^\infty$ \\
\end{tabularx}
\end{figure}
\end{notn}

\begin{notn}
We will use the following notation for various objects in (bi)category theory:
\begin{figure}[H]
\begin{tabularx}{\textwidth}{l c l}
$\Ob{\sC}$ & : & collection of objects of $\sC$ \\
$\Mor{\sC}$ & : & collection of morphisms of $\sC$ \\
$\sC(c, c')$ & : & collection of morphisms $c \to c'$ in $\sC$ \\
$\sC_{/f}$ & : & overcategory of a diagram $f : I \to \sC$ \\
$\sC_{f/}$ & : & undercategory of a diagram $f : I \to \sC$ \\
$\sC_{/c}$ & : & overcategory of the diagram $f : \set{*} \to \sC : * \mapsto c \in \Ob{sC}$ \\
$\sC_{c/}$ & : & undercategory of the diagram $f : \set{*} \to \sC : * \mapsto c \in \Ob{sC}$ \\
$F_{/f}$ & : & the functor $\sC_{/f} \to \sD_{/F \circ f}$
               induced by a (pseudo)functor $F : \sC \to \sD$ \\
$F_{f/}$ & : & the functor $\sC_{f/} \to \sD_{F \circ f/}$
               induced by a (pseudo)functor $F : \sC \to \sD$ \\
$F_{/c}$ & : & the functor $\sC_{/c} \to \sD_{/F(c)}$
               induced by a (pseudo)functor $F : \sC \to \sD$ \\
$F_{c/}$ & : & the functor $\sC_{c/} \to \sD_{F(c)/}$
               induced by a (pseudo)functor $F : \sC \to \sD$ \\
$\Fun(\sC, \sD)$ & : & (bi)category of (pseudo)functors $\sC \to \sD$ \\
$\sC^\op$ & : & opposite (bi)category of $\sC$ \\
$\sP_\sD(\sC)$ & : & $\Fun(\sC^\op, \sD)$ \\
$\sP(\sC)$ & : & $\sP_\Set(\sC)$ \\
$\Sh_\sD(\sC, J)$ & : & (2--)sheaves on the site $(\sC, J)$ valued in $\sD$ \\
$\Sh(\sC, J)$ & : & $\Sh_\Set(\sC, J)$ \\
$\PSt_{\sC}$ & : & prestacks or Grothendieck fibrations over a $1$--category $\sC$ \\
$\St_{\sC, T}$ & : & stacks over a site $(\sC, J)$ \\
$\int F \to \sC$ & : & Grothendieck construction of a pseudofunctor $F : \sC^\op \to \Cat$ \\
\end{tabularx}
\end{figure}
\end{notn}

\begin{notn}
We will suprress $J$ from $\Sh_\sD(\sC, J)$ and $\St_{\sC, J}$, when $J$ is clear from context.
We will do this, in particular, for the Grothendieck topologies fixed in
\cref{conv:standard-topologies}.
\end{notn}

\noindent\textbf{Acknowledgements.}
The first-named author was supported by
an NSERC Canada Graduate Scholarship (Doctoral),
a Canadian Department of National Defence/NSERC Supplemental Funding Award and
a University of Saskatchewan Teacher-Scholar Doctoral Fellowship during this work.
The second-named author acknowledges partial support from a Natural Sciences and Engineering
Research Council of Canada (NSERC) Discovery Grant. He thanks Syddansk Universitet (SDU),
the Centre for Quantum Mathematics at SDU, the Danish e-infrastructure Consortium (DeiC),
and the Danish Institute for Advanced Study for hosting him in August 2026 during the Niels
Bohr Quantum Summer School and for providing an inspiring and comfortable research environment.

\section{Some Results on Stacks}
\label{sec:prest-st}

The main purpose of this section is to prove a simple comparison theorem (\cref{thm:comp-thm}) which
enables us to treat stacks over $\Mfd^r$ as stacks over $\Dflg$ and vice versa.
In \cref{subsec:prest-st} we review prestacks and stacks over a site. A slight addition to the
standard theory is a seemingly stronger a definition of a stack (\cref{defn:st}) which we find to be
slightly cleaner when reasoning about descent for morphisms, but which we show to be equivalent to
the standard definition --- see \cref{prop:morphism-coll-fibrewise} and
\cref{rmk:morphism-coll-fibrewise}. In \cref{subsec:pullback-pushfwd-prest},
we give some concrete characterizations of pullbacks and pushforwards of prestacks
along fully faithful functors. These characterizations help in proving the comparison theorem.
Here, we also show that $2$--fibre products of spans where morphism is a Cartesian fibration
can be computed as strict fibre products, which gives concrete but simple descriptions of many of
the stacks of interest in \cref{sec:dflg-st} and simplifies several proofs.
In \cref{subsec:pullback-pushfwd-st}, we review continuous and cocontinuous functors and show that
pullbacks and pushforwards, respectively, along them preserve stacks.
In \cref{subsec:comp-thm}, we prove our comparison theorem, the main result of this section.
In \cref{subsec:repr-obj-morphisms}, we review the concepts of representable objects and morphisms
which we use to reason about diffeologically presentable stacks in \cref{sec:dflg-st}.

\subsection{Prestacks and Stacks}
\label{subsec:prest-st}

\begin{defn}[Cartesian Fibration]
Let $p : \sP \to \sC$ be a functor. We call a morphism $f : P \to P'$ in $\sP$ a Cartesian morphism
if for every morphism $g : P'' \to P'$ in $\sP$ and every factorization in $\sC$ of the form:
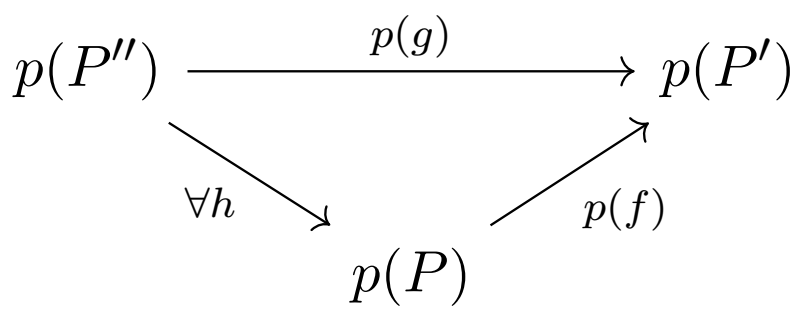
\[\begin{tikzcd}
p(P'') \ar[rd, "\forall h" below left] \ar[rr, "p(g)"] & & p(P') \\
& p(P) \ar[ru, "p(f)" below right] &
\end{tikzcd}\]
there exists a unique factorization in $\sP$ of the form:
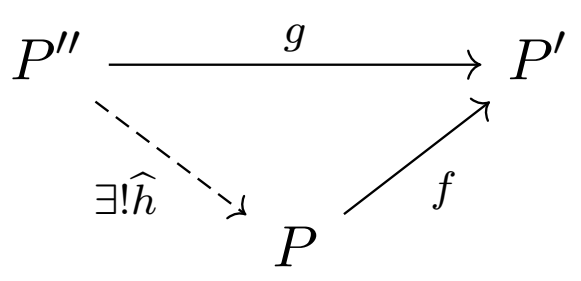
\[\begin{tikzcd}
P'' \ar[rd, dashed, "\exists ! \wh{h}" below left] \ar[rr, "g"] & & P' \\
& P \ar[ru, "f" below right] &
\end{tikzcd}\]
such that $p(\wh{h}) = h$. In this case, we call $\wh{h}$ a Cartesian lift of $h$.
The functor $p$ is called a Cartesian fibration if for every morphism $f : C \to C'$
in $\sC$ and every object $P'$ in $\sP$ with $p(P') = C'$, $f$ admits a Cartesian lift
$\wh{f} : \wh{C} \to P'$.
\end{defn}

\begin{rmk}
A Cartesian fibration is also called a fibred category in the literature
\cite[\href{https://stacks.math.columbia.edu/tag/02XJ}{Section 02XJ}]{stacks-project}.
\end{rmk}

\begin{defn}[Fibre Category]
Let $p : \sP \to \sC$ be a functor and $C \in \Ob{\sC}$. Consider the disrete subcategory
$\set{C} \subset \sC$ with the single object $C$. Then, we will call the $2$--fibre product
$\sP \times_\sC \set{C}$ in $\Cat$ the fibre category of $p$ at $\sC$ and denote it as $\sP_C$.
\end{defn}

\begin{defn}[Prestack]
A prestack on a category $\sC$ is a Cartesian fibration $p : \sP \to \sC$ whose fibre categories
are groupoids.
We will abuse notation and identify the prestack $p$ with the category $\sP = \dom(p)$.
We call $p$ the structure morphism or the projection of the prestack.

A morphism of prestacks $\sP \to \sQ$ is a functor $\sP \to \sQ$ which strictly
commutes with the projections.
We denote by $\PSt_{\sC}$ the $2$--category whose objects are prestacks over $\sC$,
whose $1$--morphisms are morphisms of prestacks, and $2$--morphisms
$\alpha : f \To g : \sP \to \sQ$ are natural transformations of functors $\alpha : f \To g$
such that the component morphisms of $\alpha$ in $\sQ$ project to identity morphisms in $\sC$.
\end{defn}

\begin{rmk}
A prestack is also called a category fibred in groupoids in the literature
\cite[\href{https://stacks.math.columbia.edu/tag/003S}{Section 003S}]{stacks-project}.
\end{rmk}

\begin{defn}[Grothendieck Topology]
A Grothendieck topology $J$ on a category $\sC$ is a choice, for each object $C \in \Ob{\sC}$,
of a set $J(C)$ of families of morphisms $\set{c_i : C_i \to C} \in J(C)$ satisfying the following
conditions:
\begin{enumerate}
\item \textbf{Isomorphisms:}
Every isomomorphism $C' \to[\cong] C$ is in $J(C)$.
\item \textbf{Compositionality:}
For each $\set{C_i \to C}_{i \in I} \in J(C)$ and each $\set{C^j_i \to C_i}_{j \in L_i} \in J(C_i)$,
the family $\set{C^j_i \to C_i \to C}_{i \in I, j \in L_i}$ is in $J(C)$.
\item \textbf{Pullback Stability:}
For each $\set{C_i \to C}_{i \in I} \in J(C)$ and each morphism $f : C' \to C$, the fibre product
$C' \times_C C_i \to C'$ exists and the family $\set{C' \times_C C_i \to C'}_{i \in I}$
is in $J(C')$.
\end{enumerate}
We will call a member $\set{C_i \to C}_{i \in I}$ of $J(C)$ for any object
$C \in \Ob{\sC}$, a covering family, when the Grothendieck topology $J$ has been fixed or is clear
from context.
\end{defn}

\begin{defn}[Site]
A site is a pair $(\sC, J)$ where $\sC$ is a category with a choice of Grothendieck topology $J$.
\end{defn}

\begin{notn}
Given a site $(\sC, J)$, a covering family $\set{C_i \to C}_{i \in I}$ any finite subset
$I' = \set{i_1, \dots, i_k} \subset I$, we will write $C_{I'}$ for the finite iterated
fibre product $C_{i_1} \times_C \cdots \times_C C_{i_k}$. For $k = 2$ or $3$, we will simply write
$C_{i_1i_2}$ or $C_{i_1i_2i_3}$ respectively.
\end{notn}

\begin{defn}[Restrictions of Objects]
Consider a site $(\sC, J)$, a prestack $p : \sP \to \sC$, a morphism $f : C \to C'$
in $\sC$, an object $P'$ in $\sP$ with $p(P') = C'$ and a Cartesian lift $\wh{f} : \wh{C} \to P'$
of $f$. If $f$ and $\wh{f}$ are clear from context, then we will often write $\wh{C}$
as $P'|C$, and call it the restriction of $P'$ to $C$.
\end{defn}

\begin{defn}[Restrictions of Morphisms]
Consider a site $(\sC, J)$, a prestack $p : \sP \to \sC$, and a commutative square in $\sC$:
\[\begin{tikzcd}
A \ar[r, "p"] \ar[d, "r" left] & B \ar[d, "q"] \\
C \ar[r, "s" below] & D
\end{tikzcd}\]
For any morphism $\wh{q} : \wh{B} \to \wh{D}$ with $p(\wh{q}) = q$, and choices of Cartesian lifts
$\wh{p} : \wh{A} \to \wh{B}$ and $\wh{s} : \wh{C} \to \wh{D}$ of $p$ and $s$ respectively, by the
Cartesian property of $\wh{s}$, there exists a unique morphism $\wh{r} : \wh{A} \to \wh{C}$ making
the following diagram commute:
\[\begin{tikzcd}
\wh{A} \ar[r, "\wh{p}"] \ar[d, "\wh{r}" left] & \wh{B} \ar[d, "\wh{q}"] \\
\wh{C} \ar[r, "\wh{s}" below] & \wh{D}
\end{tikzcd}\]
In the case $r = \id_A = \id_C$ and $q = \id_B = \id_D$, if $p = s$ is clear from context,
we will often write $\wh{q}|A = \wh{q}|C$ for $\wh{r}$, and call it the restriction of
$\wh{q}$ to $A = C$.
\end{defn}

\begin{defn}[Object Descent Datum]
Let $(\sC, J)$ be a site and $p : \sP \to \sC$, a prestack. An object descent datum
in $\sP$ consists of:
\begin{enumerate}
\item an object $C$ in $\sC$,
\item a covering family $\set{c_i : C_i \to C}_{i \in I}$ in $\sC$ with respect to $J$,
\item a family of objects $\set{P_i}_{i \in I}$ of $\sP$ such that $p(P_i) = C_i$ for each
      $i \in I$, and
\item a family of morphisms $\set{\alpha_{ij} : P_i|C_{ij} \to P_j|C_{ij}}_{i, j \in I}$ in $\sP$
      with $p(\alpha_{ij}) = \id_{C_{ij}}$ for each $i, j \in I$,
\end{enumerate}
such that the following diagram commutes for each $i, j \in I$:
\[\begin{tikzcd}
P_i|C_{ijk} \ar[rr, "\alpha_{ij}|C_{ijk}" above] \ar[rd, "\alpha_{ik}|C_{ijk}" below left] & &
P_j|C_{ijk} \ar[ld, "\alpha_{jk}|C_{ijk}" below right] \\
& P_{k}|C_{ijk} &
\end{tikzcd}\]
A collation of the object descent datum consists of:
\begin{enumerate}
\item an object $P$ in $\sP$ with $p(P) = C$, and
\item morphisms $\phi_i : P_i \to P$ with $p(\phi_i)$,
\end{enumerate}
such that the following diagram commutes for each $i, j \in I$:
\[\begin{tikzcd}
P_i|C_{ij} \ar[rr, "\alpha_{ij}"] \ar[rd, "\phi_i|C_{ij}" below left] & &
P_j|C_{ij} \ar[ld, "\phi_j|C_{ij}" below right] \\
& P|C_{ij} &
\end{tikzcd}\]
The problem of finding a collation of an object descent datum is called an object descent problem.
\end{defn}

\begin{defn}[Morphism Descent Datum]
Let $(\sC, J)$ be a site and $p : \sP \to \sC$, a prestack. A morphism descent datum
in $\sP$ consists of:
\begin{enumerate}
\item a morphism $f : C \to C'$ in $\sC$,
\item an object $P$ in $\sP$ with $p(P) = C$
\item an object $P'$ in $\sP$ with $p(P') = C'$
\item a covering family $\set{c_i : C_i \to C}_{i \in I}$ in $\sC$ with respect to $J$, and
\item a family of morphisms $\set{\phi_i : P|C_i \to P'}_{i \in I}$ of $\sP$ such that
$p(\phi_i)$ is the composite $C_i \to[c_i] C \to[f] C'$ for each $i \in I$,
\end{enumerate}
such that the following diagram commutes for each $i, j \in I$:
\[\begin{tikzcd}
P|C_{ij} \ar[r] \ar[d] & P|C_i \ar[d, "\phi_i"] \\
P|C_j \ar[r, "\phi_j" below] & P'
\end{tikzcd}\]
A collation of a morphism descent datum consists of a morphism $\phi : P \to P'$ such
that the following diagram commutes for each $i \in I$:
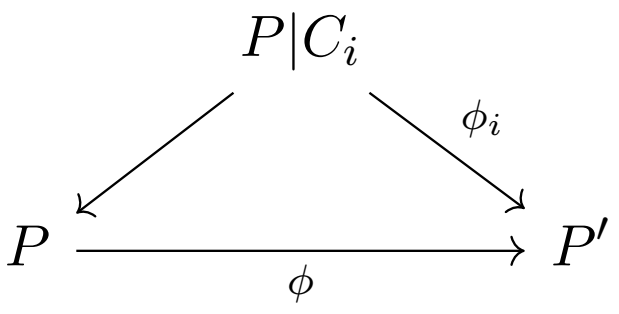
\[\begin{tikzcd}
& P|C_i \ar[ld] \ar[rd, "\phi_i"] & \\
P \ar[rr, "\phi" below] & & P'
\end{tikzcd}\]
The problem of finding a collation of a morphism descent datum is called a morphism descent problem.
\end{defn}

\begin{defn}[Stack]\label{defn:st}
A stack on a site $(\sC, J)$ is a prestack $p : \sP \to \sC$ on the underlying category such
that every object descent datum has a not-necessarily-unique collation, and
every morphism descent datum has a unique collation.
\end{defn}

\begin{prop}\label{prop:morphism-coll-fibrewise}
Suppose a prestack $p : \sP \to \sC$ over a site $(\sC, J)$ has a unique collation for every
morphism descent datum where the underlying morphism in $\sC$ is $\id_C$ for some object $C$
in $\sC$. Then, it has a unique collation for every morphism descent datum.
\end{prop}
\begin{proof}
Consider a morphism descent datum, with notation as in the definition.
Then, consider $P'|C$, where the restriction is along $f$. The Cartesian property of the morphism
$P'|C \to P'$ gives unique morphisms $\psi_i : P|C_i \to P'|C$ such that the composite
$P|C_i \to[\psi_i] P'|C \to P'$ is $\phi_i$. Then, consider the following diagram for each
$i, j \in I$:
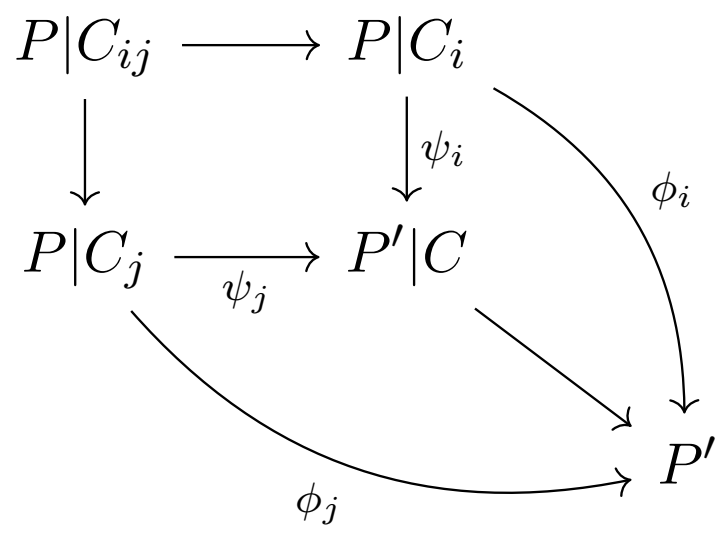
\[\begin{tikzcd}
P|C_{ij} \ar[r] \ar[d] & P|C_i \ar[d, "\psi_i"] \ar[rdd, bend left, "\phi_i" above right] & \\
P|C_j \ar[r, "\psi_j" below] \ar[rrd, bend right, "\phi_j" below left] &
P'|C \ar[rd] & \\
& & P'
\end{tikzcd}\]
Since all paths from $P|C_{ij}$ ot $P'$ in this diagram are equal, the Cartesian property
of the morphism $P'|C \to P'$ guarantees that top-left square commutes.
Thus, the $\psi_i$ form a morphism descent datum whose underlying morphism in $\sC$ is $\id_C$.
By hypothesis, we have a collation $P \to P'|C$ for this morphism descent datum which, composed
with the morphism $P'|C \to P'$, provides a collation of the original morphism descent datum.
\end{proof}

\begin{rmk}\label{rmk:morphism-coll-fibrewise}
\cref{prop:morphism-coll-fibrewise} shows that the above notion of stack is equivalent
to the standard one as in, for instance, \cite[\S 2.5]{AlperModuli}.
\end{rmk}

\begin{conv}\label{conv:standard-topologies}
We will consider the categories $\Top, \Cart^r, \Mfd^r$ to be equipped with the Grothendieck
topologies generated by families of jointly surjective open embeddings: the standard open
cover topology.
\end{conv}

\subsection{Pullback and Pushforward of Prestacks}
\label{subsec:pullback-pushfwd-prest}

\begin{defn}[Concrete $2$--Fibre Products of Categories]\label{defn:conc-fibre-prod}
Let $\sA \to[f] \sB \ot[g] \sC$ be functors. We define the concrete
$2$--fibre product $\sA \times_\sC \sB$ to be the category with the following description:
\begin{itemize}
\item Objects are tuples $(A, C, \alpha)$, where:
    \begin{itemize}
    \item $A \in \Ob{\sA}$
    \item $C \in \Ob{\sC}$
    \item $\alpha : f(A) \to g(C)$ is an isomorphism in $\sB$
    \end{itemize}
\item Morphisms $(A, C, \alpha) \to (A', C', \alpha')$ are pairs $(a, b)$, where
    $a : A \to A', b : C \to C'$ are morphisms in $\sB$ making the following
    diagram commute:
    \[\begin{tikzcd}
    f(A) \ar[r, "f(a)"] \ar[d, "\alpha" left] & f(A') \ar[d, "\alpha'"] \\
    g(C) \ar[r, "f(b)" below] & g(C')
    \end{tikzcd}\]
\end{itemize}
Composition of morphisms is defined entry-wise and identity morphisms are the pairs whose entries
are both identities.
\end{defn}

\begin{rmk}
The concrete $2$--fibre product of categories is indeed a $2$--fibre product
\cite[\href{https://stacks.math.columbia.edu/tag/0040}{Lemma 0040}]{stacks-project}.
\end{rmk}

\begin{defn}[Strict Fibre Product of Categories]
Let $\sA \to[f] \sB \ot[g] \sC$ be functors. We define the strict fibre product
$\sA \times_\sB^{str} \sC$ to be the subcategory of $\sA \times \sC$ consisting of
those objects $(A, C)$ such that $f(A) = g(C)$ and those morphisms $(a, c)$ such that
$f(a) = g(c)$.
\end{defn}

\begin{prop}\label{prop:strictify-2-fibre-prod}
For functors $\sA \to[f] \sB \ot[g] \sC$, if $g$ is a Cartesian fibration,
then $\sA \times_\sB \sC$ is equivalent as a category to the strict fibre product
$\sA \times^{str}_\sB \sC$.
\end{prop}
\begin{proof}
We have a functor $\sA \times^{str}_\sB \sC \to \sA \times_\sB \sC$ defined, on objects, by
$(A, C) \mapsto (A, C, \id_{f(A)}) = (A, C, \id_{g(C)})$ since $f(A) = g(C)$, and on morphisms,
by $(a, b) \mapsto (a, b)$. Notice that $f(a) = g(b)$ implies that the condition on morphisms
given in \cref{defn:conc-fibre-prod} is automatic. Also, notice that this functor is fully faithful
as $(a, b)$ is a morphism $(A, C, \id_{f(A)}) \to (A', C', \id_{f(A)})$ if and only if
$f(a) = g(b)$.

It remains to see essential surjectivity. Suppose $(A, C, \alpha)$ is an
object of $\sA \times_\sB \sC$. Then, choose a Cartesian lift $\wh{\alpha} : \wh{f(A)} \to C$.
This provides an object $(A, \wh{f(A)}, \id_{f(A)})$ in $\sA \times_\sB \sC$ as
$g(\wh{\alpha}) = \alpha$ also implies $g(\wh{f(A)}) = f(A)$.
Then, observe that the following diagram commutes:
\[\begin{tikzcd}[column sep=large]
f(A) \ar[r, "f(\id_A) = \id_{f(A)}"] \ar[d, "\id_{f(A)}" left] & f(A) \ar[d, "\alpha"] \\
g(\wh{f(A)}) = f(A) \ar[r, "g(\wh{\alpha}) = \alpha" below] & g(C)
\end{tikzcd}\]
Hence, we have a morphism  $(\id_A, \wh{\alpha}) : (A, \wh{f(A)}, \id_{f(A)}) \to (A, C, \alpha)$ in
$\sA \times_\sB \sC$.
Since the Cartesian lift of an isomorphism is an isomorphism
\cite[\href{https://stacks.math.columbia.edu/tag/02XL}{Lemma 02XL}]{stacks-project},
$\wh{\alpha}$ is an isomorphism
and, of course, $\id_{A}$ is an isomorphism. Thus, $(\id_A, \wh{\alpha})$ is an isomorphism,
while $(A, \wh{f(A)}, \id_{f(A)})$ is the image of $(A, \wh{f(A)})$ in $\sA \times_\sB^{str} \sC$.
Therefore, the mapping $\sA \times^{str}_\sB \sC \to \sA \times_\sB \sC$ is an equivalence.
\end{proof}

\begin{cor}\label{cor:prest-pullback-strict}
Given a prestack $p : \sP \to \sD$ and a functor $u : \sC \to \sD$, the strict fibre product
$\sC \times_\sD^{str} \sP$ is a $2$--fibre product.
\end{cor}
\begin{proof}
Apply \cref{prop:strictify-2-fibre-prod}.
\end{proof}

\begin{prop}\label{prop:prest-pb}
For a prestack $p : \sP \to \sD$ and a functor $u : \sC \to \sD$, the canonical projection of the
$2$--fibre product of categories $u^*p : u^*\sP = \sC \times_\sD \sP \to \sC$ is a prestack.
In particular, we get a strict $2$--functor of $2$--categories:
\[
u^* : \PSt_{\sD} \to \PSt_{\sC}
\]
sending prestacks over $\sD$ to their strict pullbacks along $\sC$.
\end{prop}
\begin{proof}
Consider a morphism $f : C_1 \to C_2$ in $\sC$ and choose a Cartesian lift
$\wh{u(f)} : P_1 \to P_2$ of $u(f)$ in $\sP$.
Then, by \cref{cor:prest-pullback-strict} we have a morphism
$(f, \wh{u(f)}) : (C_1, P_1) \to (C_2, P_2)$ in $\sP \times_\sD \sC$ that maps to
$f$ under the projection $u^*p$. Since composition in $u^*p$ is entrywise, it follows
that $(f, \wh{u(f)})$ is Cartesian since $\wh{u(f)}$ is Cartesian. It remains to see that
the fibre categories are groupoids. For this, consider an object $C$ in $\sC$ and the morphism
$\iota_C : \set{*} \to \sC$ sending $*$ to $C$. Then, we have the following pasting
of $2$--pullback diagrams:
\[\begin{tikzcd}
(\sC \times_\sD \sP)_C \ar[r] \ar[d, "\iota_C^*u^*p" left]
    \ar[rd, phantom, "\lrcorner" very near start] &
\sC \times_\sD \sP \ar[r] \ar[d, "u^*p" left] \ar[rd, phantom, "\lrcorner" very near start] &
\sP \ar[d, "p"] \\
\set{*} \ar[r, "\iota_C" below] &
\sC \ar[r, "u" below] & \sD
\end{tikzcd}\]
showing that the fibre category $(\sC \times_\sD \sP)_C$ is equivalent (in fact, isomorphic,
by \cref{cor:prest-pullback-strict}) to the fibre category $\sP_{u(D)}$ which is a groupoid since
$\sP$ is a prestack. By \cref{cor:prest-pullback-strict} again, the mapping $\sP \mapsto u^*\sP$
simply takes the strict pullback along $u$ and is hence a strict functor.
\end{proof}

\begin{cns}[Relative Overcategory]\label{cns:rel-over-cat}
Let $u : \sC \to \sD$ be a functor.
For each $D \in \Ob{\sD}$, we write $u(\sC)_{/D}$ for the category whose:
\begin{itemize}
\item Objects are pairs $(C, f)$ where:
    \begin{itemize}
    \item $C \in \Ob{\sC}$.
    \item $f : u(C) \to D$ is a morphism in $\sD$.
    \end{itemize}
\item Morphisms $(C, f) \to (C', f')$ are morphisms $g : C \to C'$ in $\sC$
such that the following diagram commutes in $\sD$:
\[\begin{tikzcd}
u(C) \ar[rd, "f" below left] \ar[rr, "u(g)"] & & u(C') \ar[ld, "f'" below right] \\
& D &
\end{tikzcd}\]
\end{itemize}
We have a functor $u(\sC)_{/D} \to \sC$ given by $(C, f) \mapsto C$ and $g \mapsto g$.
Given a morphism $h : D \to D'$, we have a functor $u(\sC)_{/h} : u(\sC)_{/D} \to u(\sC)_{/D'}$
defined by $(C, f) \mapsto (C, h \circ f)$ for objects and $g \mapsto g$ for morphisms.
\end{cns}

\begin{cns}[Pushforward of a Prestack]\label{cns:prest-pushfwd}
Let $p : \sP \to \sC$ be a prestack and $u : \sC \to \sD$ be a functor.
We define a category $u_*\sP$ as follows:
\begin{itemize}
\item Objects are pairs $(D, F)$ where:
    \begin{itemize}
    \item $D \in \Ob{\sD}$.
    \item $F$ is a functor $u(\sC)_{/D} \to \sP$ such that following diagram strictly commutes:
    \[\begin{tikzcd}
    u(\sC)_{/D} \ar[rr, "F"] \ar[rd] & & \sP \ar[ld, "p" below right] \\
    & \sC &
    \end{tikzcd}\]
    \end{itemize}
\item Morphisms $(D, F) \to (D', F')$ are pairs $(f, \alpha)$ where:
    \begin{itemize}
    \item $f : D \to D'$ is a morphism in $\sC$.
    \item $\alpha$ is a natural isomorphism $F' \circ u(\sC)_{/f} \To F$ whose components
    lie over the identity morphisms in $\sC$.
    \end{itemize}
\end{itemize}
Given two morphisms $(D, F) \to[(f, \alpha)] (D', F') \to[(g, \beta)] (D'', F'')$, the composite
is defined in the first entry by $g \circ f$ and in the second, by the following pasting of
$2$--cells:
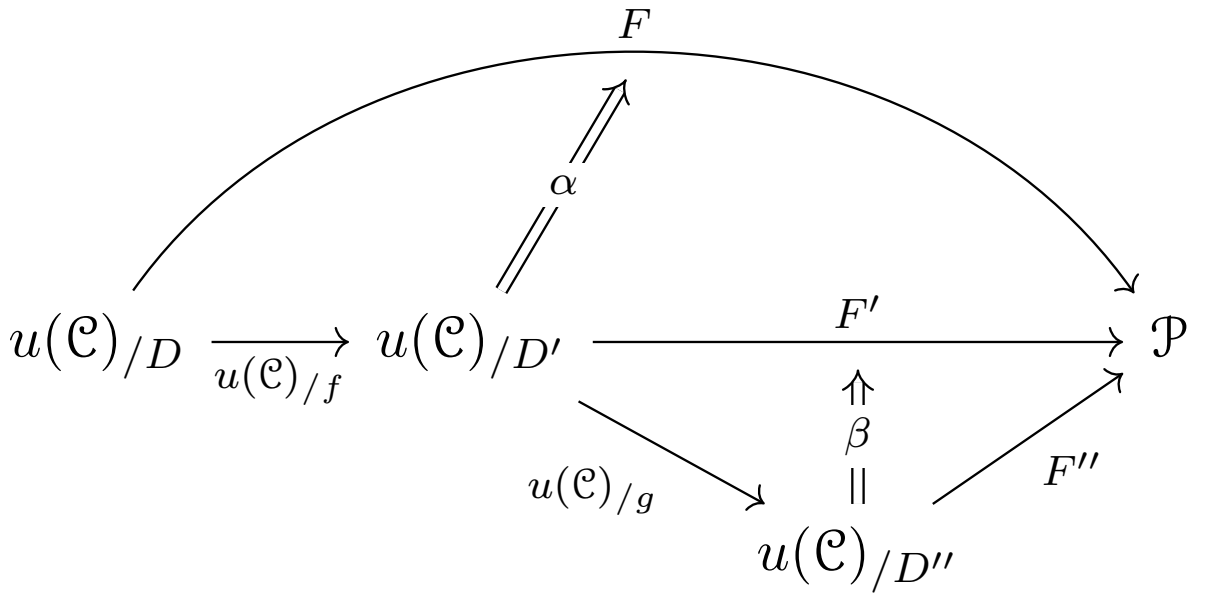
\[\begin{tikzcd}
u(\sC)_{/D}
    \ar[r, "u(\sC)_{/f}" below]
    \ar[rrr, bend left=5em, "F" above, ""{name=T1, below}] &
u(\sC)_{/D'}
    \ar[rd, "u(\sC)_{/g}" below left]
    \ar[rr, "F'" above, ""{name=T2, below}] & &
\sP \\
& & u(\sC)_{/D''} \ar[ru, "F''" below right]
\ar[from=1-2, to=T1, Rightarrow, "\alpha" description]
\ar[from=2-3, to=T2, Rightarrow, "\beta" description]
&
\end{tikzcd}\]
Given an object $(D, F)$, the identity morphism is $(\id_D, \id_F)$.
There is a functor $u_*p : u_*\sP \to \sD$ defined by $(D, F) \mapsto D$ and
$(f, \alpha) \mapsto f$.
\end{cns}

\begin{prop}
For a prestack $p : \sP \to \sC$ and a functor $u : \sC \to \sD$, $u_*p : u_*\sP \to \sD$
is a prestack.
\end{prop}
\begin{proof}
Let $f : D \to D'$ be a morphism in $\sD$ and $(D', F)$ be an object of $u_*\sP$.
Then, $(f, \id_{F \circ u(\sC)_{/f}})$ is a Cartesian morphism over $f$.
Let $D$ be an object of $\sD$ and $(\id_D, \alpha)$ be a morphism over $\id_D$ in $u_*\sP$.
Then, by definition, $\alpha$ is an natural isomorphism of functors. It follows that
$(\id_D, \alpha)$ is an isomorphism in $u_*\sP$.
\end{proof}

\begin{cns}[Unit Morphism]\label{cns:unit-morphism}
In the context of \cref{cns:prest-pushfwd}, suppose $u$ is fully faithful.
We construct a morphism of prestacks $\eta_\sP : u^*u_*\sP \to \sP$ as
follows. Let $(D, F)$ be an object of $u^*u_*\sP$ so that $D = u(C)$ for some $C \in \Ob{\sC}$.
We have $\id_{u(C)} \in u(\sC)_{/D} = u(\sC)_{/u(C)}$. We set:
\[
\eta_\sP(D, F) := F(C, \id_{u(C)})
\]
Then, let $(f, \alpha) : (D, F) \to (D', F')$ be a morphism in $u^*u_*\sP$ so that
$u_*p(f, \alpha) = u(f_0) \implies f= u(f_0)$ for some morphism $f_0 : C \to C'$ in $\sC$ with
$u(C) = D, u(C') = D'$.
We also have the isomorphism
$\alpha_{C, \id_{u(C)}} : F'(C, f \circ \id_{u(C)}) = F'(C, f) \to F(C, \id_{u(C)})$
in $\sP$. $f$ also gives a morphism in $u(\sC)_{/u(C')}$ of the form
$f \to \id_{u(C')}$. Hence, we get a morphism $F'(f_0) : F'(C, f) \to F'(C', \id_{u(C')})$ in $\sP$.
We define $\eta_\sP(f, \alpha)$ to be the composite:
\[
\eta_\sP(D, F) = F(C, \id_{u(C)})
\to[\alpha_{C, \id_{u(C)}}^{-1}] F'(C, f)
\to[F'(f_0)] F'(C', \id_{u(C')}) = \eta_\sP(D', F')
\]

We will also construct a functor $\eta'_\sP : \sP \to u^*u_*\sP$.
First, for each morphism $m : C' \to C$ in $\sC$ and object $P \in \Ob{\sP}$ above $C$,
we choose a Cartesian lift $\wh{m} : \wh{C'} \to P$.
Then, consider any $P \in \Ob{\sP}$. We construct as follows a functor
\[
\Phi_P : u(\sC)_{/u(p(P))} \to \sP
\]
For an object $(C, f)$ with $f : u(C) \to u(p(P))$,
we have $f = u(f_0)$ for some unique $f_0 : C \to p(P)$ in $\sC$ since $u$ is fully faithful.
Take the chosen Cartesian lift $\wh{f_0} : \wh{C} \to p(P)$ and define $\Phi_P(C, f) := \wh{C}$.
For a morphism $g : (C, f) \to (C', f')$ in $u(\sC)_{/u(p(P))}$, we have $f' \circ u(g) = f$
in $\sD$, but since $u$ is fully faithful, $f' = u(f'_0), f= u(f_0)$ for some unique morphisms
$f'_0 : C' \to p(P)$ and $f_0 : C \to p(P)$ in $\sC$,
and $u(f'_0 \circ g) = u(f'_0) \circ u(g) = u(f_0)$ implies $f'_0 \circ g = f_0$.
By the Cartesian property of the chosen lift $\wh{f'_0}$, we have a unique lift
$\wh{g} : \wh{C} \to \wh{C'}$ of $g$ such that $\wh{f'_0} \circ \wh{g} = \wh{f_0}$.
We set $\Phi_P(g) := \wh{g}$. The uniqueness of the lifts $\wh{g}$ guarantees that
$\Phi_P$ is a functor $u(\sC)_{/u(p(P))} \to \sP$. Thus, $(u(p(P)), \Phi_P)$ is an object of
$u_*\sP$ over $u(p(P))$.

Next, consider a morphism $r : P_2 \to P_1$ in $\sP$. We will construct a natural transformation
\[
\alpha_r : \Phi_{P_1} \circ u(\sC)_{/u(p(r))} \To \Phi_{P_2} : u(\sC)_{/u(p(P_2))} \to \sP
\]
Let $(C_2, f_2)$ be an object of $u(\sC)_{/u(p(P_2))}$. Then, $f_2 = u(f_{2, 0})$ for a unique
$f_{2, 0} : C_2 \to p(P_2)$ in $\sC$ and
$u(\sC)_{/u(p(r))}(C_2, f_2) = (C_2, u(p(r)) \circ f_2)= (C_2, u(p(r)) \circ u(f_{2, 0}))
= (C_2, u(p(r) \circ f_{2, 0}))$.
Then, the domain of the chosen Cartesian lift $\wh{p(r) \circ f_{2, 0}}$ of
$p(r) \circ f_{2, 0}$ with codomain $P_1$ is $\Phi_{P_1}(u(\sC)_{/u(p(r))}(C_2, f_2))$.
$\Phi_{P_2}(C_2, f_2)$, on the other hand, is
the domain of the chosen Cartesian lift of $f_{2, 0}$ with codomain $P_2$.
At the same time, by
\cite[\href{https://stacks.math.columbia.edu/tag/003V}{Lemma 003V}]{stacks-project},
$r \circ \wh{f_{2, 0}}$ is Cartesian.
This provides a unique isomorphism $\alpha_{r, C_2, f_2}$ over $\id_{C_2}$
making the following diagram commute:
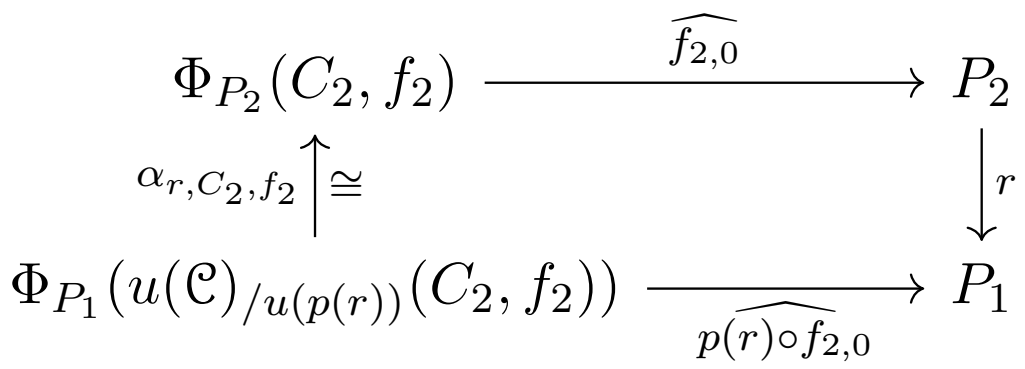
\[\begin{tikzcd}[column sep=huge]
\Phi_{P_2}(C_2, f_2) \ar[r, "\wh{f_{2, 0}}"] \ar[from=d, "\alpha_{r, C_2, f_2}" left, "\cong" right] &
P_2 \ar[d, "r"] \\
\Phi_{P_1}(u(\sC)_{/u(p(r))}(C_2, f_2)) \ar[r, "\wh{p(r) \circ f_{2, 0}}" below] & P_1
\end{tikzcd}\]
For the naturality of $\alpha_{r, C_2, f_2}$ in $(C_2, f_2)$, observe that
given any commutative diagram in $\sC$ of the form:
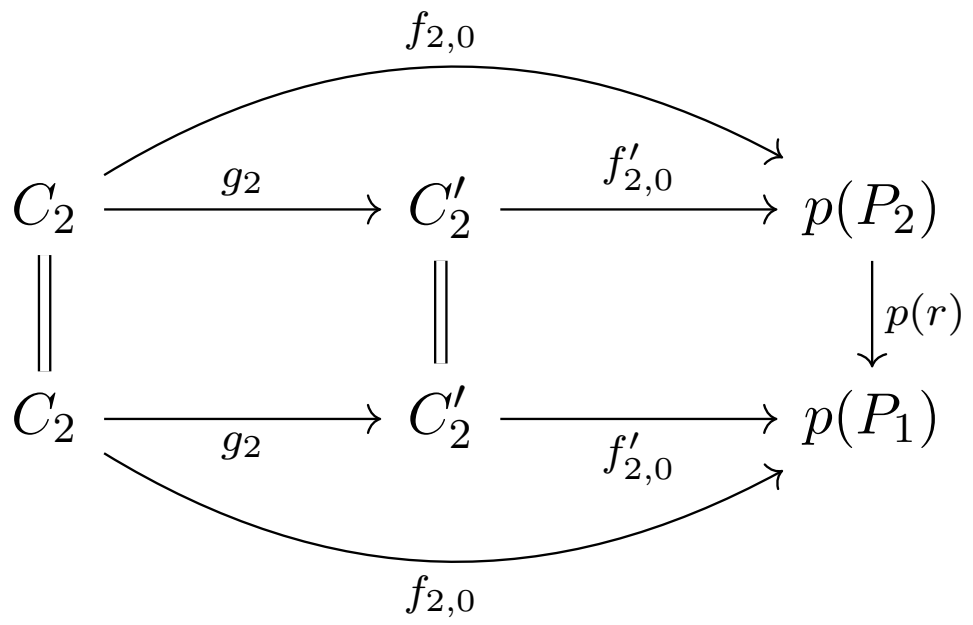
\[\begin{tikzcd}[column sep=huge]
C_2 \ar[r, "g_2"] \ar[rr, bend left, "f_{2, 0}"] \ar[d, equal] &
C_2' \ar[r, "f'_{2, 0}"] \ar[d, equal] &
p(P_2) \ar[d, "p(r)"] \\
C_2 \ar[r, "g_2" below] \ar[rr, bend right, "f_{2, 0}" below] &
C_2' \ar[r, "f'_{2, 0}" below] &
p(P_1)
\end{tikzcd}\]
by the Cartesian property of morphisms involved,
we get a unique commutative diagram in $\sP$ of the form on the left:
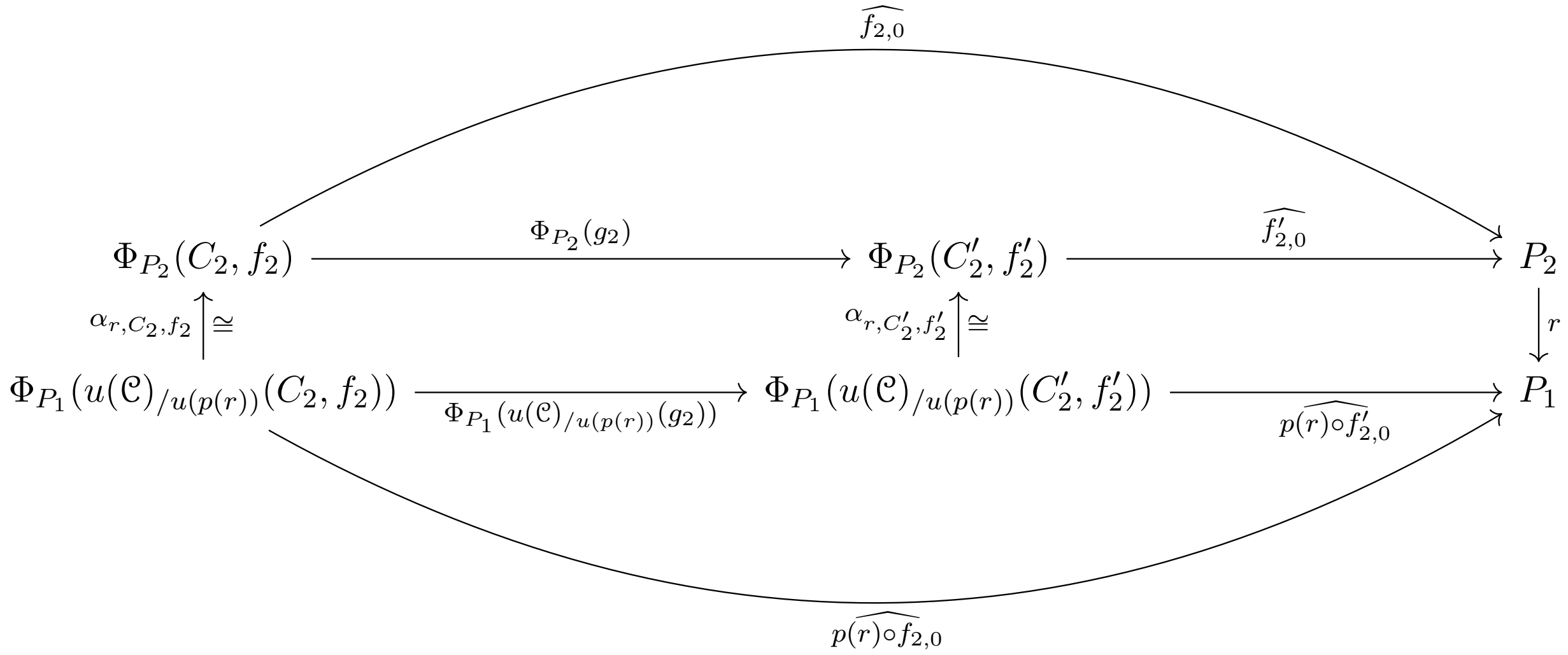
\[\begin{tikzcd}[column sep=9em]
\Phi_{P_2}(C_2, f_2) \ar[r, "\Phi_{P_2}(g_2)"] \ar[rr, bend left, "\wh{f_{2, 0}}"]
    \ar[from=d, "\alpha_{r, C_2, f_2}" left, "\cong" right] &
\Phi_{P_2}(C_2', f_2') \ar[r, "\wh{f'_{2, 0}}"] \ar[from=d, "\alpha_{r, C'_2, f'_2}", "\cong" right] &
P_2 \ar[d, "r"] \\
\Phi_{P_1}(u(\sC)_{/u(p(r))}(C_2, f_2)) \ar[r, "\Phi_{P_1}(u(\sC)_{/u(p(r))}(g_2))" below]
    \ar[rr, bend right, "\wh{p(r) \circ f_{2, 0}}" below] &
\Phi_{P_1}(u(\sC)_{/u(p(r))}(C'_2, f'_2)) \ar[r, "\wh{p(r) \circ f'_{2, 0}}" below] &
P_1
\end{tikzcd}\]

If $r : P_2 \to P_1$ is the identity morphism of an object, then $\alpha_{r, C_2, f_2}$
is necessarily the identity morphism. Next, if we have a chain of morphisms
$P_3 \to[s] P_2 \to[r] P_1$, then we have a pasting of commutative squares:
\[\begin{tikzcd}[column sep=huge]
\Phi_{P_3}(C_3, f_3)
    \ar[r, "\wh{f_{3, 0}}"]
    \ar[from=d, "\alpha_{s, C_3, f_3}" left, "\cong" right] &
P_3 \ar[d, "s"] \\
\Phi_{P_2}(u(\sC)_{u(p(s))}(C_3, f_3))
    \ar[r, "\wh{p(s) \circ f_{3, 0}}"]
    \ar[from=d, "\alpha_{r, C_2, p(s) \circ f_{3, 0}}" left, "\cong" right] &
P_2 \ar[d, "r"] \\
\Phi_{P_1}(u(\sC)_{/u(p(r))}(u(\sC)_{/u(p(s))}(C_3, f_3)))
    \ar[r, "\wh{p(r) \circ p(s) \circ f_{3, 0}}" below] &
P_1
\end{tikzcd}\]
Since, $p(r) \circ p(s) \circ f_{3, 0} = p(r \circ s) \circ f_{3, 0}$, the composite
of the left vertical arrows must be equal to $\alpha_{r \circ s, C_3, f_3}$ by the Cartesian
property of $r \circ s \circ \wh{f_{3, 0}}$. This shows that the assignment:
\[
P \mapsto \Phi_P, r \mapsto \alpha_r
\]
is a functor $\sP \to u_*\sP$. Furthermore, the following diagram of categories strictly commutes:
\[\begin{tikzcd}
\sP \ar[r] \ar[d, "p" left] & u_*\sP \ar[d, "u_*p"] \\
\sC \ar[r, "u" below] & \sD
\end{tikzcd}\]
This provides a morphism of prestacks over $\sC$: $\eta'_\sP : \sP \to u^*u_*\sP$.
\end{cns}

\begin{prop}\label{prop:prest-pull-push-equiv}
In the context of \cref{cns:unit-morphism}, $\eta_\sP$ is an equivalence
of prestacks:
\[
u^*u_*\sP \to[\simeq] \sP
\]
\end{prop}
\begin{proof}
We first show that $\eta_\sP$ is fully faithful.
Let $(f, \alpha), (g, \beta) : (D, F) \to (D', F')$ be morphisms in $u^*u_*\sP$ such that
$\eta_\sP(f, \alpha) = \eta_\sP(g, \beta)$. Since these morphisms are in $u^*u_*\sP$,
there exist morphisms $f_0, g_0 : C \to C'$ in $\sC$ with $u(f_0) = f, u(g_0) = g$. The hypothesis
is $F'(f_0) \circ \alpha_{C, \id_{u(C)}}^{-1} = F'(g_0) \circ \beta_{C, \id_{u(C)}}^{-1}$,
which implies:
\begin{align*}
& p(F'(f_0)) \circ p(\alpha_{C, \id_{u(C)}}^{-1})
    = p(F'(g_0)) \circ p(\beta_{C, \id_{u(C)}}^{-1}) \\
\implies& f_0 \circ \id_C = g_0 \circ \id_C \\
\implies& f_0 = g_0 \\
\implies& f = u(f_0) = u(g_0) = g
\end{align*}
By \cite[\href{https://stacks.math.columbia.edu/tag/003V}{Lemma 003V}]{stacks-project},
$F'(f_0) = F'(g_0)$ is a Cartesian lift of $f_0 = g_0$, and hence
we must have
\[
\alpha_{C, \id_{u(C)}}^{-1} = \beta_{C, \id_{u(C)}}^{-1}
\implies \alpha_{C, \id_{u(C)}} = \beta_{C, \id_{u(C)}}
\]

For essential surjectivity, consider an arbitrary $P \in \Ob{\sP}$.
Then, $\eta_\sP(\eta'_\sP(P)) = \Phi_P(\id_{u(p(P))})$ is the domain
of the chosen Cartesian lift $\wh{\id_{u(p(P))}}$ of $\id_{u(p(P))}$. Since $\sP$ is a prestack,
$\wh{\id_{u(p(P))}}$ is an isomorphism $\eta_\sP(\eta'_\sP(P)) \to P$ lying over $\id_{p(P)}$.
We have, in particular, shown that $\eta'_\sP$ is fibrewise fully faithful and essentially
surjective.
We can now apply \cite[\href{https://stacks.math.columbia.edu/tag/003Z}{Lemma 003Z}]{stacks-project}
to deduce that $\eta_\sP$ is an equivalence of prestacks.
\end{proof}

\subsection{Pullback and Pushforward of Stacks}
\label{subsec:pullback-pushfwd-st}

\begin{prop}\label{prop:cont-func}
Let $(\sC, J), (\sD, K)$ be sites; $p : \sP \to \sD$, a stack and $u : \sC \to \sD$, a functor.
If, for every covering family $\set{c_i : C_i \to C}_{i \in I}$ in $(\sC, J)$, $u$ satisfies:
\begin{enumerate}
\item $\set{u(c_i)}_{i \in I}$ is a covering family in $(\sD, K)$.
\item The canonical morphism $u(C_i \times_C C_j) \to u(C_i) \times_{u(C)} u(C_j)$ is an
isomorphism.
\end{enumerate}
then $u^*p : u^*\sP \to \sC$ is a stack. In particular, the strict functor of
\cref{prop:prest-pb} restricts to a strict functor:
\[
u^* : \St_{\sD, K} \to \St_{\sC, J}
\]
\end{prop}
\begin{proof}
Let $P_i, i \in I$ be objects in $u^*\sP$ over some $C_i \in \Ob{\sC}$ --- that is,
$p(P_i) = u(C_i) \in \Ob{\sD}$, $C_i \in \Ob{\sC}$ for each $i \in I$.
Suppose we have isomorphisms $\alpha_{ij} : P_i|C_{ij} \to[\cong] P_j|C_{ij}$ in
$u^*\sP$ such that $\alpha_{jk}|C_{ijk} \circ \alpha_{ij}|C_{ijk} = \alpha_{ik}|C_{ijk}$.
This means the $\alpha_{ij}$ are isomorphisms in $\sP$ over $\id_{u(C_{ij})}$.
The cocycle condition is an equality of morphisms in $\sP$ over $u(C_{ijk})$.
Taking Cartesian lifts of the isomorphisms from (ii) provides isomorphisms
$\rho_{ij} : P_{i}|u(C_i) \times_{u(C)} u(C_j) \to[\cong] P_{i}|u(C_{ij})$ in $\sP$.
Thus, there is a unique isomorphism
$\alpha'_{ij} : P'_{i}|u(C_{ij}) \to P'_{j}|u(C_{ij})$ over
$\id_{u(C_i) \times_{u(C)} u(C_j)}$ such that the following diagram commutes:
\[\begin{tikzcd}
P_{i}|u(C_i) \times_{u(C)} u(C_j) \ar[r, "\rho_{ij}"] \ar[d, "\alpha'_{ij}" left] &
P_{i}|u(C_{ij}) \ar[d, "\alpha_{ij}"] \\
P_{j}|u(C_i) \times_{u(C)} u(C_j) \ar[r, "\rho_{ji}" below] &
P_{j}|u(C_{ij})
\end{tikzcd}\]
We have isomorphisms $u(C_i) \times_{u(C)} u(C_j) \times_{u(C)} u(C_k) \to u(C_{ijk})$
induced by the binary ones from hypothesis (ii). By taking Cartesian lifts of these, we can
restrict the commutative squares above to $u(C_i) \times_{u(C)} u(C_j) \times_{u(C)} u(C_k)$
on the left and $u(C_{ijk})$ on the right, and pasting two such squares, we can show that the
$\alpha'_{ij}$ satisfy the cocycle condition.
By hypothesis (i), $\set{u(C_i) \to[u(c_i)] u(C)}_{i \in I}$ is a covering family in $\sD$ and
hence, we have a descent datum for an object in $\sP$ over $u(C)$.
Since $\sP$ is a stack by assumption, we have an object $P \in \Ob{\sP}$ over $u(C)$ and
isomorphisms $\phi_i : P|u(C_i) \to P_i$ such that
$\phi_j|u(C_i) \times_{u(C)} u(C_j) = \alpha'_{ij} \circ \phi_i|u(C_i) \times_{u(C)} u(C_j)$.
Using hypothesis (ii) again, we can show that
$\phi_j|u(C_{ij}) = \alpha_{ij} \circ \phi_i|u(C_{ij})$. Thus, $P$ is an object in
$u^*\sP$ solving the descent problem for the descent datum
$(\set{P_i}_{i \in I}, \set{\alpha_{ij}}_{i, j \in I})$. The argument for collation of morphisms
in $u^*\sP$ is similar.
\end{proof}

\begin{defn}[Continuous Functor]
A functor satisfying the conditions of \cref{prop:cont-func} is called a continuous functor.
\end{defn}

\begin{prop}\label{prop:cocont-func}
Let $(\sC, J), (\sD, K)$ be sites; $p : \sP \to \sC$, a stack and $u : \sC \to \sD$, a functor.
Suppose, for every object $C \in \Ob{\sC}$ and every covering family
$\set{d_i : D_i \to u(C)}_{i \in I}$ in $(\sD, K)$,
there exists a covering family $\set{c_q : C_q \to C}_{q \in Q}$ and a function
$\gamma : Q \to I$, such that $u(c_q)$ factors through $d_{\gamma(q)}$.
Then, $u_*p : u_*\sP \to \sD$ is a stack.
\end{prop}
\begin{proof}
Consider a covering family $\set{d_i : D_i \to D}_{i \in I}$ in $\sD$ along with
objects $(D_i, F_i)$ over $D_i$ in $u_*\sP$ and isomorphisms
$(\id_{D_i \times_D D_j}, \alpha_{ij}) :
(D_i, F_i)|D_i \times_D D_j \to[\cong] (D_j, F_j)|D_i \times_D D_j$
such that $\alpha_{jk}|D_{ijk} \circ \alpha_{ij}|D_{ijk} = \alpha_{ik}|D_{ijk}$ for each
$i, j \in I$. Then, consider $(C, f) \in \Ob{u(\sC)_{/D}}$. We have a covering family
\[
\set{d_i \times_D u(C) : D_i \times_D u(C) \to u(C)}_{i \in I}
\]
in $\sD$.
By hypothesis, we have a covering family $\set{c_q : C_q \to C}_{q \in Q}$ and a function
$\gamma : Q \to I$ such that $u(c_q)$ factors as
$u(C_q) \to[e_q] D_{\gamma(q)} \times_D u(C) \to[d_{\gamma(q)} \times_D u(C)] u(C)$.
Then,
\[
F(C_q, u(c_q)) := F_{\gamma(q)}(C_q, u(C_q) \to[e_q] D_{\gamma(q)} \times_D u(C) \to D_{\gamma(q)})
\]
is an object of $\sP$ over $C_q$ for each $q \in Q$. We have a pasting of Cartesian squares:
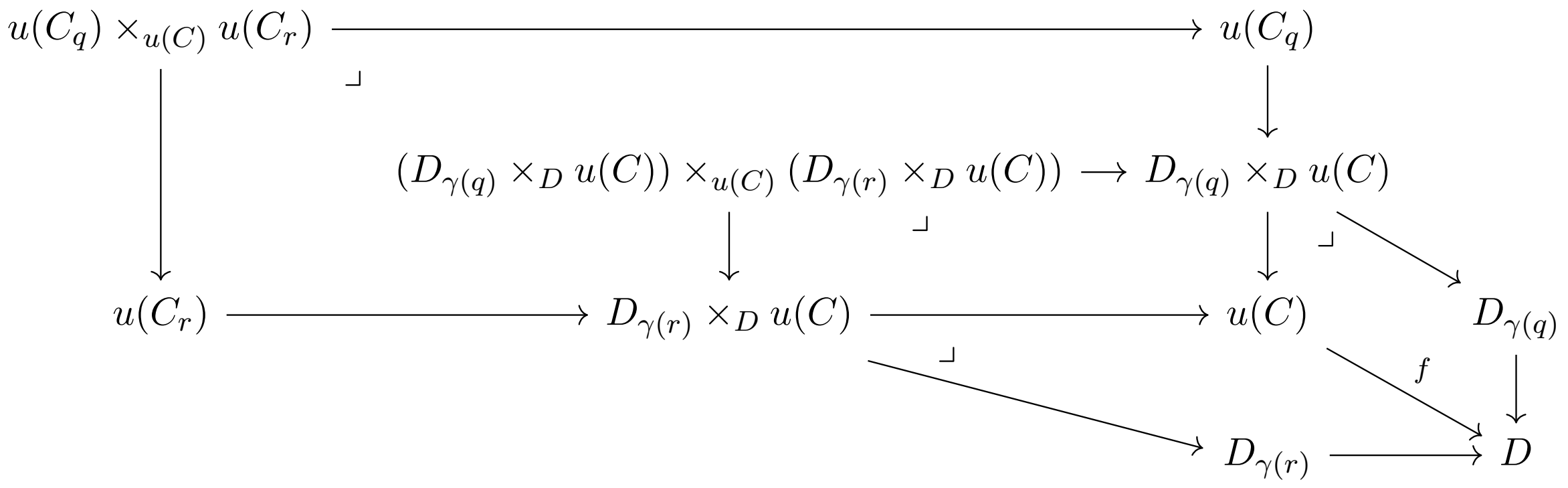
\[\begin{tikzcd}[column sep=small]
u(C_q) \times_{u(C)} u(C_r)
    \ar[rr]
    \ar[dd]
    \ar[rd, phantom, "\lrcorner" very near start] & &
u(C_q)
    \ar[d] \\ &
(D_{\gamma(q)} \times_{D} u(C)) \times_{u(C)} (D_{\gamma(r)} \times_{D} u(C))
    \ar[r]
    \ar[d]
    \ar[rd, phantom, "\lrcorner" very near start] &
D_{\gamma(q)} \times_{D} u(C)
    \ar[d]
    \ar[rd]
    \ar[rdd, phantom, "\lrcorner" very near start] \\
u(C_r)
    \ar[r] &
D_{\gamma(r)} \times_{D} u(C)
    \ar[r]
    \ar[rd, shift right]
    \ar[rrd, phantom, "\lrcorner" very near start] &
u(C)
    \ar[rd, "f"] &
D_{\gamma(q)}
    \ar[d] \\ & &
D_{\gamma(r)}
    \ar[r] &
D
\end{tikzcd}\]
The universal property of $D_{\gamma(q)} \times_D D_{\gamma(r)}$ provides a canonical morphism
$u(C_q) \times_{u(C)} u(C_r) \to D_{\gamma(q)} \times_D D_{\gamma(r)}$.
Thus, we get a morphism
$u(C_{qr}) \to u(C_q) \times_{u(C)} u(C_r) \to D_{\gamma(q)} \times_D D_{\gamma(r)}$, at which
we can take the component of the isomorphism $\alpha_{\gamma(q), \gamma(r)}$ to get an isomorphism
$\alpha'_{qr} : F(C_q, u(c_q))|_{C_{qr}} \to F(C_r, u(c_r))|_{C_{qr}}$. These isomorphisms satisfy
the cocycle condition by a similar argument as in the proof of \cref{prop:cont-func}. Since
$\sP$ is a stack by assumption, we have an object $F(C, f)$ in $\sP$ over $C$ by collating
the $F(C_q, u(c_q))$.

Now, let $g : (C, f) \to (C', f')$ be a morphism in $u(\sC)_{/D}$. Consider choices of covering
families
$\set{c_q : C_q \to C}_{q \in Q}$ and $\set{c'_s : C'_s \to C'}_{s \in S}$ in $\sC$ as in the
previous paragraph, where, in addition to $\gamma$ and the morphisms $e_q$,
we also have a function $\gamma' : S \to I$ and a factorization
$u(C'_s) \to[e'_s] D_{\gamma'(s)} \times_{D} u(C') \to[d_{\gamma'(s)} \times_D u(C)] u(C')$.
Then, $\set{c_q \times_C c'_s : C_q \times_C C'_s \to C}_{(q, s) \in Q \times S}$
is a covering family in $\sC$. Then, the composite morphisms
$u(C_q \times_C C'_s) \to u(C_q) \to D_{\gamma(q)} \times_D u(C) \to D_{\gamma(q)}$
and $u(C_q \times_C C'_s) \to u(C'_s) \to D_{\gamma'(s)} \times_D u(C') \to D_{\gamma'(s)}$
provide a morphism $u(C_q \times_C C'_s) \to D_{\gamma(q)} \times_D D_{\gamma'(s)}$ by the universal
property of fibre products. Then, we have the following commutative diagrams:
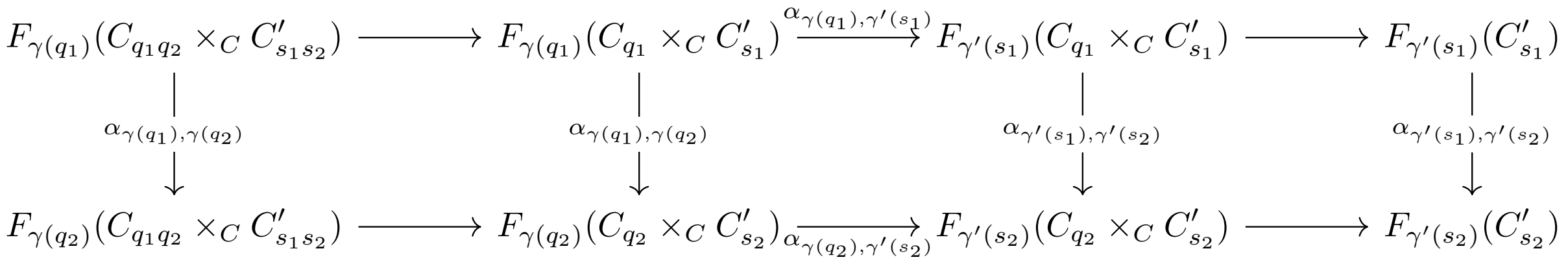
\[\begin{small}\begin{tikzcd}[row sep=huge, column sep=large]
F_{\gamma(q_1)}(C_{q_1q_2} \times_C C'_{s_1s_2}) \ar[r]
    \ar[d, "\alpha_{\gamma(q_1), \gamma(q_2)}" description] &
F_{\gamma(q_1)}(C_{q_1} \times_C C'_{s_1}) \ar[r, "\alpha_{\gamma(q_1), \gamma'(s_1)}"]
    \ar[d, "\alpha_{\gamma(q_1), \gamma(q_2)}" description] &
F_{\gamma'(s_1)}(C_{q_1} \times_C C'_{s_1}) \ar[r]
    \ar[d, "\alpha_{\gamma'(s_1), \gamma'(s_2)}" description] &
F_{\gamma'(s_1)}(C'_{s_1})
    \ar[d, "\alpha_{\gamma'(s_1), \gamma'(s_2)}" description] \\
F_{\gamma(q_2)}(C_{q_1q_2} \times_C C'_{s_1s_2}) \ar[r] &
F_{\gamma(q_2)}(C_{q_2} \times_C C'_{s_2}) \ar[r, "\alpha_{\gamma(q_2), \gamma'(s_2)}" below] &
F_{\gamma'(s_2)}(C_{q_2} \times_C C'_{s_2}) \ar[r] &
F_{\gamma'(s_2)}(C'_{s_2})
\end{tikzcd}\end{small}\]
where the unlabelled arrows are obtained by applying the appropriate functors $F_i$
to projections of fibre products, for appropriate choices of $i$. The commutativity of the
squares in this diagram is ensured by the naturality of the $\alpha_{ij}$ and the cocycle
condition satisfied by them. Next, we also have a commutative triangle:
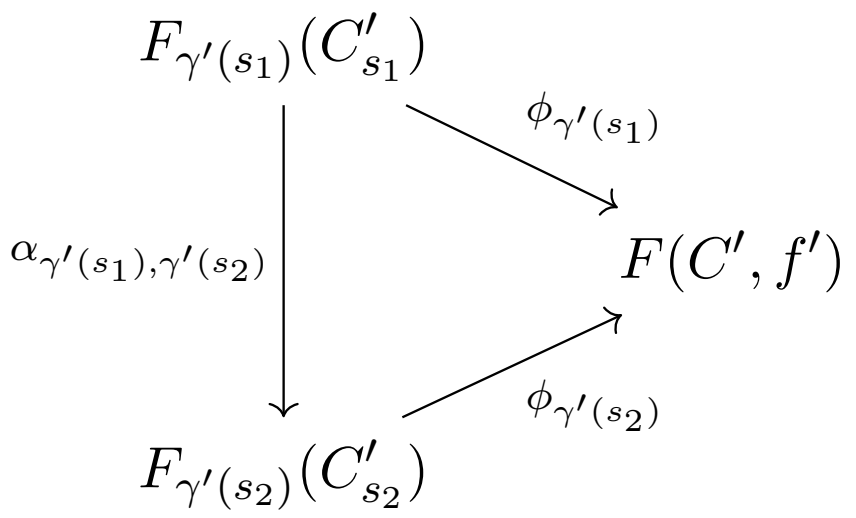
\[\begin{tikzcd}
F_{\gamma'(s_1)}(C'_{s_1})
    \ar[rd, "\phi_{\gamma'(s_1)}" above right]
    \ar[dd, "\alpha_{\gamma'(s_1), \gamma'(s_2)}" left] & \\
& F(C', f') \\
F_{\gamma'(s_2)}(C'_{s_2})
    \ar[ru, "\phi_{\gamma'(s_2)}" below right]
\end{tikzcd}\]
where the morphisms $\phi_{-}$ are provided by the collation of objects described in the previous
paragraph.
Putting these diagrams together gives a descent datum for a morphism
$F(C, f) \to F(C', f')$ which we define to be $F(g)$.

The fact that the assignment $F$ is
a functor $u(\sC)_{/D} \to \sP$ follows from pasting the commutative diagrams above
providing the values $F(g)$ on morphisms and the uniqueness of collated morphisms.
\end{proof}

\begin{defn}[Cocontinuous Functor]
A functor satisfying the conditions of \cref{prop:cocont-func} is called a cocontinuous functor.
\end{defn}

\subsection{A Comparison Theorem}
\label{subsec:comp-thm}

\begin{thm}[Comparison Theorem for Stacks]\label{thm:comp-thm}
Let $(\sC, J), (\sD, K)$ be sites and $u : \sC \to \sD$ a functor such that
the following hold:
\begin{enumerate}
\item $u$ is fully faithful.

\item For all covering families $\set{c_i : C_i \to C}_{i \in I}$ in $J$, the universal map
$u(C_i \times_C C_j) \to u(C_i) \times_{u(C)} u(C_j)$ is an isomorphism in $\sD$.

\item For every covering family $\set{c_i : C_i \to C}_{i \in I}$ in $J$,
$\set{u(c_i) : u(C_i) \to u(C)}_{i \in I}$ is a covering family in $K$.

\item For every object $C \in \Ob{\sC}$ and every covering family
$\set{d_i : D_i \to u(C)}_{i \in I}$ in $(\sD, K)$,
there exists a covering family $\set{c_q : C_q \to C}_{q \in Q}$ and a function
$\gamma : Q \to I$, such that $u(c_q)$ factors through $d_{\gamma(q)}$.

\item For every object $D$ in $\sD$, there is a covering family
$\set{d_i : u(C_i) \to D}_{i \in I}$ in $K$.
\end{enumerate}
Then, the strict $2$--functor $u^* : \PSt_\sD \to \PSt_\sC$ given by $2$--pullback of categories
along $u$ restricts to an equivalence of $2$--categories:
\[
u^* : \St_{\sD, K} \to[\simeq] \St_{\sC, J}
\]
\end{thm}
\begin{proof}
The fact that $u^*$ restricts to a $2$--functor $\St_{\sD, K} \to \St_{\sC, J}$ is
\cref{prop:cont-func} in combination with hypotheses (ii) and (iii).
We will first show that for two stacks $\sP, \sQ \in \St_{\sD, K}$ the $1$--functor of Hom
categories $u^* : \St_{\sD, K}(\sP, \sQ) \to \St_{\sC, J}(u^*\sP, u^*\sQ)$ is fully faithful.
Let $\alpha, \beta : f \To g : \sP \to \sQ$ be two $2$--morphisms of stacks, such that
$u^*\alpha = u^*\beta$. Then, let $P \in \Ob{\sP}$ over $D \in \Ob{\sD}$, and consider the
components $\alpha_P, \beta_P$. To see these are equal, use hypothesis (v) and choose a covering
family $\set{c_i : u(C_i) \to D}_{i \in I}$, where each $C_i \in \Ob{\sC}$, and consider the
restrictions $\alpha_P|u(C_i), \beta_P|u(C_i) : f(P)|u(C_i) \to g(P)|u(C_i)$. Since
$\sQ$ is a prestack and $f$ is a morphism of prestacks,
there are isomorphisms making the following diagram commute:
\begin{equation}\label{eqn:rest-nat-iso}
\begin{tikzcd}
f(P)|{u(C_i)} \ar[r, "\alpha_P|u(C_i)"] \ar[d, "\cong" left] & g(P)|u(C_i) \ar[d, "\cong"] \\
f(P|u(C_i)) \ar[r, "\alpha_{P|u(C_i)}" below] & g(P|u(C_i))
\end{tikzcd}
\end{equation}
By hypothesis, $\alpha_{P|u(C_i)} = \beta_{P|u(C_i)}$ so that $\alpha_P|u(C_i) = \beta_P|u(C_i)$,
by the commutativity of the above diagram. This shows that $\alpha_P$ and $\beta_P$ solve the
same descent problem for a morphism in $\sQ$. Since $\sQ$ is a stack, we must have
$\alpha_P = \beta_P$, by the uniqueness of collated morphisms.
Thus, the functor $u^* : \St_{\sD, K}(\sP, \sQ) \to \St_{\sC, J}(u^*\sP, u^*\sQ)$ is faithful.

Now, suppose we have a $2$--morphism $\alpha : u^*f \To u^*g : u^*\sP \to u^*\sQ$ and consider
any object $P \in \Ob{\sP}$. Choose a covering family $\set{u(C_i) \to p(P)}_{i \in I}$
using hypothesis (v).
Then, $P|u(C_i)$ is an object over $u(C_i)$ in $\sP$, and hence is an object over $C_i$ in
$u^*\sP$. Thus, we have the component morphisms $\alpha_{P|u(C_i)} : f(P|u(C_i)) \to g(P|u(C_i))$
in $u^*\sQ$ lying over $\id_{C_i}$. We have similar commutative diagrams as \ref{eqn:rest-nat-iso}
above giving morphisms $\alpha'_{P, i} : f(P)|u(C_i) \to g(P)|u(C_i)$
which, composed with the morphisms $g(P)|u(C_i) \to g(P)$ form a descent datum for a morphism
in $\sQ$ over $\id_{p(P)}$, and, since $\sQ$ is a stack, we have a collation
$\alpha'_P : f(P) \to g(P)$.
We would like to show that $\alpha'_P$ is natural in $P$: that is, for any morphism
$h : P \to P'$ in $\sP$, we need to show the following diagram commutes:
\[\begin{tikzcd}
f(P) \ar[r, "\alpha'_P"] \ar[d, "f(h)" left] & g(P) \ar[d, "g(g)"] \\
f(P') \ar[r, "\alpha'_P" below] & g(P')
\end{tikzcd}\]
Choose a covering family
$\set{u(C'_r) \to p(P')}_{r \in R}$ using hypothesis (v) which provides, by pullback, a covering
family $\set{u(C_i) \times_{p(P')} u(C'_r) \to p(P)}_{(i, r) \in I \times R}$. Using hypothesis
(v) on the $u(C_i) \times_{p(P')} u(C'_r)$,
we can produce a covering family
$\set{u(C^i_r) \to p(P)}_{(i, r) \in I \times R}$ such that we have a commutative diagram:
\[\begin{tikzcd}
u(C^i_r) \ar[r] \ar[d] & u(C'_r) \ar[d] \\
p(P) \ar[r] & p(P')
\end{tikzcd}\]
The naturality of $\alpha$ ensures the following diagram commutes:
\begin{equation}\label{eqn:rest-nat-iso-2}
\begin{tikzcd}
f(P|u(C^i_r)) \ar[r, "\alpha_{P|u(C^i_r)}"] \ar[d, "f(h|u(C^i_r))" left] &
g(P|u(C^i_r)) \ar[d, "g(h|u(C^i_r))"] \\
f(P'|u(C^i_r)) \ar[r, "\alpha_{P'|u(C^i_r)}" below] & g(P'|u(C^i_r))
\end{tikzcd}
\end{equation}
Then, we have unique isomorphisms over $\id_{u(C^i_r)}$ making all faces
of the following cube except, a priori, the front face commute:
\[\begin{tikzcd}
f(P|u(C^i_r))
  \ar[rr, "\alpha_{P|u(C^i_r)}"]
  \ar[rd, "\cong" description]
  \ar[dd, "f(h|u(C^i_r))" left] & &
g(P|u(C^i_r)) \ar[rd, "\cong" description] \ar[dd, "g(h|u(C^i_r))" near end] & \\ &
f(P)|u(C^i_r) \ar[rr, crossing over, "\alpha'_{P}|u(C^i_r)" near start] & &
g(P)|u(C^i_r) \ar[dd, "g(h)|u(C^i_r)"] \\
f(P'|u(C^i_r)) \ar[rr, "\alpha_{P'|u(C^i_r)}" below, near end] \ar[rd, "\cong" description] & &
g(P'|u(C^i_r)) \ar[rd, "\cong" description] & \\ &
f(P')|u(C^i_r)
  \ar[rr, "\alpha'_{P'}|u(C^i_r)" below]
  \ar[from=uu, crossing over, "f(h)|u(C^i_r)" left, near start] & &
g(P')|u(C^i_r) \\
\end{tikzcd}\]
However, this implies the front face also commutes. This, in turn, shows that
the two composites in diagram \ref{eqn:rest-nat-iso-2} solve the same descent problem
for a morphism in $\sQ$. By construction, we also have $u^*\alpha' = \alpha$. Thus, the functor
$u^* : \St_{\sD, K}(\sP, \sQ) \to \St_{\sC, J}(u^*\sP, u^*\sQ)$ is full.

Next, consider a $1$--morphism $f : u^*\sP \to u^*\sQ$ in $\St_{\sC, J}$. Let $P \in \Ob{\sP}$
and choose a covering family $\set{u(C_i) \to p(P)}_{i \in I}$ using hypothesis (v). Then,
$P|u(C_i)$ is an object of $\sP$ over $u(C_i)$ and hence an object of $u^*\sP$. So, we have
objects $f(P|u(C_i))$ in $u^*\sQ$, which are just objects of $\sQ$ over $u(C_i)$.
Since the $P|u(C_i)$ are restrictions of a single object on $\sP$,
the $f(P|u(C_i))$ form a descent datum for an object in $\sQ$ over $p(P)$,
whose collation we define to be $f'(P)$, using the fact that $\sQ$ is a stack.
Now, consider a morphism $h : P \to P'$ in $\sP$.
Then, similar to the previous paragraph we have morphisms
$f(h|u(C^i_r)) : f(P|u(C^i_r)) \to f(P'|u(C^i_r))$ which provide a descent datum for a morphism
$f'(P) \to f'(P')$, since the $h|u(C^i_r)$ are restrictions of a single morphism in $\sP$
and $f$ is a morphism of prestacks. We call the collated morphism $f'(h)$. For any other morphism
$k : P' \to P''$ in $\sP$, we can paste the diagrams of the descent problems for
$f'(h)$ and $f'(k)$ and use the functoriality of $f$ to deduce that $f'(k) \circ f'(h)$
and $f'(k \circ h)$ solve the same descent problem for a morphism in $\sQ$ and hence must be equal,
since $\sQ$ is a stack. We can similarly show that $f'(\id_{P})$ and $\id_{f'(P)}$ solve the same
descent problem for a morphism in $\sQ$, and hence must be equal. Thus $f'$ is a functor
$\sP \to \sQ$. It is not difficult to see that $f'$ is a morphism of prestacks
(and hence, of stacks). Finally, $u^*f' = f$ by construction. Thus,
$u^* : \St_{\sD, K}(\sP, \sQ) \to \St_{\sC, J}(u^*\sP, u^*\sQ)$ is essentially surjective.

Given an object $\sP \in \Ob{\St_{\sC, J}}$, \cref{prop:cocont-func} in combination
with hypotheses (i) and (iv) shows that $u_*\sP$ is in $\St_{\sD, K}$.
\cref{prop:cont-func} in combination with hypothesis (i) and (ii) shows that
$u^*u_*\sP$ is in $\St_{\sC, J}$.
\cref{prop:prest-pull-push-equiv} in combination with hypothesis (i) shows that
$u^*u_*\sP$ and $\sP$ are equivalent as prestacks and hence also as stacks.
Thus, $u^* : \St_{\sD, K} \to \St_{\sC, J}$ is essentially surjective.
\end{proof}

\subsection{Representable Objects and Morphisms}
\label{subsec:repr-obj-morphisms}

\begin{defn}[Prestack Represented by an Object]\label{defn:rep-prest}
Let $\sC$ be a category. Given an object $C \in \Ob{\sC}$, the slice category
$\sC_{/C}$ admits a canonical functor to $\sC$. Then, $\sC_{/C}$ is called the canonical prestack
represented by $C$. If $\sP \to \sC$ is any other prestack such that there exists an equivalence
$\sP \simeq \sC_{/C}$ of categories over $\sC$,
we say that $\sP$ is represented by $C$. If no choice of
$C$ is specified but one exists, we say that $\sP$ is a representable prestack.
\end{defn}

\begin{prop}
A representable prestack is indeed a prestack.
\end{prop}
\begin{proof}
Use notation as in \cref{defn:rep-prest}. Consider an object $f : C' \to C$ in
$\sC_{/C}$ over $C' \in \Ob{\sC}$, and a morphism $g : C'' \to C'$ in $\sC$.
Then, $f \circ g : C'' \to C$ is an object in $\sC_{/C}$ over $C'' \in \Ob{\sC}$
and $g : C'' \to C'$ is a morphism $f \circ g \to f$ in $\sC_{/C}$. Now, suppose
we have another object $k : D \to C$ in $\sC$ and a morphism $l : k \to f$ in $\sC_{/C}$
so that $l$ is a morphism $D \to C'$ in $\sC$ with $f \circ l = k$. Also, assume
that we have a morphism $m : D \to C''$ in $\sC$ such that $f \circ g \circ m = k$
in $\sC$. Then, $m$ is a morphism $k \to f \circ g$ in $\sC_{/C}$ as well and
$f \circ g \circ m = f \circ l$ implies $g \circ m = l$ as morphisms in $\sC_{/C}$.
Any other morphism $m'$ satisfying these conditions must project to $m$ in $\sC$
and hence must be equal to $m$. Thus, $g : f \circ g \to f$ is a Cartesian morphism
in $\sC_{/C}$ lifting $g : C'' \to C'$ in $\sC$.
Thus, we have shown that $\sC_{/C} \to \sC$ satisfies the conditions for being a prestack.
\end{proof}

\begin{thm}\label{thm:Yoneda-prest}
Let $\sC$ be a category let $f : C \to C'$ be a morphism in $\sC$.
Then, there is morphism of prestacks $f_* : \sC_{/C} \to \sC_{/C'}$ given by post-composition
with $f$. The assignment:
\[\begin{array}{ccccc}
y &:& \sC &\to& \PSt_{\sC} \\
  &:& C &\mapsto& \sC_{/C} \\
  &:& (C \to[f] C') &\mapsto& (\sC_{/C} \to[f_*] \sC_{/C'})
\end{array}\]
is strict fully faithful functor.
\end{thm}
\begin{proof}
The facts that $f_*$ is a morphism of prestacks and that $y$ is a strict functor are not difficult
to check. To see that it is fully faithful,
we need to show that for two objects $C, C' \in \sC$, the mapping:
\[
y : \sC(C, C') \to \PSt_{\sC}(y(C), y(C'))
\]
is an equivalence of categories. On the left, we have a discrete category,
and hence the only morphisms are identities. Suppose we have a $2$--morphism of prestacks
$\alpha : F \To G : y(C) \to y(C')$. Then, the components of $\alpha$
must lie over identity morphisms in $\sC$. However, if a morphism in $\sC_{/C'}$ lies
over the identity morphism in $\sC$, then it is the identity morphism in $\sC_{/C'}$.
This means $\alpha$ must be the identity $2$--morphism of $f = g$. Thus,
$\PSt_{\sC}(y(C), y(C'))$ is a discrete category as well, and we must show that $y$ is a bijection
onto its set of objects from $\Ob{\sC(C, C')} = \sC(C, C')$.
So, suppose $f, g : C \to C'$ are morphisms in
$\sC$ such that $f_* = g_*$. Then, $f_*(\id_C) = g_*(\id_C) \implies f = g$.
Thus, $y : \sC(C, C') \to \PSt_{/\sC}(y(C), y(C'))$ is injective.
Next, suppose we have a
$1$--morphism of prestacks $F : \sC_{/C} \to \sC_{/C'}$ for objects $C, C' \in \sC$.
Viewing $\id_C : C \to C$ as an object in $\sC_{/C}$, since $F$ is a morphism of prestacks,
$F(\id_C)$ is an object in $\sC_{/C'}$ lying over $C$ in $\sC$, and hence $F(\id_C)$ must be
a morphism $C \to C'$ in $\sC$. Now, suppose we have any other object $f : C'' \to C$
in $\sC_{/C}$. Then, viewing $f$ as a morphism $f \to \id_C$ in $\sC_{/C}$,
$F(f)$ is a morphism $F(f) \to F(\id_C)$ in $\sC_{/C'}$ lying over $f : C'' \to C$:
that is, we must have $F(f) = F(\id_C) \circ f$. That is, $F$ is given on objects by
post-composition by $F(\id_C)$. It is also not difficult to show that $F$ is given
on morphisms by post-composition by $F(\id_C)$. Thus, $F = F(\id_C)_*$.
Thus, $y : \sC(C, C') \to \PSt_{/\sC}(y(C), y(C'))$ is surjective.
\end{proof}

\begin{rmk}
The fully faithful functor of \cref{thm:Yoneda-prest} is the Yoneda embedding for prestacks over
$\sC$, and it corresponds to the usual Yoneda embedding under the Grothendieck construction.
\end{rmk}

\begin{notn}
In light of \cref{thm:Yoneda-prest}, we will abuse notation slightly write
$C$ for $y(C)$.
\end{notn}

\begin{defn}[Representable Morphism]\label{defn:rep-morphism}
Let $f : \mcX \to \mcY$ be a $1$--morphism of prestacks over a category $\sC$.
We say that $f$ is a representable morphism if for each object $C \in \Ob{\sC}$,
and each $1$--morphism $g : C \to \mcY$, the $2$--fibre product of prestacks
$C \times_\mcY \mcX$ is representable.
\end{defn}

\section{Diffeological Stacks}
\label{sec:dflg-st}

This section contains the main results of this work.
In \cref{subsec:dflg}, we review diffeological spaces and show how they capture $C^r$ manifolds for
all $r = 0, 1, \dots, \infty$ simultaneously.
We prove several elementary but useful facts about diffeological spaces and their morphisms
including stability properties of subductions and the fact that they form a Grothendieck topology
on the category $\Dflg$ of diffeological spaces. We also give concrete
characterizations of limits and coproducts of diffeological spaces.
These culminate in \cref{thm:St-Cr-Cinf-equiv}: stacks over $\Mfd^r$
are equivalent to stacks over $\Dflg$ equipped with the Grothendieck topology generated by
subductions.

In \cref{subsec:dflg-st}, we review morphisms of stacks representable by diffeological spaces and
diffeologically presentable, or simply, diffeological stacks.
We prove several important properties of these objects that are then used to prove that
the equivalences of stacks over $\Dflg$, $\Mfd^r$ and $\Cart^r$ from \cref{thm:St-Cr-Cinf-equiv}
preserve diffeological stacks --- see \cref{prop:dflg-st-pullback-iff} --- allowing us to treat
diffeological stacks over all these sites interchangeably.
We also use properties of diffeologically representable morphisms and diffeological stacks
extensively in the remaining three subsections to show diffeological presentability of several
moduli stacks important in complex geometry and non-Abelian Hodge theory.

In \cref{subsec:M-Vect-r}, we construct the moduli stack $\sM_{\Vect, r}(X)$ of $C^r$ vector bundles
on a smooth manifold $X$ in the expected way: we equip $\Gl_n(\bC)$ with the diffeology of $C^r$
functions $U \to \Gl_n(\bC)$ and define
$\sM_{\Vect, r}(X) := \coprod_{n \in \bN} \HHom(X, B\Gl_n(\bC))$.
The results of the previous two subsections along with the fact the $\HHom(X, BG)$ is diffeological
for any diffeological group $G$ (\cref{thm:map-X-BG-dflg}) then imply that $\sM_{\Vect, r}(X)$ is
diffeological --- see \cref{thm:mod-st-Vec-Bun-dflg}.

In \cref{subsec:M-Flt-r}, we construct a bidirectional version of the non-Abelian Hodge filtration
stack in $C^r$ regularity, denoted by $\sM_{\Flt, r}(X)$, and equipped with a map
to $\bC^2$ (\cref{cns:M-Flt-r}). This stack has a natural forgetful map
to $\sM_{\Vect, r}(X) \times \bC^2$ and by computing the fibres of this map explicitly, we show
the stack to be diffeologically presentable.
The fibres of the map $\sM_{\Flt, r}(X) \to \bC^2$ above $(0, 1)$ and $(1, 1)$ are
the moduli stacks $\sM_{\Dol, r}(X)$ and $\sM_{\dR, r}$ parametrizing $C^r$ families of Higgs
bundles and flat bundles respectively over $X$.
The results of the previous three sections immediately yield diffeological presentations for these
stacks. The pullback of $\sM_{\Flt, r}(X)$ along the map $\bC \to \bC^2 : a \mapsto (1, a)$
yields the usual non-Abelian Hodge filtration stack $\sM_{\Hodge, r}(X)$ in $C^r$ regularity,
and similar arguments show that it is diffeological. See \cref{rmk:M-flt-r-concrete}.
We should note that while the families have
$C^r$ regularity, each individual member of such a family is an honest Higgs bundle or flat bundle
and hence a $C^\infty$ object. Hence, we find it more appropriate to refer to these as
$C^{r, \infty}$ families.

In \cref{subsec:M-Diag-r}, we construct simplicial stacks
$S_\bullet\sM_{W, r}$ for $W = \Vect, \Flt, \Dol, \dR, \Hodge$ whose $n$--th level
parametrizes $C^{r, \infty}$ families of diagrams of or quiver representations valued in the
respective objects. We first show that the zeroth and first levels are diffeological which
implies that all levels are diffeological by the same arguments in \cite[\S 4]{ModStQuivBun}
--- see \cref{thm:simplicial-filtration-dflg}, \cref{cor:simplicial-Dol-dR-Hodge-dflg}
and \cref{thm:simplicial-Vect-dflg}.

\subsection{Diffeology}
\label{subsec:dflg}

\begin{defn}[Diffeological Space {\cite[Definition 1.5]{Diff}}]\label{defn:dflg}
A diffeology on a set $X$ is a collection $D \subset \Ob{\Set_{/X}}$ of functions $f : U \to X$,
where $U$ is an open subset of $\bR^n$ for some $n \in \bN$, such that:
\begin{enumerate}
\item \textbf{Covering:} every constant function $U \to[!] \set{*} \to[x] X$ is in $D$.
\item \textbf{Compositionality:} for any smooth map $f : V \to U$ and
$p : U \to X \in D$, $p \circ f : V \to X \in D$.
\item \textbf{Descent:} for any function $p : U \to X$, if there exists an open cover
$U = \bigcup_{i \in I} U_i$ with $p|U_i \in D$ for all $i \in I$, then $p \in D$.
\end{enumerate}
We call the pair $(X, D)$ a diffeological space. We call the elements of $D$, the plots of the
diffeological space $(X, D)$. We will simply write $X$ for the diffeological space when the
diffeology $D$ is clear from context.
A morphism of diffeological spaces $f : (X, D) \to (X', D')$ is a function $f : X \to X'$ that
post composes plots to give plots: that is, for each $p \in D$,
$f \circ p \in D'$.
Diffeological spaces form a category with this notion of morphism, which we denote by $\Dflg$.
\end{defn}

\begin{rmk}
There is an alternative characterization of diffeological spaces as concrete sheaves
on $\Cart$ \cite{BH11}: that is, there is an equivalence of categories
between $\Dflg$ and the category of concrete sheaves on $\Cart$. We will use this characterization
when convenient.
\end{rmk}

\begin{cor}
We have a fully faithful functor $\frgt : \Dflg \to \Set$ that sends a diffeological space
to its underlying set.
\end{cor}

\begin{exm}[$C^r$ Manifolds]
Let $M$ be a $C^r$ manifold, $r = 0, 1, \dots, \infty$ and let $D^r_M$ be the collection of all
$C^r$ maps $U \to M$ where $U$ is an open subset of $\bR^n$ for some $n \in \bN$.
Then, each constant map $U \to M$ is $C^r$ and is, hence, in $D^r_M$.
Post composing $C^r$ maps with smooth maps gives $C^r$ maps.
Finally, $r$--times differentiability is a local property.
Hence, $(M, D^r_M)$ is a diffeological space.
We will call $D^r_M$ the $C^r$ diffeology.
Given a morphism $f : M \to N$ of $C^r$ manifolds and any $C^r$ map $p : U \to M$, $f \circ p$
is again $C^r$, and hence $f$ is a morphism of diffeological spaces $(M, D^r_M) \to (N, D^r_N)$.
\end{exm}

\begin{exm}
The set $\bR^0 = \set{0}$ is a $C^r$ manifold for every $r = 0, 1, 2, \dots, \infty$.
Then, a plot of the $C^r$ diffeology on $\bR^0$ with domain $U \in \Ob{\Cart^r}$
is the unique map $U \to \bR^0$. However, this diffeology is the same regardless of
the choice of $r$.
\end{exm}

\begin{exm}
Given two diffeological spaces $(X, D), (X', D')$, the set $\Dflg((X, D), (X', D'))$
of morphisms of diffeological spaces $(X, D) \to (X', D')$ can be given a diffeology
where a function $p : U \to \Dflg((X, D), (X', D')) \subset \Set(X, X')$
is a plot if the function $U \times X \to X' : (u, x) \mapsto p(u)(x)$ is a morphism of
diffeological spaces. This is called the functional diffeology and makes $\Dflg$ a Cartesian closed
category \cite[\S 1.60]{Diff}.
\end{exm}

\begin{prop}\label{prop:plot-dflg-morphism}
Let $(X, D)$ be a diffeological space and $p : U \to X \in D$ a plot. Then,
$p$ is a morphism of diffeological spaces, where $U$ is equipped with the $C^\infty$ diffeology.
\end{prop}
\begin{proof}
This is precisely the compositionality property of the diffeology $D$.
\end{proof}

\begin{prop}\label{prop:Cr-Dflg-incl}
The mapping
\[\begin{array}{ccccc}
(-, D^r_{-}) &:& \Mfd^r &\to& \Dflg \\
             &:& M &\mapsto& (M, D^r_M) \\
             &:& f &\mapsto& f
\end{array}\]
is a fully faithful functor. In particular, we have a sequence of fully faithful functors:
\[
\Cart^r \hto \Mfd^r \hto \Dflg
\]
\end{prop}
\begin{proof}
$(-, D^r_{-})$ is clearly a faithful functor. It remains to see that it is full.
Let $M, N$ be $C^r$ manifolds equipped with the $C^r$ diffeologies $D_M, D_N$ respectively.
Let $f : (M, D^r_M) \to (N, D^r_N)$ be a morphism of diffeological spaces.
Then, consider any $C^r$ charts $p : U \to M$ and $q : V \to N$. This implies $p \in D^r_M$ and
$q^{-1} \in D^r_V$. Since $f$ is a morphism of diffeological spaces, $f \circ p \in D^r_N$.
Finally, $q^{-1} \circ (f \circ p) \in D^r_V$, which means it is $C^r$. Since this holds for all
charts, $f$ is a $C^r$ map. Thus, our functor is full.
\end{proof}

\begin{warn}
It is tempting to think that the functor $\Cart^r \to \Dflg$ factors through the faithful functor
$\Cart^r \to \Cart^{s}$ for all $0 \leq s < r$, but this is not the case:
if $M$ is a $C^r$ manifold, then it is also a $C^s$ manifold, but the $C^r$
diffeology and $C^s$ diffeology on $M$ are different. More precisely, the identity
map $M \to M$ is a morphism of diffeological spaces $(M, D^r_M) \to (M, D^s_M)$
but it admits no inverse in $\Dflg$. In particular, we have a $2$--cell of categories:
\[\begin{tikzcd}
\Cart^r \ar[rd, "{(-, D^r_{-})}" below left] \ar[rr, ""{below, name=S}] & &
\Cart^s \ar[ld, "{(-, D^s_{-})}" below right] \\
& \Dflg &
\ar[from=S, to=2-2, Rightarrow]
\end{tikzcd}\]
\end{warn}

\begin{prop}\label{prop:lim-dflg}
Given a small diagram $F : I \to \Dflg$, $\lim F$ can be computed as the diffeological
space whose underlying set is $\lim (\frgt \circ F)$ and whose diffeology $D_{\lim F}$
is given by those functions $p : U \to \lim (\frgt \circ F)$, where $U \in \Cart$, that make the
composite $U \to[p] \lim (\frgt \circ F) \to[\pi_i] \frgt(F(i))$ a plot in $F(i)$ for each
$i \in \Ob{i}$, where $\pi_i$ is the canonical projection of the limit in $\Set$.
\end{prop}
\begin{proof}
We will abuse notation slightly and write the underlying set $\lim \frgt \circ F$ also as $\lim F$.
We first check that $D_{\lim F}$ is indeed a diffeology on $\lim F$.
Consider a constant function $U \to \set{*} \to \lim F$.
Then, composing with the projection $\pi_i$ again gives a constant map which
is a plot in $F(i)$, for each $i \in \Ob{I}$. Thus, such a constant map is a plot
of $\lim F$.
Let $q : V \to \lim F$ be a plot in $D_{\lim F}$ and $f : U \to V$, a morphism in
$\Cart$. Then, $\pi_i \circ p$ is a plot of $F(i)$ and hence $\pi_i \circ (p \circ f)$ is also a
plot of $F(i)$ for each $i \in \Ob{I}$. Thus, $p \circ f$ is a plot in $\lim F$.
Next, let $U = \bigcup_{j \in J} U_j$ be an open cover of an object of $\Cart$ with
plots $p_j : U_j \to \lim F$ for each $j \in J$ such that
$p_j|(U_j \cap U_{j'}) = p_j'|(U_j \cap U_{j'})$ for each $j, j' \in J$.
Then, $p^i_{j} := \pi_i \circ p_j : U_j \to F(i)$ is a plot in $F(i)$ for each $i \in I$,
and $p^i_{j}|(U_j \cap U_{j'}) = p^i_{j'}|(U_j \cap U_{j'})$ for each $j, j' \in J$.
Thus, there exists a unique plot $p^i : U \to F(i)$ with $p^i|U_j = p^i_{j}$.
Then, for any morphism $\alpha : i \to i'$ in $I$, we can pass to the $U_j$ to see that
$F(\alpha) \circ p^{i}$ and $p^{i'}$ agree locally:
\[
(F(\alpha) \circ p^i)|U_j
= F(\alpha) \circ (p^i|U_j)
= F(\alpha) \circ p^i_j
= F(\alpha) \circ \pi_i \circ p_j
= \pi_{i'} \circ p_j
= p^{i'}_j
= p^{i'}|U_j
\]
Thus, $F(\alpha) \circ p^i = p^{i'}$ so that the $\set{p^i}_{i \in I}$ make $U$ a cone over
$\frgt \circ F$ in $\Set$. This provides a unique map $p : U \to \lim F$ such that
$\pi_i \circ p = p^i$ for all $i \in I$, but the $p^i$ were plots in $F(i)$, by construction.
Hence, $p$ is also a plot in $\lim F$. Thus, $D_{\lim F}$ is a diffeology on $\lim F$.

We now check that $\lim F$ is indeed a limit in $\Dflg$.
Let $(X, D)$ be diffeological space forming the apex of a cone over $F$ given by morphisms
$f_i : X \to F(i)$.
Then, we get a canonical function of sets $f : X \to \lim \frgt \circ F$ by the universal property
of limits such that $\pi_i \circ f = f_i$.
It remains to check that this function is a morphism of diffeological spaces.
If $p : U \to X$ is a plot then $\pi_i \circ (f \circ p) = f_i \circ p$ is a plot in
$F(i)$ for each $i \in I$. Thus, $f \circ p \in D_{\lim F}$, and we are done.
\end{proof}

\begin{prop}\label{prop:lim-pres-Cr-Mfdr-Dflg}
The fully faith functors $\Cart^r \hto \Dflg$ and $\Mfd^r \hto \Dflg$
preserve all small limits that exist.
\end{prop}
\begin{proof}
Let $F : I \to \Cart^r$ be a small diagram such that $\lim F$ exists in $\Cart^r$.
Then, consider the limits of $F$ in $\Cart^r$ and $\Dflg$, which we will denote by
$X$ and $(Y, D)$ respectively, where $D$ is the limit diffeology from \cref{prop:lim-dflg}.
$X$ has the same underlying set $Y$ as $(Y, D)$. It remains to see that the $C^r$ diffeology on $X$
is the same as the limit diffeology $D$ on $Y$. This follows from the definition of the limit
diffeology and what it means to be a limit in $\Cart^r$. This shows that the functor
$\Cart^r \hto \Dflg$ preserves all small limits that exist in $\Cart^r$. The argument for
$\Mfd^r$ is identical.
\end{proof}

\begin{prop}\label{prop:coprod-dflg}
Let $\set{(X_i, D_i)}_{i \in I}$ be a collection of diffeological spaces.
Let $X = \coprod_{i \in I} X_i$ and $\iota_i : X_i \to \coprod_i X_i$, the canonical inclusion.
Let $D$ be the set of maps $p : U \to X$, such that $U$ is an open subset of $\bR^n, n \in \bN$,
and there exists an open cover $U = \bigcup_j U_j$ with $p|U_j = \iota_{i_j} \circ p^j_0$
for some $i_j \in I$.
Then, $\br{X, D}$ is a diffeological space and is the coproduct in $\Dflg$.
\end{prop}
\begin{proof}
We first verify the axioms of a diffeology.
Let $x \in X = \coprod_{i \in I} X_i$ so that $x \in X_i$ for some specific $i \in I$.
Then, the constant map $c_x : U \to X_i$ valued in $x$ is in $D_i$ so that
$\iota_i \circ c_x$, the constant map $U \to X$ valued in $x$, is in $D$.
Suppose $p : U \to D \in D$ so that there is an open cover $U = \bigcup_{j \in J} U_j$ with
$p|U_{j} = \iota_{i_j} \circ p^j_0$ for some $p^j_0 \in D_{i_j}$ for some $i_j \in I$.
Let $f : V \to U$ be any smooth map with $V \subset \bR^m$ open.
Then, $V = \bigcup_{j \in J} f^{-1}(U_j)$ is an open cover, and
$(p \circ f)|f^{-1}(U_j) = \iota_{i_j} \circ p^j_0 \circ f|f^{-1}(U_j)$.
This shows that $p \circ f \in D$.
Now, suppose $p : U \to X$ is any function
for which there is an open cover $U = \bigcup_j U_j$ with $p|U_j = \iota_{i_j} \circ p^j_0$
for some $p^j_0 \in D_{i_j}$ for some $i_j \in I$, for each $j$. Then, $p$ is in $D$
by definition. Thus, $D$ is a diffeology on $X$.

To show that $(X, D)$ is the coproduct, we first see that each $\iota_{i} : X_i \to X$
is a morphism of diffeological spaces. Let $p : U \to X_i \in D_i$. Then, we can take $U$ as
its own open cover and observe that $p \in D_i$ so that $\iota_i \circ p \in D$.
Thus, $\iota_i$ is a morphism of diffeological spaces.
Now, suppose we have morphisms of diffeological spaces
$f_i : (X_i, D_i) \to (X' D')$ for each $i \in I$. Then, we have a function
$f : X \to X'$ defined on $X_i$ by $f_i$. We wish to show that $f$ is a morphism
of diffeological spaces. Let $p : U \to X \in D$ be a plot. Then, we have an open
cover $U = \bigcup_{j \in J} U_j$ with $p|U_j = \iota_{i_j} \circ p^j_0$ for some
$p^j_0 \in D_{i_j}$ for some $i_j \in I$. Then, for each $j \in J$,
$(f \circ p)|U_j = f_{i_j} \circ p^j_0$ is in $D'$ as $p^j_0$ is in $D_{i_j}$ and $f_{i_j}$ is a
morphism of diffeological spaces. This means $f \circ p \in D'$. Thus, $f$ is a morphism
of diffeological spaces.
Finally, consider a morphism of diffeological spaces $f : (X, D) \to (X', D')$.
Then, $\iota_i$ being a morphism of diffeological spaces implies $f \circ \iota_i$
is a morphism of diffeological spaces. This shows that the bijection:
\[
\Set\br{\coprod_j X_i, X'} \to[\cong] \prod_j \Set\br{X_i, X'}
\]
restricts to a bijection:
\[
\Dflg\br{\coprod_j X_i, X'} \to[\cong] \prod_j \Dflg\br{X_i, X'}
\]
\end{proof}

\begin{defn}[Subduction]\label{defn:subd}
A morphism of diffeological spaces $f : (X, D) \to (X', D')$ is called a subduction
if for every plot $p : U \to X' \in D$ and every $u \in U$, there exists
an open neighbourhood $u \in V \subset U$ and a plot $q : V \to X \in D$ with $p|V = f \circ p$.
The map $q$ is called a local lift of $p$ at $u$ along $f$.
\end{defn}

\begin{rmk}
For any point $x' \in X'$, we have the
morphism $p : \bR^0 \to X'$ sending $0$ to $x'$. This is a morphism of diffeological spaces,
and $f : X \to X'$ being a subduction implies $0$ has an open neighbourhood $U$ in $\bR^0$
such that $p|U$ factors through $f$. However, $U$ must be equal to $\set{*}$. This implies that
$f$ is a surjective function. Thus, \cref{defn:subd} recovers \cite[\S 1.48]{Diff}.
\end{rmk}

\begin{prop}\label{prop:subd-iso}
Every isomorphism of diffeological spaces is a subduction.
\end{prop}
\begin{proof}
Let $f : (X, D) \to (X', D')$ be an isomorphism of diffeological spaces.
Since it is an isomorphism of the underlying sets, it is surjective.
Consider any plot $p : U \to X' \in D'$.
Then, $f^{-1} \circ p$ is a local lift of $p$ with respect to the cover of $U$ with
$U$ itself as the only open set.
\end{proof}

\begin{prop}\label{prop:subd-coprod}
Let $\set{f_i : (Y_i, E_i) \to (X_i, D_i)}_{i \in I}$ be a family of subductions. Then,
\[
\coprod_{i \in I} f_i : \coprod_{i \in I} (Y_i, E_i) \to \coprod_{i \in I} (X_i, D_i)
\]
is a subduction.
\end{prop}
\begin{proof}
$\coprod_{i \in I} f_i$ is surjective since each $f_i$ is.
Let $p : U \to \coprod_{i \in I} X_i$ be a plot such that there exists a cover
$U = \bigcup_{j \in J} U_j$ with $p|U_j = \iota_{i_j} \circ p^j_0$ for some plot
$p^j_0 \in D_{i_j}$ for some $i_j \in I$. Then $p^j_0$ admits local lifts along
$f_{i_j}$ since $f_{i_j}$ is a subduction, by hypothesis. These local lifts composed with
the inclusions $Y_{i_j} \hto \coprod_{i \in I} Y_i$ provide local lifts of $p$
along $\coprod_{i \in I} f_i$, with respect to a refined open cover of $U$.
\end{proof}

\begin{prop}\label{prop:subd-comp}
A composition of subdunctions is a subduction.
\end{prop}
\begin{proof}
Consider a sequence of subductions $(X, D) \to[f] (X', D') \to[g] (X'', D'')$. Then,
$g \circ f$ is surjective since $f$ and $g$ are individually. Next, consider a plot
$p : U \to X'' \in D''$.
We have local lifts of $p$ along $g$ some cover of $U$, each of which again admit local lifts along
$f$ along some cover of each open set in the chosen cover of $U$.
Thus, the local lifts along $f$ of the local lifts along $g$ of $p$ are local lifts of $p$ along
$g \circ f$, with respect to the refined cover of $U$.
\end{proof}

\begin{prop}\label{prop:subd-pullback}
Let $f : (X, D) \to (X', D')$ be a subduction and $g : (X'', D'') \to (X', D')$
be a morphism of diffeological spaces. Then, the map
$g^*f : (X, D) \times_{(X', D')} (X'', D'') \to (X'', D'')$ is a subduction.
\end{prop}
\begin{proof}
The pullback of surjection of sets is a surjection and hence $g^*f$ is surjective
since $f$ is.
Consider a plot $p : U \to X'' \in D''$.
Then $g \circ h$ admits local lifts along $f$ with respect to some cover of $U$, since $f$ is a
subduction. The universal property of the pullback, along with \cref{prop:plot-dflg-morphism}, gives
local lifts of $h$ along $g^*f$.
\end{proof}

\begin{cor}
Given a diffeological space $(X, D)$, let $J_{subd}(X, D)$ be the set of families of morphisms
$\set{c_i : (X_i, D_i) \to (X, D)}_{i \in I}$  such that the canonical map
$c : \coprod_{i \in I} (X_i, D_i) \to (X, D)$ corresponding to the family is a subduction.
Then, $J_{subd}$ is a Grothendieck topology on $\Dflg$.
\end{cor}
\begin{proof}
Apply \cref{prop:subd-iso} for property (i), \cref{prop:subd-coprod} and \cref{prop:subd-comp}
for property (ii) and \cref{prop:subd-pullback} for property (iii) in the definition of a
Grothendieck topology.
\end{proof}

\begin{prop}\label{prop:Cr-man-op-cov-subd}
Let $\set{c_i : M_i \to M}_{i \in I}$ be an open cover of a $C^r$ manifold $M$ by charts.
Then, the function $\coprod_{i \in I} c_i : \coprod_{i \in I} M_i \to M$ is
a subduction.
\end{prop}
\begin{proof}
$\coprod_{i \in I} c_i$ is clear surjective onto $M$.
Let $p : U \to M$ be a plot: that is, a $C^r$ map, where $U$ is some open subset of
$\bR^n, n \in \bN$. Then, $p^{-1}(\im(c_i))$ is an open subset of $U$, and
$p|p^{-1}(\im(c_i))$ factors through $\im(c_i)$, which is homeomorphic to $M_i$ via
$c_i^{-1} : \im(c_i) \to M_i$. Let $\iota_j : M_j \to \coprod_i M_i$ be the canonical inclusion.
Thus, $\iota_i \circ c_i^{-1} \circ p|p^{-1}(\im(c_i))$ is a local
lift of $p$ along $\coprod_{i \in I} c_i$. Furthermore, the $p^{-1}(\im(c_i))$ as $i$ ranges over
$I$ cover $U$.
\end{proof}

\begin{prop}\label{prop:Cr-Dflg-cocont}
The functors $\Cart^r \hto \Dflg$ and $\Mfd^r \hto \Dflg$ are cocontinuous.
\end{prop}
\begin{proof}
Suppose we have a family of morphisms of diffeological spaces
$\set{d_i : D_i \to M}$ where $M$ is a $C^r$ manifold and the canonical morphism
$d : \coprod_{i \in I} D_i \to M$ associated with this family is a subduction.
Consider any open cover $\set{c_q : M_q \to M}_{q \in Q}$ by charts, and consider an open
cover $\set{c_{q, r} : M_{q, r} \to M_q}_{r \in R_q}$ and local
lifts $c_{q, r} : M_{q, r} \to \coprod_{i} D_i$ of $c_q$ along $d$.
By passing to a further refinement if necessary,
we may assume $c_{q, r}$ factors through $D_{\gamma(q, r)}$
for some $\gamma(q, r) \in I$, and we are done.
\end{proof}

\begin{prop}\label{prop:Dflg-cov-by-Cr}
Let $(X, D)$ be a diffeological space, then there exists a covering family
\[
\set{c_i : (U_i, D^r_{U_i}) \to (X, D)}_{i \in I}
\]
where each $U_i$ is an open subset of some Euclidean space equipped with the $C^r$ diffeology.
\end{prop}
\begin{proof}
Consider the family $\set{p : \dom(p) \to X}_{p \in D}$: this is a family of morphisms
in $\Dflg$ by \cref{prop:plot-dflg-morphism} whose domains are open subsets of
Euclidean spaces.
We will show that $\coprod_{p \in D} p : \coprod_{p \in D} \dom(p) \to X$ is a subduction.
Let $q : U \to X \in D$ be a plot. Then, we have a canonical inclusion
$\iota : U = \dom(q) \to \coprod_{p \in D} \dom(p)$. This is a local lift
of $q$ along $\coprod_{p \in D} p$ with the open cover of $U$ where $U$ is only open set
in the cover. Then, we observe that $\dom(p) \in \Cart^\infty \subset \Cart^r \subset \Mfd^r$
for every $r = 0, 1, 2, \dots, \infty$. Thus, we are done.
\end{proof}

\begin{thm}\label{thm:St-Cr-Cinf-equiv}
The mappings $(-, D^r_{-}) : \Cart^r \to \Dflg$, $(-, D^r_{-}) : \Mfd^r \to \Dflg$
induces, by pullback, an equivalence of $2$--categories:
\[
\St_{\Cart^r} \ot[\simeq] \St_{\Dflg, J_{subd}} \to[\simeq] \St_{\Mfd^r}
\]
for all $r = 0, 1, 2, \dots, \infty$.
In particular, for all such $r$, we have equivalences of $2$--categories:
\[
\St_{\Cart^r} \simeq \St_{\Mfd^r} \simeq \St_{\Dflg, J_{subd}} \simeq \St_{\Cart}
\]
\end{thm}
\begin{proof}
\cref{prop:Cr-Dflg-incl}, \cref{prop:lim-pres-Cr-Mfdr-Dflg}, \cref{prop:Cr-man-op-cov-subd}
\cref{prop:Cr-Dflg-cocont} and \cref{prop:Dflg-cov-by-Cr} show that the hypotheses
of \cref{thm:comp-thm} hold when we take $u$ to be the functors
$(-, D^r_{-}) : \Cart^r \to \Dflg$ and $(-, D^r_{-}) : \Mfd^r \to \Dflg$.
\end{proof}

\begin{rmk}
This is very similar to the argument in \cite{CarchediComparison} showing stacks over
$\Dflg$ are equivalent to stacks over $\Mfd$, but we have extended this to the finite regularity
cases.
\end{rmk}

\begin{rmk}
\cref{thm:St-Cr-Cinf-equiv} can be interpreted as the statement: stacky $C^r$ geometry
can be reduced to stacky smooth geometry. Moreover, this is mediated by diffeological spaces,
and requires no analysis explicitly. Note that for $r = 0$, we are not saying that stacks over
$\Top$ can be reduced to stacks over $\Cart$, but that stacks over $\Mfd^0$ can be.
\end{rmk}

\subsection{Diffeological Stacks}
\label{subsec:dflg-st}

\begin{defn}[Prestack Represented by a Diffeological Space]
Let $X$ be a diffeological space. We have the representable prestack $\Dflg_{/X} \to \Dflg$.
We can pull this prestack back along the fully faithful functors $\Mfd^r \to \Dflg$,
$\Cart^r \to \Dflg$ for each $r = 0, 1, 2, \dots, \infty$.
We will call $\Dflg_{/X}$ and its pullbacks to $\Cart^r, \Mfd^r$ the prestacks represented by $X$
over the respective sites. Any prestack equivalent to such a prestack will be called diffeologically
representable.
We will simply write $X$ for all of these prestacks, by slight abuse of notation.
\end{defn}

\begin{cor}\label{cor:dflg-rep-prest-pullback}
Diffelogically representable prestacks are preserved by pullback along
the fully faithful functors $\Cart^r \to \Mfd^r, \Mfd^r \to \Dflg, \Cart^r \to \Dflg$.
\end{cor}

\begin{prop}\label{prop:dflg-sp-st}
Given a diffeological space $Y$, the prestack $\Dflg_{/Y} \to \Dflg$ is a stack.
In particular, the prestack represented by $Y$ over $\Mfd^r, \Cart^r$ is a stack for all
$r = 0, 1, 2, \dots, \infty$. Thus,
The Yoneda embeddings $\Dflg \to \PSt_{\Dflg}, \Mfd^r \to \PSt_{\Mfd^r}, \Cart^r \to \PSt_{\Cart^r}$
factor through the inclusions $\St_{\Dflg} \hto \PSt_{\Dflg}, \St_{\Mfd^r} \hto \PSt_{\Mfd^r},
\St_{\Cart^r} \hto \PSt_{\Cart^r}$.
\end{prop}
\begin{proof}
Suppose we have a family of morphisms of diffeological spaces $\set{\pi_i : X_i \to X}_{i \in I}$
such that the corresponding morphism $\pi : \coprod_{i \in I} X_i \to X$ is a subduction.
and morphisms of diffeological spaces $f_i : X_i \to Y$ such that the following diagram commutes
for each $i, j \in I$:
\[\begin{tikzcd}
X_i \times_X X_j \ar[r] \ar[d] & X_i \ar[d, "f_i"] \\
X_j \ar[r, "f_j" below] & Y
\end{tikzcd}\]
Then, consider any $x \in X$ so that $x = \pi_i(x_i)$ for some $x_i \in X_i$ and $i \in I$,
since $\pi$ is a surjection. At the same time if $x = \pi_j(x_j)$ for some $j \in I$,
then $f_i(x_i) = f_j(x_j)$, by the commutativity of the above diagram.
We set $f(x) = f_i(x_i)$ giving a function $f : X \to Y$. Now, suppose we have a plot
$p : U \to X$. Since $\pi$ is a subduction, we have an open cover $U = \bigcup_{k \in K} U_k$
such that $p|U_k$ factors through $\coprod_{i \in I} X_i$. By \cref{prop:coprod-dflg},
refining covers if needed, we may assume $p|U_k$ factors as through
$U_k \to[p_k] X_{\gamma(k)} \to[\pi_{\gamma(k)}] X$ for some plot
$p_k : U_k \to X_{\gamma(k)}$, $\gamma(k) \in I$.
Thus, $f \circ p|U_k = f_{\gamma(k)} \circ p_k$. The descent property of the diffeology
of $Y$ now implies that $f \circ p$ is plot of $Y$. That is, $f$ is morphism of diffeological
spaces. This shows that $\Dflg(-, Y)$ is a sheaf. Since $\Dflg_{/Y}$ is the Grothendieck
construction of a sheaf $\Dflg(-, Y)$, it is a stack.
\end{proof}

\begin{prop}\label{prop:pullback-dflg-rep-implies-dflg-rep-ob}
If a stack pulls back along the fully faithful functors
$\Cart^r \to \Mfd^r, \Mfd^r \to \Dflg, \Cart^r \to \Dflg$ to a diffeologically representable one,
then it is representable by diffeological spaces.
\end{prop}
\begin{proof}
Suppose we have a prestack $p : \sP \to \Dflg$ such that its pullback $\sP'$ to $\Mfd^r$ is
representable by a diffeological space $X$. Then, $\sP'$ is a stack on $\Mfd^r$ by
\cref{prop:dflg-sp-st}. Choose a $2$--inverse $i : \St_{\Mfd^r} \to \St_{\Dflg}$
of the pullback functor $\St_{\Dflg} \to \St_{\Mfd^r}$ by \cref{thm:St-Cr-Cinf-equiv}. Then,
$\Dflg_{/X} \simeq i(\Mfd^r_{/X}) \simeq i(\sP') \simeq \sP$. The proofs for the other cases
are similar.
\end{proof}

\begin{warn}
The argument of \cref{prop:pullback-dflg-rep-implies-dflg-rep-ob} does not work
for general prestacks as we do not have an equivalence between
$\PSt_{\Dflg}$ and $\PSt_{\Mfd^r}$.
\end{warn}

\begin{prop}\label{prop:subfunc-dflg}
Let $X$ be a diffeological space and $\iota : F \hto \Dflg_{/X}$, a fully faithful morphism of
stacks. Then, $F$ is representable by a diffeological space.
In particular, $F|\Mfd^r, F|\Cart^r$ are representable by diffeological spaces for all
$r = 0, 1, 2, \dots, \infty$
\end{prop}
\begin{proof}
The inverse Grothendieck construction $\st(F)$ of $F$ is a sheaf since $F$ is a stack
by hypothesis. Furthermore, it is a subfunctor
of the functor $y(X) : \Dflg^\op \to \Set : Y \mapsto \Dflg(Y, X)$ determined by $X$.
This means $\st(F)|\Mfd^\op$ is a subfunctor of $y(X)|\Mfd^\op$. However, $X$ being a diffeological
space means it is a concrete sheaf which implies $\st(F)$ is also a concrete sheaf and hence
is representable by a diffeological space. This implies $F$ is representable by a diffeological
space, by applying the Grothendieck construction. The last claim follows from
\cref{prop:pullback-dflg-rep-implies-dflg-rep-ob}.
\end{proof}

\begin{defn}[Diffeologically Representable Morphisms]
Let $\sC = \Cart^r, \Mfd^r$ or $\Dflg$, $r = 0, 1, 2, \dots, \infty$.
A $1$--morphism of prestacks $f : \sP \to \sQ$ over $\sC$ is said to be diffeologically
representable or representable by morphisms of diffeological spaces if for each $C \in \Ob{\sC}$ and
each $1$--morphism of prestacks $g : C \to \sQ$, the $2$--fibre product of prestacks
$\sP \times_\sQ C$ is representable by a diffeological space.
Note that since the Yoneda embedding $\Dflg \to[y] \St_{\Dflg}$ is fully faithful
(\cref{thm:Yoneda-prest})
and we have an equivalence $\St_{\Dflg} \to[\simeq] \St_{\sC}$ given by pullback
(\cref{thm:St-Cr-Cinf-equiv}),
the morphism $\sP \times_\sQ C \to C$ is (equivalent in the arrow $2$--category of $\St_{\Cart^r}$)
to a morphism of diffeological spaces.
If $M$ is a class of morphisms in $\Dflg$, then $g$ is said to be
in $M$ if for each $f$ as above, the morphism $\sP \times_\sQ C \to C$ is (equivalent to) a
morphism in the image of $M$ under the Yoneda embedding.
\end{defn}

\begin{exm}
A $1$--morphism of prestacks $p : \sP \to \sQ$ over $\sC = \Cart^r, \Mfd^r$ or $\Dflg$
is a subduction if it is diffeologically representable and its pullback along every $1$--morphism
of prestacks $C \to \sQ$, where $C$ is representable by a diffeological space, is a subduction.
\end{exm}

\begin{prop}\label{prop:subd-st-comp}
A composition of diffeologically representable morphisms is diffeologically representable.
The same holds for subductions of stacks.
\end{prop}
\begin{proof}
The first statement follows from pasting $2$--fibre product diagrams.
The second statement follows from this same observation in combination with \cref{prop:subd-comp}.
\end{proof}

\begin{prop}\label{prop:subd-st-pullback}
A pullback of a diffeologically representable morphism is diffeologically representable.
The holds for subductions of stacks.
\end{prop}
\begin{proof}
This again follows from pasting $2$--fibre product diagrams and
\cref{prop:subd-pullback}.
\end{proof}

\begin{prop}\label{prop:subd-st-coprod}
Countable coproducts of diffeologically representable morphisms are diffeologically representable.
The same holds for subductions.
\end{prop}
\begin{proof}
The first statement follows from the fact that $2$--fibre products commute with countable
coproducts. The second statement follows from the same observation in combination with
\cref{prop:subd-st-pullback}.
\end{proof}

\begin{prop}\label{prop:diff-rep-morphism-pullback-iff}
A $1$--morphism in $\St_{\Dflg}$ is diffeologically representable if and only if its
its image along the pullback functor $\St_{\Dflg} \to \St_{\Mfd^r}$ is one.
The same holds for $\St_{\Dflg} \to \St_{\Cart^r}$ and
$\St_{\Mfd^r} \to \St_{\Cart^r}$.
\end{prop}
\begin{proof}
Suppose $f : \mcX \to \mcY$ is a diffeologically representable morphism of stacks over $\Dflg$.
Then, consider the pullback $f' : \mcX' \to \mcY'$ of $f$ along $\Mfd^r \to \Dflg$ and
any $1$--morphism $g : M \to \mcY'$ in $\St_{\Mfd^r}$, where $M \in \Ob{\Mfd^r}$.
Choose a $2$--inverse $i : \St_{\Mfd^r} \to \St_{\Dflg}$ of the pullback functor by
\cref{thm:St-Cr-Cinf-equiv}.
Then, we have a morphism $i(g) : i(M) \to i(X') \to[\simeq] \mcX$. Since
$i(M) \simeq \Dflg_{/M}$, this gives a morphism $M \to \mcX$.
Then, $\mcX' \times_{\mcY'} M$ is equivalent to the image of
$\mcX \times_{\mcY} M$ under the pullback functor $\St_{\Dflg} \to \St_{\Mfd^r}$, since
it is an equivalence of $2$--categories. However, $\mcX \times_\mcY M$ is diffelogically
representable by hypothesis, and hence so is $\mcX' \times_{\mcY'} M$ by
\cref{cor:dflg-rep-prest-pullback}.

Now, suppose $f' : \mcX' \to \mcY'$ is diffeologically representable.
Then, $\mcX' \times_{\mcY'} M$ is diffeologically representable. It is also the pullback
of $\mcX \times_{\mcY} M$, which must, hence, be diffeologically representable by
\cref{prop:pullback-dflg-rep-implies-dflg-rep-ob}.

The other cases are similar.
\end{proof}

\begin{prop}\label{prop:subd-pullback-st-iff}
A $1$--morphism of stacks is a subduction if and only if its image along the pullback functor
$\St_{\Dflg} \to \St_{\Mfd^r}$ is a subduction. The same holds for the pullback functors
$\St_{\Mfd^r} \to \St_{\Cart^r}$ and $\St_{\Dflg} \to \St_{\Cart^r}$.
\end{prop}
\begin{proof}
Let $f : \mcX \to \mcY$ be a subduction in $\St_{\Dflg}$. Consider the pullback
$f' : \mcX' \to \mcY'$ to $\St_{\Mfd^r}$, a $1$--morphism $g : M \to \mcY'$,
and a morphism of diffelogical spaces $h : X \to M$.
Choose a $2$--inverse $i$ of the pullback functor $\St_{\Dflg} \to \St_{\Mfd^r}$.
Then, $i(\mcX' \times_{\mcY'} M) \simeq \mcX \times_\mcY M$ in $\St_{\Dflg}$ by
\cref{thm:St-Cr-Cinf-equiv}. $\mcX \times_\mcY M \to M$ is a subduction by hypothesis
and hence the map $X \simeq i(X) \to[i(h)] i(M) \simeq M$ admits local lifts
along $\mcX \times_\mcY M \to M$. The pullbacks of these local lifts to $\St_{\Mfd^r}$
provide local lifts of $h$ along $\mcX' \times_{\mcY'} M \to M$, showing that it is a subduction.

Suppose, now, that $f' : \mcX' \to \mcY'$ is a subduction. Consider a $1$--morphism
$g : M \to \mcY$ in $\St_{\Dflg}$ and morphism of diffeological spaces $h : X \to M$.
Then, we get $1$--morphisms $X \to[h] M \to[g'] \mcY'$ in $\St_{\Mfd^r}$ by pullback.
We get local lifts of $h$ along $\mcX' \times_{\mcY'} M \to M$, since the hypothesis implies it
is a subduction. The images of these local lifts under $i$ yield local lifts of
$h$ along $\mcX \times_{\mcY} M \to M$ because it is equivalent in the arrow $2$--category
of $\St_{\Dflg}$ to $i(\mcX' \times_{\mcY'} M \to M)$.

The other cases are similar.
\end{proof}

\begin{defn}[Diffeological Stack]\label{defn:dflg-st}
Let $\sC = \Cart^r, \Mfd^r$ or $\Dflg$, $r = 0, 1, 2, \dots, \infty$.
A stack $p : \mcX \to \sC$ is called diffeological or diffeologically presentable or presentable by
diffeological groupoids if there exists
a diffeological space $A$ and a subduction of stacks $a : A \to \mcX$.
The morphism $a$ is called a diffeological atlas, in this case.
\end{defn}

\begin{prop}\label{prop:subd-st-epi}
A subduction of stacks $f : \mcX \to \mcY$ over $\Dflg$ is an epimorphism of stacks.
The same holds for $\Mfd^r$ or $\Cart^r$ in place of $\Dflg$.
\end{prop}
\begin{proof}
Consider a diffeological space $S$ and an object $y$ of $\mcY$ over $S$.
By the Yoneda lemma for prestacks \cite[Lemma 2.4.21]{AlperModuli}, we have
a $1$--morphism of (pre)stacks $\ul{y} : \Dflg_{/S} \to \mcY$ corresponding to $y$.
Using \cref{prop:Dflg-cov-by-Cr}, choose a covering family $\set{U_i \to S}_{i \in I}$
where each $U_i$ is an object of $\Cart^r$.
Passing to a finer cover if necessary and using the hypothesis that $f$ is a subduction, we may
assume $U_i \to S$ factors through the map $\mcX \times_\mcY S \to S$. This, in turn, gives
a map $\ul{x_i} : U_i \to \mcX$ such that $U_i \to U \to[\ul{y}] \mcY$ factors up to
$2$--isomorphism in $\St_{\Dflg}$ as $U_i \to \mcX \to[f] \mcY$.
Under the $2$--Yoneda lemma, this means precisely that the objects $x_i$ corresponding
to the maps $\ul{y_i}$ satisfy $f(x_i) \cong y|U_i$ in $\mcY$.
This is the characterization of epimorphisms of stacks given in \cite[Proposition 2.13(ii)]{BX11}.
The other cases are similar since every $C^r$ manifold admits a cover by objects of $\Cart^r$
by definition.
\end{proof}

\begin{rmk}
\cref{prop:subd-st-epi} shows that \cref{defn:dflg-st} recovers
\cite[Definition 2.9 and following paragraphs]{RV18}.
\end{rmk}

\begin{prop}\label{prop:dflg-st-pullback-iff}
A stack over $\Dflg$ is diffeological if and only if its pullback along the fully faithful functor
$\Mfd^r \to \Dflg$ is diffeological. The same holds for the fully faithful functors
$\Cart^r \to \Mfd^r$ and $\Cart^r \to \Dflg$.
\end{prop}
\begin{proof}
Apply \cref{cor:dflg-rep-prest-pullback}, \cref{prop:pullback-dflg-rep-implies-dflg-rep-ob},
\cref{prop:diff-rep-morphism-pullback-iff} and \cref{prop:subd-pullback-st-iff}.
\end{proof}

\begin{cor}\label{cor:dflg-rep-morphism-dflg-target}
If $\mcX \to \mcY$ is a diffeologically representable morphism and $\mcY$ is diffeologically
representable, then so is $\mcX$.
\end{cor}
\begin{proof}
Choose a diffeological atlas $A \to \mcY$, observe that $A \times_\mcY \mcX$ is representable
by diffeological spaces, and then apply \cref{prop:subd-st-pullback}
to see that the map $A \times_\mcY \mcX \to \mcX$ is a diffeological atlas.
\end{proof}

\begin{thm}\label{thm:fin-lim-dflg-st}
Fibre products of diffeological stacks are diffeological: that is, for a diagram
$\mcX \to \mcY \ot \mcZ$ in $\St_\Dflg$ where $\mcX, \mcY$ and $\mcZ$ are
diffeological, $\mcX \times_\mcY \mcZ$ is diffeological.
The same results holds with $\Cart^r, \Mfd^r, r = 0, 1, \dots, \infty$ in place of
$\Dflg$.
In particular, finite limits of diffeological stacks are diffeological.
\end{thm}
\begin{proof}
Choose diffeological atlases $U \to \mcX, V \to \mcZ$ and $A \to \mcY$ and consider the following
pasting of Cartesian diagrams:
\[\begin{tikzcd}
U \times_\mcY V \ar[r] \ar[d] &
U \times_\mcX (\mcX \times_\mcY \mcZ) \ar[r] \ar[d] &
U \ar[d] \\
(\mcX \times_\mcY \mcZ) \times_\mcZ V \ar[r] \ar[d] &
\mcX \times_\mcY \mcZ \ar[r] \ar[d] &
\mcX \ar[d] \\
V \ar[r] & \mcZ \ar[r] & \mcY
\end{tikzcd}\]
Notice that all arrows in the above diagram outside the bottom right square are subductions
by \cref{prop:subd-st-pullback}. This implies that the map
$U \times_\mcY V \to \mcX \times_\mcY \mcZ$ is a subduction. However, we are not
done yet as we cannot deduce that $U \times_\mcY V$ is representable by diffeological space.
Then, consider the following pasting of Cartesian diagrams:
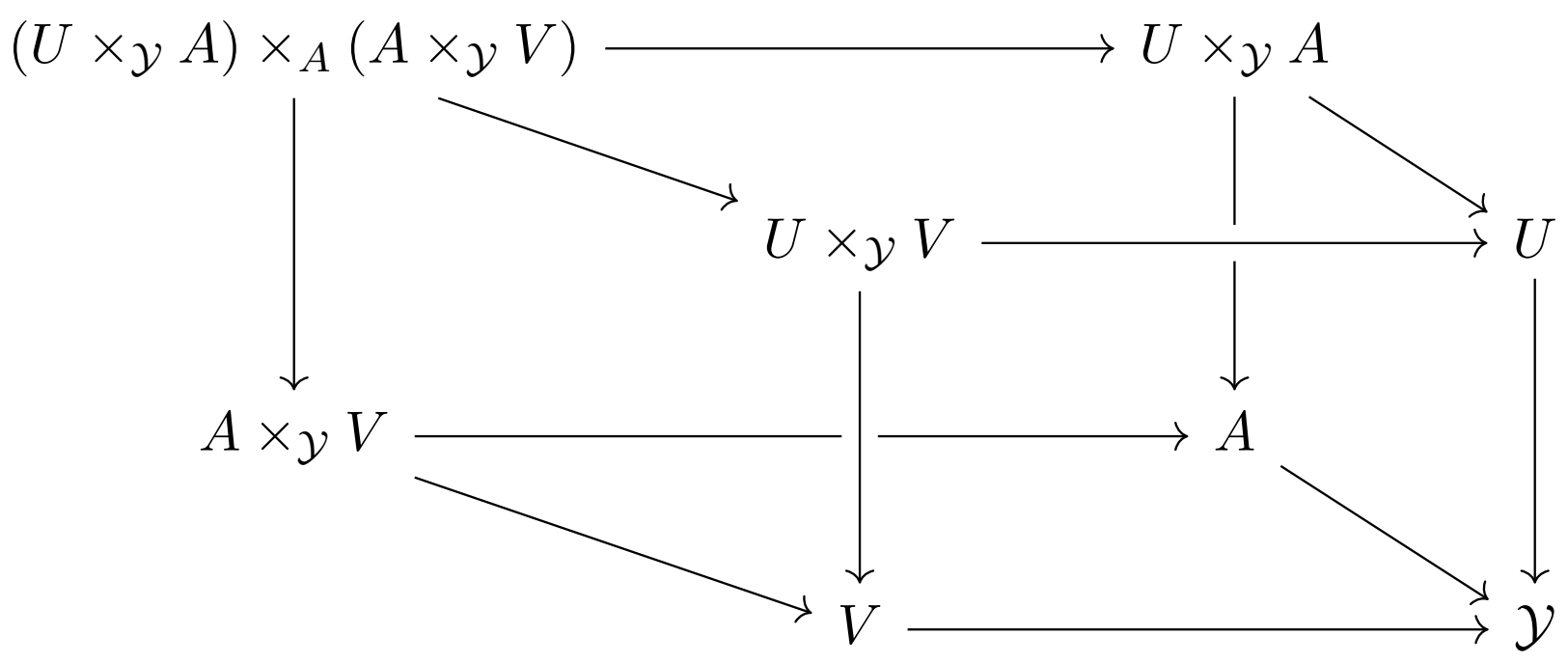
\[\begin{tikzcd}
(U \times_\mcY A) \times_A (A \times_\mcY V) \ar[rr] \ar[rd] \ar[dd] & &
U \times_\mcY A \ar[rd] \ar[dd] & \\
& U \times_\mcY V \ar[rr, crossing over] & & U \ar[dd] \\
A \times_\mcY V \ar[rr] \ar[rd] & & A \ar[rd] & \\
& V \ar[rr] \ar[from=uu, crossing over] & & \mcY
\end{tikzcd}\]
Note that since $A \to \mcY$ is an atlas while $U$ and $V$ are representable by diffeological
spaces, so are $U \times_\mcY A$ and $A \times_\mcY V$.
This implies that $(U \times_\mcY A) \times_A (A \times_\mcY V)$ is representable by
a diffeological space.
Again, since $A \to \mcY$ is an atlas, $U \times_\mcY A \to U$ is a subduction and this, in turn,
implies that $(U \times_\mcY A) \times_A (A \times_\mcY V) \to U \times_\mcY V$ is a subduction, by
\cref{prop:subd-st-pullback}. The composite
\[
(U \times_\mcY A) \times_A (A \times_\mcY V) \to U \times_\mcY V \to \mcX \times_\mcY \mcZ
\]
is, therefore, a diffeological atlas.

Next, observe that any finite limit can be computed using
fibre products and a terminal object.
The terminal object is $\Mfd$ and is representable by the singleton diffeological space
$\pt = \bR^0$. The final statement now follows.
\end{proof}

\subsection{Moduli Stack of $C^r$ Vector Bundles}
\label{subsec:M-Vect-r}

\begin{defn}[Diffeological Group]
A diffeological group is a group object in $\Dflg$.
\end{defn}

\begin{cor}
The fully faithful functors $\Mfd^r \to \Dflg$ preserve group objects.
\end{cor}
\begin{proof}
Any functor that preserves products preserves group objects. The result thus follows from
\cref{prop:lim-pres-Cr-Mfdr-Dflg}.
\end{proof}

\begin{exm}
Let $G$ be any Lie group. Then $(G, D^\infty_G)$ is a diffeological group.
\end{exm}

\begin{exm}
Let $G$ be any Lie group so that it is also a group object in $\Mfd^r$ for all $0 \leq r < \infty$
as well. Then, $(G, D^r_G)$ is a diffeological group.
\end{exm}

\begin{defn}[Classifying Stack of a Diffeological Group]
Given any diffeological group $G$, we have a unique action $G \cong G \times \bR^0 \to \bR^0$ of
$G$ on $\bR^0$ equipped with its unique diffeology.
We will denote the $2$--colimit of the diagram $G \rightrightarrows \bR^0$ in $\St_{\Dflg}$ or the
quotient stack $[\bR^0/G]$ as $BG$.
\end{defn}

\begin{cor}\label{cor:BG-Cr-Dflg}
Let $G$ be a group object in $\Mfd^r$. Then, the $2$--colimits of the diagram
$G \rightrightarrows \bR^0$ in $\St_{\Mfd^r}$ and $\St_{\Cart^r}$ are equivalent to the pullbacks of
$BG$ under the fully faithful functors $\Mfd^r \to \Dflg, \Cart^r \to \Dflg$ respectively.
Furthermore, the $2$--colimit in $\St_{\Cart^r}$ is also equivalent to the pullback of the
$2$--colimit in $\St_{\Mfd^r}$.
\end{cor}
\begin{proof}
$2$--equivalences of $2$--categories preserve $2$--colimits, and the result follows from
\cref{thm:St-Cr-Cinf-equiv}.
\end{proof}

\begin{notn}
Given \cref{cor:BG-Cr-Dflg}, we will write $BG$ for the $2$--colimits in
$\St_{\Dflg}, \St_{\Mfd^r}, \St_{\Cart^r}$, where there is no chance of confusion.
\end{notn}

\begin{thm}\label{thm:BG-dflg}
Given any diffeological group $G$, $BG$ is a diffeological stack.
\end{thm}
\begin{proof}
The quotient stack of a groupoid object in $\Dflg$ is a diffeological stack by
\cite[Theorem 2.10]{RV18}.
\end{proof}

\begin{defn}[Mapping Prestack]
Let $\sC$ be a category and $p : \sP \to \sC, q : \sQ \to \sC$, two prestacks on $\sC$.
We define a category $\HHom(\sP, \sQ)$ as follows:
\begin{itemize}
\item Objects are tuples $(C, F)$, where:
  \begin{itemize}
  \item $C \in \Ob{\sC}$
  \item $F$ is a $1$--morphism of prestacks $\sP \times_\sC y(C) \to \sQ$
  \end{itemize}
\item Morphisms $(C, F) \to (C', F')$ are tuples $(f, \alpha)$ where:
  \begin{itemize}
  \item $f : C \to C'$ is a morphism in $\sC$
  \item $\alpha$ is a $2$--morphism of prestacks:
  \[\begin{tikzcd}
  \sP \times_\sC y(C)
    \ar[rd, "F" below left]
    \ar[rr, "\id_\sP \times_\sC y(f)" above, ""{name=S, below}]  & &
  \sP \times_\sC y(C') \ar[ld, "F'" below right] \\
  & \sQ &
  \ar[from=S, to=2-2, Rightarrow, "\alpha" left]
  \end{tikzcd}\]
  \end{itemize}
\end{itemize}
\end{defn}

\begin{prop}\label{prop:mapping-prest-st}
Let $\sP$ be a prestack and $\mcX$, a stack over a site $(\sC, J)$. Then,
the mapping prestack $\HHom(\sP, \mcX)$ is a stack.
\end{prop}
\begin{proof}
Let $C$ be an object of $\sC$; $\set{c_i : C_i \to C}_{i \in I}$, a covering family and
$(C_i, F_i)$. We will first show collation of objects in $\HHom(\sP, \mcX)$.
Consider objects of $\HHom(\sP, \mcX)$ over $C_i$ for each $i \in I$.
Suppose, we also have isomorphisms
$(\id_{C_i \times_C C_j}, \alpha_{ij}) : (C_i, F_i)|C_i \times_C C_j
\to (C_j, F_j)|C_i \times_C C_j$ satisfying the cocycle condition.
By \cref{prop:strictify-2-fibre-prod}, we may assume objects of $\sP \times_{\sC} y(C)$
are of the form $(P, a)$ where $P \in \Ob{sP}$, $a : C' \to C \in \Mor{\sC}$
and $P$ lies over $C'$. Then, we have the objects $P|C' \times_C C_i$ over $C' \times_C C_i$
and the pullbacks $a_i$ of $a$ along the morphisms $c_i$,
so that $(P|C' \times_C C_i, a_i)$ is an object of $\sP \times_\sC y(C' \times_C C_i)$.
Applying the composite functor
$\sP \times_\sC y(C' \times_C C_i) \to \sP \times_\sC y(C_i) \to[F_i] \mcX$
to this object provides an object of $X_i$ of $\mcX$. Due to the natural isomorphisms
$\alpha_{ij}$, the $X_i$ form a descent datum for an object $F(P, a)$ in $\mcX$.
Similarly, a morphism $(P_1, a_1) \to (P_2, a_2)$ in $\sP \times_\sC y(C)$ is of the form
$(p, f)$ where $f : C'_1 = \dom(a_1) \to C'_2 = \dom(a_2)$ is a morphism in $y(C)$
and $p : P_1 \to P_2$ is a morphism in $\sP$ over $f$. Then,
$p|C_i : P_1|C'_1 \times_C C_i \to P_2|C_2' \times_C C_i$ is a morphism over
$f \times_C C_i : C'_1 \times_C C_i \to C'_2 \times_C C_i$ so that $(p|C_i, f \times_C C_i)$
is a morphism in $\sP \times_{y(C)} y(C_i)$, and we can apply the functor
$\sP \times_{y(C)} y(C_i) \to[F_i] \mcX$ to get a morphism
$F_i(P_1|C_1' \times_C C_i, a_1|C_i)
\to[F_i(p|C_i, f \times_C C_i)] F_i(P_2|C_2' \times_C C_i, a_2|C_i)$ in $\mcX$.
From the collation of objects, we have a morphism
$F_i(P_2|C_2' \times_C C_i, a_2 \times_C C_i) \to F(P_2, a_2)$.
Thus, we get a morphism in $\mcX$ of the form:
\[
F_i(P_1|C_1' \times_C C_i, a_1|C_i)
    \to[F_i(p|C_i, f \times_C C_i)] F_i(P_2|C_2' \times_C C_i, a_2|C_i)
    \to F(P_2, a_2)
\]
The naturality and cocycle condition satisfied by the $\alpha_{ij}$ can then be used to
show that these morphisms form a descent datum for a morphism $F(P_1, a_1) \to F(P_2, a_2)$ in
$\mcX$. Furthermore, the uniqueness of collation of morphisms shows that the assignment is
$F(p, f)$ functorial.

Next, we will show collation of morphisms in $\HHom(\sP, \mcX)$.
Consider two objects $(C, F)$ and $(C', F')$ of $\HHom(\sP, \mcX)$, a morphism $f : C \to C'$ and
morphisms $(f \circ c_i, \alpha_i) : (C, F)|C_i \to (C', F')$ forming a descent datum for a morphism
in $\HHom(\sP, \mcX)$. Consider an object $(P, a)$ in $\sP \times_{\sC} y(C)$ so that
$P$ is an object in $\sP$ over $C'' = \dom(a : C'' \to C)$ in $\sC$. Then, we have the object
$(P|C'' \times_C C_i, a|C_i)$ in $\sP \times_\sC y(C_i)$, and the morphism
$\alpha_i : F'(P|C'' \times_C C_i, f \circ c_i \circ a|C_i)
\to F(P|C'' \times_C C_i, c_i \circ a|C_i)$
in $\mcX$. This yields a morphism:
\[
F'(P|C'' \times_C C_i, f \circ c_i \circ a|C_i)
    \to[\alpha_i] F(P|C'' \times_C C_i, c_i \circ a|C_i)
    \to F(P, a)
\]
These form a descent datum for a morphism $\alpha : F'(P, f \circ a) \to F(P, a)$ in $\mcX$
by the hypothesis that the $(f \circ c_i, \alpha_i)$ form a descent datum for a morphism
in $\HHom(\sP, \mcX)$. The naturality of the $\alpha$ can be shown from that of the
$\alpha_i$. Finally, $(f, \alpha)$ can be shown to be the collation of the
$(f \circ c_i, \alpha_i)$.
\end{proof}

\begin{prop}
Let $(\sC, J)$ and $(\sD, K)$ be sites, and $u : \sC \to \sD$.
Then, for every object $\sP, \sQ$ in $\PSt_{\sD}$, there is a morphism of prestacks:
\[
u^*\HHom(\sP, \sQ) \to \HHom(u^*\sP, u^*\sQ)
\]
which is an equivalence when $\sP$ and $\sQ$ are stacks and $u$ satisfies the hypotheses of
\cref{thm:comp-thm}.
\end{prop}
\begin{proof}
By \cref{prop:strictify-2-fibre-prod}, an object of $u^*\HHom(\sP, \sQ)$ is of the form
$(C, f)$, where $f : \sP \times_\sD y(u(C)) \to \sQ$ is a $1$--morphism of prestacks over $\sD$,
which pulls back along $u$ to a $1$--morphism $u^*f : u^*\sP \times_\sD y(u(C))) \to u^*(\sQ)$.
Then, we have a canonical equivalence
$u^*(\sP \times_\sD y(u(C))) \simeq u^*\sP \times_\sC u^*(y(u(C)))$
since $2$--limits commute with one another and $u^*$ is given by a $2$--fibre product.
We have a further equivalence $y(C) \simeq u^*(y(u(C)))$.
Thus, we have a $1$--morphism:
\[
u^*\sP \times_\sC y(C) \simeq u^*(\sP \times_\sD y(u(C))) \to[u^*f] u^*\sQ
\]
The same chain of equivalences allows us to pull back the $2$-cells forming morphisms in
$u^*\HHom(\sP, \sQ)$ to morphisms in $\HHom(u^*\sP, u^*\sQ)$. It is not hard to show that
this mapping $u^*\HHom(\sP, \sQ) \to \HHom(u^*\sP, u^*\sQ)$ is a functor.

When $u$ satisfies the hypotheses of \cref{thm:comp-thm}, $u^*$ is an equivalence
$\St_{\sD} \to[\simeq] \St_{\sC}$.
We thus have a chain of equivalences:
\[
\St_{\sD}(\sP \times_\sD y(u(C)), \sQ)
\to[u^*] \St_{\sC}(u^*(\sP \times_\sD y(u(C))), u^*\sQ)
\to[\simeq] \St_{\sC}(u^*\sP \times_\sC y(C), u^*\sQ)
\]
for every choice of $C \in \Ob{\sC}$. This simultaneously shows that the functor constructed
in the previous paragraph is fully faithful and essentially surjective.
\end{proof}

\begin{thm}\label{thm:map-X-BGl-mod-st-vec-bun}
Let $X$ be a $C^\infty$ manifold. Then, $\HHom(X, B\Gl_n(\bC))$, where $\Gl_n(\bC)$ is seen as a
group object in the category of $C^r$ manifolds (or equipped with the $C^r$ diffeology), is the
moduli stack of rank $n$ $C^r$ vector bundles on $X$. In particular, the moduli stack parametrizing
$C^r$ families of vector bundles on $X$ is $\coprod_{n \in \bN} \HHom(X, B\Gl_n(\bC))$.
\end{thm}
\begin{proof}[Proof sketch]
If we unwrap the definition of the colimit of stacks (the stackification
\cite[\href{https://stacks.math.columbia.edu/tag/02ZN}{Lemma 02ZN}]{stacks-project}
of the $2$--colimit of prestacks), we observe that, for any group object $G$ in the category of
$C^r$ manfiolds, the objects of $\HHom(X, BG)$ over a $C^r$ manifold $U$ consist of precisely the
cocyle data of a $C^r$ principal $G$--bundle on $U \times X$. When $G$ is $\Gl_n(\bC)$, this is
equivalent to the data of a $C^r$ vector bundle of rank $n$ on $U \times X$, via the associated
vector bundle construction with respect to the defining representation of $\Gl_n(\bC)$.
The same holds for morphisms.
We omit further details, as this is a direct adaptation of similar standard results in algebraic
geometry or $C^\infty$ differential geometry.
\end{proof}

\begin{defn}[Moduli Stack of $C^r$ Vector Bundles]
We will denote the moduli stack of $C^r$ vector bundles on $X$ as $\sM_{\Vect, r}(X)$.
\end{defn}

\begin{thm}\label{thm:map-X-BG-dflg}
For every $r = 0, 1, 2, \dots, \infty$, every $C^\infty$ manifold $X$, and every
diffeological group $G$, $\HHom(X, BG)$ is a diffeological stack when considered as an object of
$\St_{\Mfd^r}, \St_{\Cart^r}$ or $\St_{\Dflg}$.
\end{thm}
\begin{proof}
It is a stack by \cref{prop:mapping-prest-st}.
By \cref{prop:dflg-st-pullback-iff}, it suffices to show that it is diffeological
when considered as an object of $\St_{\Mfd^r}$. This follows from a combination
of \cref{thm:BG-dflg} and the argument of \cite[Theorem 4.2, and the following paragraph]{RV18}.
\end{proof}

\begin{rmk}\label{rmk:Frechet-Lie-groupoid-weak-presentation}
$\HHom(X, BG)$ is presented as a diffeological stack by the Fr\'echet Lie groupoid
$\Map(X, BG)_1 \rightrightarrows \Map(X, B\Gl_n(\bC))_0$ of
\cite[2]{RV18HomStack}. This does not mean that $\HHom(X, BG)$ is a Fr\'echet differentiable
stack, because the mapping $\Map(X, B\Gl_n(\bC))_0 \to \HHom(X, BG)$, while representable by
diffeological spaces, might fail to be representable by Fr\'echet manifolds, as the Fr\'echet Lie
groupoid only weakly presents the stack.
\end{rmk}

\begin{thm}\label{thm:mod-st-Vec-Bun-dflg}
For a smooth manifold $X$, the moduli stack $\sM_{\Vect, r}(X)$ of $C^r$ vector bundles on $X$
is diffeological.
\end{thm}
\begin{proof}
Choose diffeological atlases $A_n \to \HHom(X, B\Gl_n(\bC))$, where $\Gl_n(\bC)$ is equipped with
the $C^r$ diffeology, and apply \cref{thm:map-X-BG-dflg} with \cref{prop:subd-st-coprod} to
obtain a diffeological atlas for $\sM_{\Vect, r}(X) = \coprod_{n \in \bN} \HHom(X, B\Gl_n(\bC))$.
\end{proof}

\begin{rmk}\label{rmk:mod-st-Vec-Bun-Frechet-Lie-groupoid-weak-presentation}
For the reason discussed in \cref{rmk:Frechet-Lie-groupoid-weak-presentation},
it is not clear if one can upgade \cref{thm:mod-st-Vec-Bun-dflg} to a statement about Fr\'echet
differentiable stacks.
\end{rmk}

\subsection{$C^{r, \infty}$ Moduli Stacks of Vector Bundles with $\lambda$--Connections}
\label{subsec:M-Flt-r}

\begin{notn}
For two sets $A, B$, any subset $C \subset A \times B$, and any element $a \in A$, we write
$C_a := C \cap (\set{a} \times B)$.
\end{notn}

\begin{defn}[$C^{r, \infty}$ Family of Vector Bundles]
Let $X$ be a smooth manifold and $U$, a $C^r$ manifold.
A $C^{r, \infty}$ family of vector bundles on $X$ parametrized by $U$
is $C^r$ vector bundle $E \to U \times X$ whose restrictions $E|X_u \to X_u$, for each $u \in U$,
are smooth vector bundles.
\end{defn}

\begin{defn}[$C^{r, \infty}$ Sections]\label{defn:Cr-inf-sec}
Let $E \to U \times X$ be a $C^{r, \infty}$ family of vector bundles on a smooth manifold $X$
parametrized by a $C^r$ manifold $U$. Then, for an open $W \subset U \times X$,
we say a $C^r$ section $s : W \to E|W$ is $C^{r, \infty}$ if for each $u \in U$,
the restriction $s_u : W_u \to E|W_u$ is smooth.
We denote the set of $C^{r, \infty}$ sections over $W$ as $\Gamma_{r, \infty}^E(W)$.
A $C^{r, \infty}$ complex function is a $C^{r, \infty}$ section of the trivial bundle
$U \times X \times \bC \to U \times X$.
\end{defn}

\begin{cor}
In the context of \cref{defn:Cr-inf-sec}, the mapping
$\Gamma^E_{r, \infty} : W \mapsto \Gamma^E_{r, \infty}(W)$, with the
usual restriction of functions, is a sheaf. Furthermore, it is a finite rank locally free
$\sA^0_{X/U, r, \infty}$--module.
\end{cor}

\begin{defn}[$C^{r, \infty}$ Differential Forms]
Let $X$ be a smooth manifold and $U$, a $C^r$ manifold.
Let $\pi_X : U \times X \to X$ denote the projection of the product.
We will write $\Gamma^{\pi_X^*\br{\bigwedge^k T^*X_\bC}}_{r, \infty}$ as
$\sA^k_{X/U, r, \infty}$. We will call sections of this sheaf $C^{r, \infty}$ partial differential
$k$--forms.
For another $C^r$ complex vector bundle $E \to U \times X$, we will write
$\Gamma^{E \otimes \pi_X^*\br{\bigwedge^k T^*X_\bC}}_{r, \infty}$ as
$\sA^{k, E}_{X/U, r, \infty}$. Sections of this sheaf will be called 
$C^{r, \infty}$ partial differential $k$--forms valued in $E$.
If $X$ is, in addition, a complex manifold, we will write
$\Gamma^{\pi_X^*(T^*X_\bC^{p, q})}_{r, \infty}$ as $\sA^{p, q}_{X/U, r, \infty}$ and call
its sections $C^{r, \infty}$ partial $(p, q)$--forms.
We will also write
$\Gamma^{E \otimes \pi_X^*\br{T^*X^{p, q}_\bC}}_{r, \infty}$ as
$\sA^{p, q, E}_{X/U, r, \infty}$ and call its sections
its sections $C^{r, \infty}$ partial $(p, q)$--forms valued in $E$.
In particular, $\sA^0_{X/U, r, \infty} = \Gamma^{U \times X \times \bC}_{r, \infty}$
is the sheaf of $C^{r, \infty}$ complex functions.
\end{defn}

\begin{cns}[Pullback of $C^{r, \infty}$--Forms]
Consider a $C^r$ map of $C^r$ manifolds $f : V \to U$ a smooth manifold $X$ and an open
$W \subset U \times X$.
Given a section $s \in \sA^k_{X/U, r, \infty}(W)$, define, for each $(v, x) \in V \times X$:
\[
f^\sharp s(v, x) := s(f(v), x)
\]
It is not difficult to verify that this defines a section
$f^\sharp s \in \sA^k_{X/U, r, \infty}((f \times \id_X)^{-1}W)$, respects
$(p, q)$--decompositions, when $X$ is a complex manifold, and commutes with restriction.
Thus, we get $\bC$--linear sheaf maps:
\[
\sA^k_{X/U, r, \infty} \to (f \times \id_X)_*\sA^k_{X/V, r, \infty},
\sA^{p, q}_{X/U, r, \infty} \to (f \times \id_X)_*\sA^{p, q}_{X/V, r, \infty}
\]
\end{cns}

\begin{defn}[Partial Differentials]\label{defn:prt-diff}
Let $X$ be a smooth manifold and $U$, a $C^r$ manifold.
Suppose we have a $C^{r, \infty}$ section $s \in \sA^k_{X/U, r, \infty}(W)$ over an open
$W \subset U \times X$.
Then, we define:
\[
d_{X/U}^k(s)(u, x) := d(s_u)(x)
\]
where $s_u$ is the restriction of $s$ to $W_u$ and $d$ is the usual exterior
differential of $X$.
It is not difficult to check that this defines a sheaf map
$\sA^k_{X/U, r, \infty} \to \sA^{k + 1}_{X/U, r, \infty}$.
Then, the collection of maps
\[
d_{X/U}^\bullet : \sA^\bullet_{X/U, r, \infty} \to \sA^{\bullet + 1}_{X/U, r, \infty}
\]
will be called the $C^{r, \infty}$ partial exterior differential of $U \times X$ along $X$.
Given a complex structure $\oprt$ on $X$, we define:
\[
\oprt^k_{X/U}(s)(u, x) := \oprt(s_u)(u, x), \prt^k_{X/U}(s)(u, x) := \prt(s_u)(u, x)
\]
We get similar sheaf maps:
\[
\oprt_{X/U}^{\bullet} : \sA^{\bullet}_{X/U, r, \infty}
    \to \sA^{\bullet + 1}_{X/U, r, \infty},
\prt_{X/U}^{\bullet} : \sA^{\bullet}_{X/U, r, \infty}
    \to \sA^{\bullet + 1}_{X/U, r, \infty}
\]
which decompose to sheaf maps as follows for each $p, q$:
\[
\oprt_{X/U}^{p, \bullet} : \sA^{p, \bullet}_{X/U, r, \infty}
    \to \sA^{p, \bullet + 1}_{X/U, r, \infty},
\prt_{X/U}^{\bullet, q} : \sA^{\bullet, q}_{X/U, r, \infty}
    \to \sA^{\bullet + 1, q}_{X/U, r, \infty}
\]
Furthermore, we have:
\[
d^\bullet_{X/U} = \oprt_{X/U}^{\bullet} + \prt_{X/U}^{\bullet}
\]
We will supress the $\bullet$'s whenever there is no chance of confusion.
\end{defn}

\begin{cor}
For a smooth manifold $X$ and a $C^r$ manifold $U$, we have the following differential
$\bC$--spaces \cite[Definition 2.2.1]{AR26a}:
\begin{enumerate}
\item
$(U \times X, \sA^0_{X/U, r, \infty}, \sA^\bullet_{X/U, r, \infty}, d_{X/U}^\bullet)$
\item
$(U \times X, \sA^0_{X/U, r, \infty}, \sA^\bullet_{X/U, r, \infty}, \oprt_{X/U}^\bullet)$
\item
$(U \times X, \sA^0_{X/U, r, \infty}, \sA^\bullet_{X/U, r, \infty}, \prt_{X/U}^\bullet)$
\end{enumerate}
For a $C^r$ map $f : V \to U$, pullback of differential forms along $f \times \id_X$ yields
morphisms of differential $\bC$ spaces of the form above similar to \cite[\S 2]{AR26a}.
\end{cor}
\begin{proof}
This can be seen by computations in trivializations.
\end{proof}

\begin{rmk}
The difference with \cite[\S 2]{AR26a} is we are now working with $C^r$ manifolds and
$C^{r, \infty}$ sections in place of smooth ones.
\end{rmk}

\begin{cns}\label{cns:M-Flt-r}
Let $\fH$ be the category with the following description:
\begin{itemize}
\item Objects are tuples $(Y, R, K, d_{K, 1}, d_{K, 2}, \lambda_1, \lambda_2)$ where
    \begin{enumerate}
    \item $(Y, R, K, d_{K, i}), i = 1, 2$ are a differential $\bC$--spaces.
    \item $\lambda \in \ker{d_{K, 1}} \subset R(Y), \lambda' \in \ker{d_{K, 2}} \subset R(Y)$.
    \end{enumerate}
\item Morphisms $(Z, S, L, d_{L, 1}, d_{L, 2}, \tau_1, \tau_2) \to
(Y, R, K, d_{K, 1}, d_{K, 2}, \lambda_1, \lambda_2)$ are morphisms of differential spaces
with $(f, f^\sharp) : (Z, S, L, d_{L, i}) \to (Y, R, K, d_{K, i}), i = 1, 2$, with the 
same underlying maps $f : Z \to Y$ and $f^\sharp : K^{\wedge \bullet} \to K^{\wedge \bullet + 1}$,
and such that $f^\sharp \lambda_i = \tau_i, i = 1, 2$.
\end{itemize}
Composition and identities are given by those of morphisms of differential spaces.
$(f, f^\sharp)$ is also a morphism of differential spaces
\[
(Z, S, L, \tau_1d_{L, 1} + \tau_2d_{L, 2}) \to (Y, R, K, \lambda_1d_{K, 1} + \lambda_2d_{K, 2})
\]
This allows us to define a pseudofunctor:
\[
\Conn : \fH^\op \to \bC\dg-\Cat : (Y, R, K, d_{K, 1}, d_{K, 2}, \lambda_1, \lambda_2)
\mapsto \Conn_{1, \lambda_1d_{K, 1} + \lambda_2d_{K, 2}}(Y)
\]
as in \cite[Theorem 3.3.4]{AR26a}.
We also have a functor
\[\begin{array}{ccccc}
\Flt_X &:& \Mfd^r_{/\bC^2} &\to& \fH \\
       &:& ((\lambda_1, \lambda_2) : U \to \bC^2) &\mapsto&
(U \times X, \sA^0_{X/U, r, \infty}, \sA^\bullet_{X/U, r, \infty}, \prt_{X/U}, \oprt_{X/U},
\pi_U \circ \lambda_1, \pi_U \circ \lambda_2)
\end{array}\]
We denote the composite pseudofunctor:
\[
\Mfd^{r, \op}_{/\bC^2} \to[\Flt^\op_X] \fH^\op \to[\Conn] \bC\dg-\Cat
\]
by $\mcM_{\Flt, r}(X)$.
By taking homotopy categories, followed by cores and then applying the Grothendieck construction, we
get a prestack over $\Mfd^r_{\bC^2}$:
\[
\sM_{\Flt, r}(X) := \int \br{\Mfd^{r, \op}_{/\bC^2} \to[\Flt^\op_X] \fH^\op \to[\Conn] \bC\dg-\Cat
\to[h] \Cat \to[(-)^\simeq] \Grpd} \to \Mfd^r_{/\bC^2}
\]
Of course, this also makes the composite $\sM_{\Flt, r}(X) \to \Mfd^r_{/\bC^2} \to \Mfd^r$
a prestack. We will call $\sM_{\Flt, r}(X)$ the bidirectional $C^r$ non-Abelian Hodge filtration.
\end{cns}

\begin{notn}
Given a morphism of $C^r$ manifolds $f : V \to U$, we will denote the pullback of various
objects along $f \times \id_X : V \times X \to U \times X$ by $f^<$ instead of
$(f \times \id_X)^*$ for brevity.
\end{notn}

\begin{rmk}\label{rmk:M-flt-r-concrete}
First, observe that for any $C^r$ manifold $U$ for $r = 0, 1, 2, \dots, \infty$,
and smooth manifold, $X$, we have an equivalence of categories between finite rank locally free
$\sA_{X/U, r, \infty}^0$--modules and $C^{r, \infty}$ families of complex vector bundles over $X$
parametrized by $U$. This follows from the same cocycle data argument giving the equivalence of
smooth vector bundles with finite rank locally free $C^\infty_X$--modules.
Under this equivalence, $\sM_{\Flt, r}(X)$ has the following concrete description:
\begin{itemize}
\item Objects are tuples $(U, \lambda_1, \lambda_2, E, \nabla)$ where:
    \begin{itemize}
    \item $U$ is a $C^r$ manifold.
    \item $\lambda_i : U \to \bC, i = 1, 2$ are $C^r$ functions.
    \item $E \to U \times X$ is a $C^{r, \infty}$ family of complex vector bundles on $X$
    parametrized by $U$.
    \item $\nabla$ is a flat
    $1$--$((\pi_U \circ \lambda_1)\prt_{X/U} + (\pi_U \circ \lambda_2)\oprt_{X/U})$--connection on
    $E$.
    \end{itemize}
\item Morphisms $(V, \tau_1, \tau_2, E', \nabla') \to (U, \lambda_1, \lambda_2, E, \nabla)$ are
    tuples $(f, \alpha)$ where:
    \begin{itemize}
    \item $f : V \to U$ is a $C^r$ function.
    \item $\tau_i = \lambda_i \circ f, i = 1, 2$.
    \item $\alpha : (f^<E, f^<\nabla) \to (E', \nabla')$ is an isomorphism of
    $1$--$((\pi_U \circ \tau_1)\prt_{X/V} + (\pi_U \circ \tau_2)\oprt_{X/V})$--connections.
    \end{itemize}
\end{itemize}
\end{rmk}

\begin{defn}\label{defn:M-Dol-dR-Hodge}
Consider the map $P : \bC \to \bC^2 : a \mapsto (1, a)$.
Also, consider the map $P_{a, b} : \pt \to \bC^2$ given by sending $0 \in \pt = \bR^0$
to $(a, b) \in \bC^2$. 
These provide maps
$y(P) : \Mfd_{/\bC} \to \Mfd_{/\bC^2}$ and
$y(P_{a, b}) : \Mfd \cong \Mfd_{/\pt} \to \Mfd_{/\bC^2}$.
Then, it is not hard to see from the above description of $\sM_{\Flt, r}$ along with
\cref{rmk:M-flt-r-concrete} and
\cref{prop:strictify-2-fibre-prod} that
the pullback of $\sM_{\Flt, r}$ along the maps $y(P_{0, 1})$ and $y(P_{1, 1})$
are the moduli prestacks parametrizing $C^{r, \infty}$ families of Higgs bundles and flat bundles
respectively. We will denote these by $\sM_{\Dol, r}(X)$ and $\sM_{\dR, r}(X)$.
Similarly, the pullback with respect to the map $y(P)$ is the $C^r$ non-Abelian Hodge filtration
stack $\sM_{\Hodge, r}(X)$.
\end{defn}

\begin{notn}
For brevity, we will abuse notation slightly and use $\lambda$ to denote $\pi_U \circ \lambda$
as well.
\end{notn}

\begin{prop}\label{prop:M-flt-r-stack}
$\sM_{\Flt, r}(X)$ is a stack over $\Mfd^r$.
\end{prop}
\begin{proof}
An object descent datum consists of an open cover $U = \bigcup_{j \in J} U_j$, objects
$(U_i, \lambda_{1, j}, \lambda_{2, j}, E_j, \nabla_j)$ such that
$\lambda_{i, j}|(U_j \cap U_k) = \lambda_{i, k}|(U_j \cap U_k)$, and isomorphisms of
$1$--$(\lambda_{1, j}|(U_j \cap U_k)\prt_{X/U}
+ \lambda_{2, j}|(U_j \cap U_k)\oprt_{X/U})$--connections
$\alpha_{jk} : (E_j, \nabla_j)|(U_j \cap U_k) \to[\cong] (E_k, \nabla_k)|(U_j \cap U_k)$
that satisfy the cocycle condition.
This implies there exist $\lambda_1, \lambda_2 : U \to \bC$ with $\lambda_1|U_j = \lambda_{1, j}$
and $\lambda_2|U_j = \lambda_{2, j}$.
Then, we observe that
$(\lambda_1 \prt_{X/U} + \lambda_2 \oprt_{X/U})|U_j = \lambda_{1, j}\prt_{X/U}
+ \lambda_{2, j}\oprt_{X/U}$.
Collation then follows by an argument similar to the proof of descent in $\sM_\dR(X)$ given in
\cite[Theorem 3.4.1]{AR26a}.

A morphism descent datum consists of a $C^r$ function $f : U \to V$, and open cover
$V = \bigcup_{j \in J} V_j$, objects $B = (U, \tau_1, \tau_2, E', \nabla')$ and
$A = (V, \lambda_1, \lambda_2, E, \nabla)$, morphisms $(f|U_j, \alpha_j) : A|V_j \to B$
that agree on $V_j \cap V_k$ for all $j, k \in J$. Then, this gives a morphism descent datum
in $\sM_{\Vect, r}(X)$ which collate to a morphism $(f, \alpha) : (V, E') \to (U, E)$.
Local computations then show that $\alpha$ is also a morphism of connections showing
that $(f, \alpha)$ is also a collation of morphisms in $\sM_{\Flt, r}(X)$.
\end{proof}

\begin{thm}\label{thm:M-Flt-r-dflg}
The forgetful map
\[
\sM_{\Flt, r}(X) \to \Mfd^r_{/\bC^2} \times \sM_{\Vect, r}(X)
: (U, \lambda_1, \lambda_2, E, \nabla) \mapsto (U, \lambda_1, \lambda_2, E)
\]
is representable by diffeological spaces. In particular, $\sM_{\Flt, r}(X)$ is diffeologically
presentable.
\end{thm}
\begin{proof}
Consider a $C^r$ manifold $U$ and a morphism $U \to \Mfd^r_{/\bC^2} \times \sM_{\Vect, r}(X)$
corresponding to an object
$(U, \lambda_1, \lambda_2, E)$ of $\Mfd^r_{/\bC} \times \sM_{\Vect, r}(X)$.
We will show that the fibre product stack
\[
P_{E, \lambda_1, \lambda_1}
:= \sM_{\Flt, r}(X) \times_{\Mfd^r_{/\bC^2} \times \sM_{\Vect, r}(X)} U
\]
is representable by a diffeological space.
This stack has the following concrete description by
\cite[\href{https://stacks.math.columbia.edu/tag/0040}{Tag 0040}]{stacks-project}:
\begin{itemize}
\item Objects are $(v, F, \nabla_F, \alpha)$ where:
    \begin{itemize}
    \item $v : V \to U$ is a $C^r$ map from a smooth manifold $V$,
    \item $(F, \nabla_F)$ is a flat
    $((\lambda_1 \circ v) \prt_{X/V} + (\lambda_2 \circ v)\oprt_{X/V})$--connection, and
    \item $\alpha : v^<E \to F$ is an isomorphism of vector bundles.
    \end{itemize}
\item Morphisms $(w, G, \nabla_G, \beta) \to (v, F, \nabla_F, \alpha)$ are tuples $(f, a)$ where:
    \begin{itemize}
    \item $f : W \to V$ is a $C^r$ map commuting with the maps to $U$, and
    \item $a : f^<(F, \nabla_F) \to (G, \nabla_G)$ is an isomorphism of
    $((\lambda_1 \circ w)\prt_{X/V} + (\lambda_2 \circ w)\oprt_{X/V})$--connections,
    such that the following diagram of vector bundles commutes:
    \[\begin{tikzcd}
    f^<v^<E \ar[r, "f^<\alpha"] \ar[d, "\cong" left] & f^<F \ar[d, "a"] \\
    w^<E \ar[r, "\beta" below] & G
    \end{tikzcd}\]
    \end{itemize}
\end{itemize}
In the above diagram, since all maps except $a$ are isomorphisms, there is a unique such $a$.
Thus, $P_{E, \lambda_1, \lambda_2}$ is a prestack in setoids.

We will show that this is equivalent to
$\Mfd^r_{/T}$ for a diffeological subspace $T \subset S_{E, \lambda_1, \lambda_2}$ of the fibre
product of diffeological spaces:
\[\begin{tikzcd}
S_{E, \lambda_1, \lambda_2} \ar[r] \ar[d] \ar[rd, phantom, "\lrcorner" very near start] &
\Dflg(X, \End(E) \otimes \sA^1_{X/U, r, \infty}) \ar[d, "p_*"] \\
U \ar[r, "\hat{\id}_{U \times X}" below] & \Dflg(X, U \times X)
\end{tikzcd}\]
where the bottom map is the unique map corresponding to $\id_{U \times X}$ and the right vertical
map is induced by the vector bundle projection
$p : \End(E) \otimes \sA_{X/U, r, \infty}^1 \to U \times X$.
First, fix a $(\lambda_1\prt_{X/U} + \lambda_2\oprt_{X/U})$--connection $\nabla$, not necessarily
flat, on $E \to U \times X$ --- this exists by a partition of unity argument identical to the case
of ordinary connections.
Let $(v, F, \nabla_F, \alpha)$ be any object of $P_{E, \lambda_1, \lambda_2}$.
Then, we have the map $v : V \to U$ along with the $\sA_{X/U, r, \infty}^1$--linear map
\[
v^<\nabla - (\alpha^{-1} \otimes \id) \circ \nabla_F \circ \alpha :
v^<E \to v^<E \otimes \sA_{X/V, r, \infty}^1
\]
This corresponds to a section $s : V \times X \to \End(v^<E) \otimes \sA_{X/V, r, \infty}^1$.
We have a canonical isomorphism
$\sA_{X/V, r, \infty}^1 = \pi_X^*\sA^1_X \cong v^<\sA_{X/U, r, \infty}^1$.
Thus, $s$ gives a map of diffeological spaces:
\[
s'_{\nabla_F} : V \times X \to[s] \End(v^<E) \otimes v^<\sA_{X/U, r, \infty}^1
\to \End(E) \otimes \sA_{X/U, r, \infty}^1
\]
such that $p \circ s' = v \times \id_X$. By the definition of the fibre product,
$(s'_{\nabla_F}, v)$ yields a morphism of diffeological spaces
$\Psi(v, F, \nabla_F, \alpha) : V \to S_{E, \lambda_1, \lambda_2}$.
Then, consider a morphism $(f, a) : (w, G ,\nabla_G, \beta) \to (v, F, \nabla_F, \alpha)$
in $P_{E, \lambda_1, \lambda_2}$. Let $\gamma : f^<v^<E \to w^<E$ denote the canonical isomorphism.
Then, we observe that:
\begin{align*}
 & w^<\nabla - (\beta^{-1} \otimes \id)\nabla_G\beta \\
=& \gamma \circ f^<v^<\nabla \circ \gamma^{-1}
   - ((\gamma f^<\alpha^{-1} a^{-1}) \otimes \id) \nabla_G (a f^<\alpha \gamma^{-1}) \\
=& \gamma \circ (f^<v^<\nabla
   - f^<(\alpha^{-1} \otimes \id) (f^<\nabla_F) (f^<\alpha)) \circ \gamma^{-1} \\
=& \gamma \circ f^<(v^<\nabla
   - (\alpha^{-1} \otimes \id) \nabla_F \alpha) \circ \gamma^{-1}
\end{align*}
where the last step can be shown by a local computation. This implies that the diagram on the
left below commutes, which, in turn, implies the diagram on the right commutes:
\[\begin{tikzcd}[column sep=small]
W \times X \ar[rr, "f \times \id_X"] \ar[rd, "s'_{\nabla_G}" below left] & &
V \ar[ld, "s'_{\nabla_F}" below right] \\
& \End(E) \otimes \sA_{X/U, r, \infty}^1 &
\end{tikzcd}
\hspace{1.5em}
\begin{tikzcd}[column sep=small]
W \ar[rr, "f"] \ar[rd, "{\Psi(w, G, \nabla_G, \beta)}" below left] & &
V \ar[ld, "{\Psi(v, F, \nabla_F, \beta)}" below right] \\
& S_{E, \lambda_1, \lambda_2} &
\end{tikzcd}\]
We then set $\Psi(f, a) = f$. It is immediate that $\Psi$ is a morphism of prestacks
$P_{E, \lambda_1, \lambda_2} \to S_{E, \lambda_1, \lambda_2}$. By the uniqueness of morphisms
in $P_{E, \lambda_1, \lambda_2}$ covering a morphism in $\Mfd^r$ discussed in the previous
paragraph, $\Psi$ is fully faithful.
The result now follows from \cref{prop:subfunc-dflg}.
The last claim follows from \cref{cor:dflg-rep-morphism-dflg-target}.
\end{proof}

\begin{thm}\label{thm:M-Dol-dR-Hodge-dflg}
The prestacks $\sM_{\Dol, r}(X)$, $\sM_{\dR, r}(X)$ and $\sM_{\Hodge, r}(X)$ are stacks.
The forgetful maps
\[
\sM_{\Dol, r}(X) \to \sM_{\Flt, r}(X),
\sM_{\dR, r}(X) \to \sM_{\Flt, r}(X)
\text{ and } \sM_{\Hodge, r}(X) \to \sM_{\Flt, r}(X)
\]
are diffeologically representable.
In particular, their domains are diffeological stacks.
\end{thm}
\begin{proof}
Descent can be deduced from the fact that these prestacks are fibre products of stacks by
\cref{prop:M-flt-r-stack} and \cref{defn:M-Dol-dR-Hodge}.
We then observe that the forgetful maps
\[
\sM_{\Dol, r}(X) \to \sM_{\Flt, r}(X),
\sM_{\dR, r}(X) \to \sM_{\Flt, r}(X)
\text{ and } \sM_{\Hodge, r}(X) \to \sM_{\Flt, r}(X)
\]
are pullbacks of the maps $y(P_{0, 1}), y(P_{1, 1}) : \Mfd \to \Mfd_{/\bC^2}$ and
$y(P) : \Mfd_{/\bC} \to \Mfd_{/\bC^r}$ respectively
from \cref{defn:M-Dol-dR-Hodge}, and are hence, representable by diffeological spaces by
\cref{prop:subd-st-pullback}. \cref{cor:dflg-rep-morphism-dflg-target} along with
\cref{thm:M-Flt-r-dflg} then yields the final statement.
\end{proof}

\subsection{$C^{r, \infty}$ Moduli Stacks of Bundle and Connection Diagrams}
\label{subsec:M-Diag-r}

\begin{defn}
We will denote by $[n]$ the the totally ordered set
$\set{0 < 1 < \cdots < n}$. The category of these sets and order preserving morphisms
is the simplicial indexing category $\Delta$.
There is a functor $\Delta \to \Cat$ defined by sending $[n]$ to category whose objects
are $0, 1, \dots, n$ with a unique morphism $i \to j$ if $i < j$.
We will denote the image of $[n]$ under this functor by $\Delta^n$.
\end{defn}

\begin{defn}
For any two categories $\sC$ and $\sD$, we will denote the constant presheaf of categories
on $\sC$ valued in $\sD$ as $\ul{\sD} : \sC^\op \to \Cat$.
We will also denote its Grothendieck construction as $\ul\sD \to \sC$.
\end{defn}

\begin{defn}
Let $(\sC, J)$ be a site and $\sP \in \PSt_{\sC, J}$.
We denote by $S_n\sP$ the Hom prestack $\HHom(\ul{\Delta^n}, \sP)$.
\end{defn}

\begin{cor}\label{cor:quiver-prestack-simplicial}
If $\sP$ is a prestack over a category $\sC$, then
$S_\bullet\sP := \HHom(\ul{\Delta^\bullet}, \sP)$ is a simplicial prestack.
\end{cor}
\begin{proof}
The mapping $[n] \mapsto S_n\sP$ is a composite of strict pseudofunctors:
\[
\Delta^\op \to[\Delta^{-, \op}] \Cat^\op \to[\ul{(-)}] \PSt_\sC^\op \to[\HHom(-, \sP)] \PSt_\sC
\]
where the first map sends a totally ordered set to its associated category.
\end{proof}

\begin{cor}
If $\mcX$ is a stack over a site, then $S_n\mcX$ is stack.
In particular, $S_\bullet\mcX$ is a simplicial stack.
\end{cor}
\begin{proof}
Apply \cref{prop:mapping-prest-st} and \cref{cor:quiver-prestack-simplicial}.
\end{proof}

\begin{notn}
We will write
$S_\bullet\sM_{\Flt, r}(X)$, $S_\bullet\sM_{\Vect, r}(X)$, $S_\bullet\sM_{\Dol, r}(X)$,
$S_n\sM_{\Hodge, r}(X)$ and $S_n\sM_{\dR, r}(X)$ as
$S_\bullet\Flt_r(X), S_\bullet\Vect_r(X), S_\bullet\Dol_r(X)$, $S_\bullet\dR_r(X)$
and $S_n\Hodge_r(X)$ respectively, for brevity.
\end{notn}

\begin{thm}\label{thm:simplicial-filtration-dflg}
For a complex manifold $X$, the morphism of stacks
\[\begin{tikzcd}
S_1\Flt_r(X) = \HHom(\ul{\Delta^1}, \sM_{\Flt, r}(X)) \ar[d, "{(\ev_0, \ev_1)}"] \\
S_0\Flt_r(X) \times S_0\Flt_r(X) \cong \sM_{\Flt, r}(X) \times \sM_{\Flt, r}(X)
\end{tikzcd}\]
given by pre-composition with the two inclusions
$\set{i} \hto \sbr{1} = \set{0 < 1}, i = 0 , 1$ is diffeologically representable.
In particular, $S_n\Flt_r(X)$ is diffeological for all $n = 0, 1, 2, \dots$ --- that is,
$S_\bullet\Flt_r(X)$ is a simplicial diffeological stack.
\end{thm}
\begin{proof}
Consider a map $U \to \sM_{\Flt, r}(X) \times \sM_{\Flt, r}(X)$ corresponding to
$C^{r}$ functions $\lambda_1, \lambda_2, \lambda'_1, \lambda'_2 : U \to \bC$,
$C^{r, \infty}$ vector bundles $E, F \to U \times X$, a
$(\lambda_1\prt_{X/U} + \lambda_2\oprt_{X/U})$--connection $\nabla$ on $E$ and
a $(\lambda'_1\prt_{X/U} + \lambda'_2\oprt_{X/U})$--connection $\nabla'$ on $F$.
Then, an argument similar to \cref{thm:M-Flt-r-dflg} shows that
the fibre product $S_1\Flt_r(X) \times_{\sM_{\Flt, r}(X)^2} U$ is equivalent to the stack
represented by a diffeological subspace of the following fibre product of diffeological spaces:
\[\begin{tikzcd}
S_{E, F} \ar[r] \ar[d] & \Dflg(X, E^\vee \otimes F) \ar[d, "p_*"] \\
U \ar[r, "\hat{\id}_{U \times X}" below] & \Dflg(X, U \times X)
\end{tikzcd}\]
where the right vertical map is given by post-composition with the vector bundle projection
$p : E^\vee \otimes F \to U \times X$ and the bottom map is induced by $\id_{U \times X}$.
The fact that $S_1\Flt_r(X)$ is diffeological follows from
\cref{cor:dflg-rep-morphism-dflg-target}.
By arguments similar to \cite[\S 4]{ModStQuivBun}, we can show that $S_n\Flt_r(X)$ can be computed
using finite limits involving $S_k\Flt_r(X), k = 0, 1$ and $\bR^0$.
The final statement then follows from \cref{thm:fin-lim-dflg-st}.
\end{proof}

\begin{cor}\label{cor:simplicial-Dol-dR-Hodge-dflg}
For a complex manifold $X$, the simplicial stacks
$S_\bullet\Dol_r(X)$, $S_\bullet\dR_r(X)$ and $S_n\Hodge_r(X)$
are diffeological at all levels.
\end{cor}
\begin{proof}
For each $n$, we have a map $S_n\Flt_r(X) \to \Mfd^r_{/\bC^{2(n + 1)}}$.
Let $D : \bC^2 \to \bC^{2(n + 1)}$ be the diagonal map and recall the maps
$P_{0, 1}, P_{1, 1} : \pt \to \bC^2$, $P : \bC \to \bC^2$ from \cref{defn:M-Dol-dR-Hodge}.
This yields maps $y(D \circ P_{0, 1}), y(D \circ P_{1, 1}) : \Mfd^r \to \Mfd^r_{/\bC^{2(n + 1)}}$.
The fibres products
\[
S_n\Flt_r(X) \times_{\Mfd^r_{/\bC^{2(n + 1)}}} \Mfd^r
\]
with respect to these two maps are $S_n\Dol_r(X)$ and $S_n\dR_r(X)$ respectively.
The result now follows from \cref{thm:simplicial-filtration-dflg} and \cref{thm:fin-lim-dflg-st}.
The same argument with $P$ in place of $P_{0, 1}$ or $P_{1, 1}$ yields the result for
$S_\bullet\Hodge_r(X)$.
\end{proof}

\begin{thm}\label{thm:simplicial-Vect-dflg}
For a smooth manifold $X$, the morphism of stacks
\[\begin{tikzcd}
S_1\Vect_r(X) = \HHom(\ul{\Delta^1}, \sM_{\Vect, r}(X)) \ar[d, "{(\ev_0, \ev_1)}"] \\
S_0\Vect_r(X) \times S_0\Vect_r(X) \cong \sM_{\Vect, r}(X) \times \sM_{\Vect, r}(X)
\end{tikzcd}\]
given by pre-composition with the two inclusions
$\set{i} \hto \sbr{1} = \set{0 < 1}, i = 0 , 1$ is diffeologically representable.
In particular, $S_n\Vect_r(X)$ is diffeological for all $n = 0, 1, 2, \dots$ --- that is,
$S_\bullet\Vect_r(X)$ is a simplicial diffeological stack.
\end{thm}
\begin{proof}
Similar to \cref{thm:simplicial-filtration-dflg}, but the fibres
$S_1\Vect_r(X) \times_{\sM_{\Vect, r}(X)} U$ are equivalent to the stack represented by all of
$S_{E, F}$ in this case.
\end{proof}

\begin{rmk}
Let $I$ be a finite simplicial set. We can define moduli stacks
$S_I\Vect_r(X)$, $S_I\Flt_r(X)$, $S_I\Dol_r(X)$ and $S_I\dR_r(X)$ parametrizing the
$I$--shaped diagrams of the respective objects similarly. These are diffeological
by the observation that $I$ can be computed as a finite limit involving the $S_n$ stacks
above, which we showed to be diffeological, and applying \cref{thm:fin-lim-dflg-st}.
We omit the details.
\end{rmk}

\printbibliography

\vspace{5em}

\end{document}

%% file: preamble.tex
\usepackage{amsmath}
\usepackage{amsthm}
\usepackage{amsfonts}
\usepackage{amssymb}
\usepackage{amsthm}
\usepackage{bbm}
\usepackage{enumitem}
\usepackage[mathscr]{euscript}
\usepackage{float}
\usepackage{ifthen}
\usepackage{listings}
\usepackage{mathdots}
\usepackage{mleftright}
\usepackage[new]{old-arrows}
\usepackage{rotating}
\usepackage{sansmath}
\usepackage{stmaryrd}
\usepackage{tabularx}
\usepackage{tikz}
  \usetikzlibrary{positioning}
\usepackage{tikz-cd}
  \tikzcdset{
    cells={font=\everymath\expandafter{\the\everymath\displaystyle}},
  }
  \usetikzlibrary{patterns}
\usepackage{todonotes}
\usepackage{xcolor}
\usepackage{xspace}

\usepackage[capitalize]{cleveref}

\definecolor{codegreen}{rgb}{0,0.6,0}
\definecolor{codegray}{rgb}{0.5,0.5,0.5}
\definecolor{codepurple}{rgb}{0.58,0,0.82}
\definecolor{backcolour}{rgb}{0.95,0.95,0.92}
\lstdefinestyle{mystyle}{
    xleftmargin={2.5em},
    backgroundcolor=\color{backcolour},
    commentstyle=\color{codegreen},
    keywordstyle=\color{magenta},
    numberstyle=\tiny\color{codegray},
    stringstyle=\color{codepurple},
    basicstyle=\ttfamily\footnotesize,
    breakatwhitespace=false,
    breaklines=true,
    captionpos=b,
    keepspaces=true,
    numbers=left,
    numbersep=5pt,
    showspaces=false,
    showstringspaces=false,
    showtabs=false,
    tabsize=2
}
\usepackage{subfiles}

\newtheorem{thm}{Theorem}
\numberwithin{thm}{subsection}

\newtheorem{prop}[thm]{Proposition}

\newtheorem{cor}[thm]{Corollary}

\theoremstyle{definition}

\newtheorem{cns}[thm]{Construction}
\newtheorem{conv}[thm]{Convention}
\newtheorem{exm}[thm]{Example}
\newtheorem{defn}[thm]{Definition}
\newtheorem{notn}[thm]{Notation}
\newtheorem{rmk}[thm]{Remark}
\newtheorem{warn}[thm]{Warning}

\numberwithin{qsn}{section}

\newcommand{\bb}{\mathbb}

\newcommand{\f}{\mathfrak}
\newcommand{\s}{\mathscr}

\newcommand{\mc}{\mathcal}

\newcommand{\mcM}{\mc{M}}

\newcommand{\mcO}{\mc{O}}

\newcommand{\mcX}{\mc{X}}
\newcommand{\mcY}{\mc{Y}}
\newcommand{\mcZ}{\mc{Z}}

\newcommand{\fH}{\f{H}}

\newcommand{\sA}{\s{A}}
\newcommand{\sB}{\s{B}}
\newcommand{\sC}{\s{C}}
\newcommand{\sD}{\s{D}}

\newcommand{\sM}{\s{M}}

\newcommand{\sP}{\s{P}}
\newcommand{\sQ}{\s{Q}}

\newcommand{\bC}{\bb{C}}

\newcommand{\bN}{\bb{N}}

\newcommand{\bR}{\bb{R}}

\newcommand{\wh}{\widehat}
\newcommand{\br}[1]{\mleft( #1 \mright)}
\newcommand{\sbr}[1]{\mleft[ #1 \mright]}

\newcommand{\set}[2][]{
  \ifthenelse{\equal{#1}{}}{
    \mleft\{ #2 \mright\}
  }{
    \mleft\{ #1\ :\ #2 \mright\}
  }
}

\renewcommand{\to}[1][]{
  \ifthenelse{\equal{#1}{}}{
    \longrightarrow
  }{
    \stackrel{#1}{\longrightarrow}
  }
}

\newcommand{\To}[1][]{
  \ifthenelse{\equal{#1}{}}{
    \Longrightarrow
  }{
    \stackrel{#1}{\Longrightarrow}
  }
}

\renewcommand{\mapsto}[1][]{
    \ifthenelse{\equal{#1}{}}{
      \longmapsto
    }{
      \stackrel{#1}{\longmapsto}
    }
}

\newcommand{\ot}[1][]{
  \ifthenelse{\equal{#1}{}}{
  \longleftarrow
  }{
    \stackrel{#1}{\longleftarrow}
  }
}

\newcommand{\ol}{\overline}
\newcommand{\ul}{\underline}

\renewcommand{\ker}[1]{\text{ker}\br{#1}}
\newcommand{\im}{\mathrm{im}\xspace}

\newcommand{\Gl}{\mathrm{Gl}}

\newcommand{\End}{\mathrm{End}}

\newcommand{\hto}[1][]{\stackrel{#1}{\longhookrightarrow}}

\newcommand{\id}{\mathrm{id}}
\newcommand{\ev}{\mathrm{ev}}
\newcommand{\prt}{\partial}
\newcommand{\oprt}{\ol{\prt}}

\renewcommand{\Ref}[2][]{\ifthenelse{\equal{#1}{}}{\ref{#2}}
                      {\hyperref[#2]{\ref*{#1}(\ref*{#2})}}}
\newcommand{\Aref}[2][]{\ifthenelse{\equal{#1}{}}{\autoref{#2}}
                      {\hyperref[#2]{\autoref*{#1}\ref*{#2}}}}
\newcommand{\Sref}[1]{\hyperref[#1]{\S \ref*{#1}}}
\newcommand{\dom}{\mathrm{dom}}

\newcommand{\Set}{\mathrm{Set}}

\newcommand{\Grpd}{\mathrm{Grpd}}
\newcommand{\op}{\mathrm{op}}

\newcommand{\Sh}{\mathrm{Sh}}
\newcommand{\Fun}{\mathrm{Fun}}
\newcommand{\Top}{\mathrm{Top}}

\newcommand{\HHom}{\mathcal{H}\mathrm{om}}

\newcommand{\Map}{\mathrm{Map}}

\newcommand{\Sc}[1][]{\ifthenelse{\equal{#1}{}}{\mathrm{Sch}}{\mathrm{Sch}_{/#1}}}

\newcommand{\Cat}{\s{C}\mathrm{at}}

\newcommand{\Vect}{\mathrm{Vect}}

\newcommand{\Conn}{\mathrm{Conn}}

\newcommand{\Hodge}{\mathrm{Hodge}}

\newcommand{\St}{\mathrm{St}}
\newcommand{\PSt}{\mathrm{PSt}}

\newcommand{\Dol}{\mathrm{Dol}}
\newcommand{\Flt}{\mathrm{Flt}}
\newcommand{\dR}{\mathrm{dR}}

\newcommand{\dg}{\mathrm{dg}}

\newcommand{\Ob}[1]{\mathrm{Ob}\br{#1}}
\newcommand{\Mor}[1]{\mathrm{Mor}\br{#1}}

\newcommand{\pt}{\mathrm{pt}}

\newcommand{\Dflg}{\mathrm{Dflg}}

\newcommand{\Cart}{\mathrm{Cart}}
\newcommand{\Mfd}{\mathrm{Mfd}}

\newcommand{\frgt}{\mathrm{frgt}}
\newcommand{\st}{\mathrm{st}}

%% file: refs.bib
@misc{AR26a,
title={A diffeological perspective on non-Abelian Hodge theory}, 
author={Mahmud Azam and Steven Rayan},
year={2026},
eprint={2606.16772},
archivePrefix={arXiv},
primaryClass={math.DG},
url={https://arxiv.org/abs/2606.16772}, 
note={In review (``Journal of Differential Geometry'', since June 15, 2026)}
}

@misc{AR26b,
title={Diffeological non-Abelian Hodge theory: relative harmonic metrics and deformation theory}, 
author={Mahmud Azam and Steven Rayan},
year={2026},
eprint={2607.18989},
archivePrefix={arXiv},
primaryClass={math.DG},
url={https://arxiv.org/abs/2607.18989}, 
note={In review (``Annals of Mathematics'', since July 23, 2026)}
}

@article{ModStQuivBun,
title = {Moduli stacks of quiver bundles with applications to Higgs bundles},
author = {Mahmud Azam and Steven Rayan},
journal = {Journal of Algebra},
volume = {708},
pages = {333-379},
year = {2026},
issn = {0021-8693},
doi = {https://doi.org/10.1016/j.jalgebra.2026.05.018},
url = {https://www.sciencedirect.com/science/article/pii/S0021869326003157}
}

@Inbook{GR15,
author="Garc{\'i}a-Raboso, Alberto
and Rayan, Steven",
editor="Laza, Radu
and Sch{\"u}tt, Matthias
and Yui, Noriko",
title="Introduction to Nonabelian Hodge Theory",
bookTitle="Calabi-Yau Varieties: Arithmetic, Geometry and Physics: Lecture Notes on Concentrated Graduate Courses",
year="2015",
publisher="Springer New York",
address="New York, NY",
pages="131--171",
isbn="978-1-4939-2830-9",
doi="10.1007/978-1-4939-2830-9_5",
url="https://doi.org/10.1007/978-1-4939-2830-9_5"
}

@article{RV18,
author       = {Roberts, David Michael and Vozzo, Raymond F.},
title        = {Smooth loop stacks of differentiable stacks and gerbes},
journaltitle = {Cahiers de Topologie et G{\'e}om{\'e}trie Diff{\'e}rentielle Cat{\'e}goriques},
volume       = {59},
number       = {2},
year         = {2018},
pages        = {95--141},
issn         = {1245-530X},
url          = {https://cahierstgdc.com/wp-content/uploads/2018/04/LIX-2_Roberts-Vozzo.pdf},
eprint       = {1602.07973},
eprinttype   = {arxiv},
eprintclass  = {math.DG}
}

@incollection{RV18HomStack,
author    = {Roberts, David Michael and Vozzo, Raymond F.},
title     = {The Smooth Hom-Stack of an Orbifold},
booktitle = {2016 MATRIX Annals},
editor    = {de Gier, Jan and Praeger, Cheryl E. and Tao, Terence},
series    = {MATRIX Book Series},
volume    = {1},
pages     = {43--47},
publisher = {Springer},
date      = {2018-04-11},
doi       = {10.1007/978-3-319-72299-3_3},
url       = {https://link.springer.com/chapter/10.1007/978-3-319-72299-3_3},
urldate   = {2026-04-07}
}

@article{BH11,
  ISSN = {00029947},
  URL = {http://www.jstor.org/stable/41307457},
  author = {John C. Baez and Alexander E. Hoffnung},
  journal = {Transactions of the American Mathematical Society},
  number = {11},
  pages = {5789--5825},
  publisher = {American Mathematical Society},
  title = {Convenient categories of smooth spaces},
  urldate = {2026-08-14},
  volume = {363},
  year = {2011}
}

@MISC{CarchediComparison,
TITLE = {Stacks over diffeologies},
AUTHOR = {David Carchedi},
HOWPUBLISHED = {MathOverflow},
NOTE = {URL:https://mathoverflow.net/q/147425 (version: 2013-11-09)},
EPRINT = {https://mathoverflow.net/q/147425},
URL = {https://mathoverflow.net/q/147425}
}

@article{BX11,
author       = {Behrend, Kai and Xu, Ping},
title        = {Differentiable stacks and gerbes},
journaltitle = {Journal of Symplectic Geometry},
year         = {2011},
volume       = {9},
number       = {3},
pages        = {285--341},
doi          = {10.4310/JSG.2011.v9.n3.a2},
eprint       = {math/0605694},
eprinttype   = {arxiv},
url          = {https://projecteuclid.org/journalArticle/Download?urlId=jsg%2F1310388899}
}

@misc{NonAbHodgeFilt,
title={The Hodge filtration on nonabelian cohomology}, 
author={Simpson, Carlos T.},
year={1996},
eprint={alg-geom/9604005},
archivePrefix={arXiv},
primaryClass={alg-geom},
url={https://arxiv.org/abs/alg-geom/9604005}, 
}

@book{AlperModuli,
title={Stacks and Moduli},
author={Jarod Alper},
year={2025},
url={https://sites.math.washington.edu/~jarod/moduli.pdf}
}

@misc{stacks-project,
  author       = {The {Stacks project authors}},
  title        = {The Stacks project},
  howpublished = {\url{https://stacks.math.columbia.edu}},
  year         = {2026},
}

@book{Diff,
title={Diffeology},
year={2013},
author={Patrick Iglesias-Zemmour},
series={Mathematical Surveys and Monographs},
volume={185},
publisher={AMS}
}
